\documentclass[leqno,12pt]{article}
\usepackage{graphicx,amsthm}
\usepackage{verbatim,amsmath,amscd,amssymb,color}
\usepackage{enumerate}
\usepackage[all]{xy}
\xyoption{all}
\numberwithin{equation}{section}

\font\germ=eufm10

\newcommand{\cB}{{\mathcal B}}

\newcommand{\cV}{{\mathcal V}}

\newcommand{\frg}{\mathfrak g}

\newcommand{\frt}{\mathfrak t}

\newcommand{\bbC}{\mathbb C}

\newcommand{\bbP}{\mathbb P}

\newcommand{\bbZ}{\mathbb Z}

\newcommand{\wt}{\mathrm{wt}}

\newcommand{\lbr}{\begin{bmatrix}}
\newcommand{\rbr}{\end{bmatrix}}
\newcommand{\for}{\bigcirc\kern-2.6ex \because}
\newcommand{\forb}{\bigcirc\kern-2.8ex \because}
\newcommand{\forbb}{\bigcirc\kern-3.0ex \because}
\newcommand{\forbbb}{\bigcirc\kern-3.1ex \because}
\newcommand{\cd}{commutative diagram }

\newcommand\C{\mathbb C}

\def\ge{\frg}

\def\al{\alpha}

\def\beneme{\begin{enumerate}}
\def\beq{\begin{equation}}
\def\beqn{\begin{eqnarray}}
\def\beqnn{\begin{eqnarray*}}

\def\bbra#1,#2,#3{\left\{\begin{array}{c}\hspace{-5pt}
#1;#2\\ \hspace{-5pt}#3\end{array}\hspace{-5pt}\right\}}
\def\cd{\cdots}

\def\del{\delta}
\def\Del{\Delta}

\def\Delre{\Delta^{\rm re}}

\def\eit{\tilde{e}_i}
\def\eneme{\end{enumerate}}

\def\eeq{\end{equation}}
\def\eeqn{\end{eqnarray}}
\def\eeqnn{\end{eqnarray*}}
\def\fit{\tilde{f}_i}

\def\gau#1,#2{\left[\begin{array}{c}\hspace{-5pt}#1\\
\hspace{-5pt}#2\end{array}\hspace{-5pt}\right]}

\def\ify{\infty}
\def\io{\iota}

\def\lan{\langle}
\def\lar{\longrightarrow}
\def\max{{\rm max}}
\def\lm{\lambda}
\def\Lm{\Lambda}
\def\mapright#1{\smash{\mathop{\longrightarrow}\limits^{#1}}}

\def\nd{\noindent}
\def\nn{\nonumber}

\def\ot{\otimes}

\def\osigma{\ovl\sigma}
\def\ovl{\overline}

\def\qq{\qquad}
\def\q{\quad}
\def\qed{\hfill\framebox[2mm]{}}

\def\ra{\rightarrow}
\def\ran{\rangle}

\def\til{\tilde}
\def\tm{\times}
\def\tt{\frt}
\def\TY(#1,#2,#3){#1^{(#2)}_{#3}}

\def\uq{U_q(\ge)}
\def\uqp{U'_q(\ge)}

\def\uqm{U^-_q(\ge)}

\def\vep{\varepsilon}
\def\vp{\varphi}

\def\W1{W(\varpi_1)}

\def\wt{{\rm wt}}
\def\wtil{\widetilde}
\def\what{\widehat}

\def\ZZ{\mathbb Z}

\def\m@th{\mathsurround=0pt}

\def\fsquare(#1,#2){v_{#2}}

\newtheorem{thm}{Theorem}[section]
\newtheorem{pro}[thm]{Proposition}
\newtheorem{lem}[thm]{Lemma}
\newtheorem{ex}[thm]{Example}
\newtheorem{cor}[thm]{Corollary}
\newtheorem{conj}[thm]{Conjecture}
\theoremstyle{definition}
\newtheorem{df}[thm]{Definition}

\newcommand{\seteq}{\mathbin{:=}}
\newcommand{\cl}{\colon}
\newcommand{\be}{\begin{enumerate}}
\newcommand{\ee}{\end{enumerate}}
\newcommand{\bnum}{\be[{\rm (i)}]}
\newcommand{\enum}{\ee}
\newcommand{\ro}{{\rm(}}
\newcommand{\rf}{{\rm)}}
\newcommand{\set}[2]{\left\{#1\,\vert\,#2\right\}}
\newcommand{\sbigoplus}{{\mbox{\small{$\bigoplus$}}}}
\newcommand{\ba}{\begin{array}}
\newcommand{\ea}{\end{array}}
\newcommand{\on}{\operatorname}
\newcommand{\eq}{\begin{eqnarray}}
\newcommand{\eneq}{\end{eqnarray}}
\newcommand{\bysame}{\makebox[3em]{\hrulefill}\thickspace}
\newcommand{\hs}{\hspace*}

\begin{document}
\font\germ=eufm10
\def\bl{\bullet}
\def\aaa{@}
\title{\Large\bf Affine Geometric Crystals 
and \\
Limit of Perfect Crystals}
\vskip0.6cm
\author{ Masaki K\textsc{ashiwara}
\thanks{Research Institute for Mathematical Sciences, 
Kyoto University, Kitashiwakawa, Sakyo-ku, Kyoto 606, Japan. 
E-mail: masaki{\aaa}kurims.kyoto-u.ac.jp
}
\and 
Toshiki N\textsc{akashima}
\thanks{Department of Mathematics, 
Sophia University, Kioicho 7-1, Chiyoda-ku,
Tokyo 102-8554, Japan.\hfill\break
\qquad E-mail: toshiki{\aaa}mm.sophia.ac.jp
}
\and Masato O\textsc{kado}
\thanks{
Department of Mathematical Science, Graduate School of Engineering
Science, Osaka University, Toyonaka, Osaka 560-8531, Japan.
E-mail: okado{\aaa}sigmath.es.osaka-u.ac.jp}}

\date{}
\maketitle


\section{Introduction}

The theory of perfect crystals is introduced in 
\cite{KMN1}, \cite{KMN2}
in order to study certain physical models, called
vertex models, which are associated with
quantum $R$-matrices.
A perfect crystal is a  (pseudo-)crystal base
of a finite-dimensional module of a
quantum affine algebra $\uq$
(for more details, see \ref{def-perfect}), and 
it is labeled by 
a positive integer $l$, called the level
(Definition~\ref{perfect-def}). 
One of the most important 
properties of perfect crystals is as follows:
for a perfect crystal  $B$ of level $l$ and 
the crystal $B(\lm)$ of the irreducible 
$\uq$-module with a dominant integral weight $\lm$ of level $l$
as a highest weight,
there exist a unique dominant integral 
weight $\mu$ of 
level $l$ and an isomorphism of crystals:
$$ B(\lm)\cong B(\mu)\ot B.$$
By iterating this isomorphism,
we obtain the path realization of $B(\lm)$, which 
plays a crucial role in solving the vertex models.
In \cite{KMN2}, for each $\ge=A_n^{(1)}$, $B_n^{(1)}$,
$C_n^{(1)}$, $D_n^{(1)}$, $A_{2n-1}^{(2)}$,
$D_{n+1}^{(2)}$, $A_{2n}^{(2)}$, 
perfect crystals are constructed explicitly 
by patching two classical crystals.
It is conjectured that 
certain ``Kirillov-Reshetikhin modules'' have
 perfect crystal bases and all perfect crystals 
are obtained as a tensor product of those perfect 
crystal bases (\cite{HKOTY}, \cite{HKOTT}).
But, it is still far from the 
complete classification of perfect crystals.

Let $\{B_l\}_{l\geq1}$ be a family of perfect crystals.
If it satisfies certain conditions
(Definition~\ref{def-limit}), there exists
a {\it limit} $B_\ify$ of $\{B_l\}_{l\geq1}$.
In such a case the family
$\{B_l\}_{l\geq1}$ is called a 
{\it coherent family} of perfect crystals
(\cite{KKM}).
Let $B(\ify)$ be the crystal of the lower 
nilpotent subalgebras $\uqm$ of the quantum 
algebra $\uq$. Then, 
similarly to perfect crystals, 
we have an isomorphism of crystals:
\[
 B(\ify)\ot B_\ify\cong B(\ify).
\]
By iterating this isomorphism, 
the path realization of $B(\ify)$ is obtained.
In \cite{KKM}, for 
$\ge=A_n^{(1)}$, $B_n^{(1)}$,
$C_n^{(1)}$, $D_n^{(1)}$, $A_{2n-1}^{(2)}$,
$D_{n+1}^{(2)}$, $A_{2n}^{(2)}$,
explicit forms of the limit of 
perfect crystals are described, which are also
reviewed in \S \ref{limit} below.  

A. Berenstein and D. Kazhdan
introduced the theory of geometric crystals (\cite{BK}),
which is a structure  on geometric objects
analogous to crystals.
%
Let $I$ be the index set of simple roots of $\ge$.
A geometric crystal 
consists of a variety $X$, 
a rational $\bbC^\times$-
action $e_i\cl\bbC^\times\times X\longrightarrow 
X$ and rational functions
$\gamma_i,\vep_i \cl X\longrightarrow\bbC$ which 
satisfy certain conditions 
(see Definition~\ref{def-gc}).
This structure resembles
the one of crystals, for instance, 
we have the tensor product of a pair of geometric crystals
similarly to
the tensor product of crystals.
There is a more direct and remarkable relation 
between crystals and geometric crystals,
that is, {\it ulrta-discretization} functor $\cal UD$.
This is a functor from the category of 
{\it positive }geometric 
crystals to that of crystals, where a positive 
geometric crystal is a geometric crystal equipped with 
a positive
structure (see \S \ref{positive-str}).
Applying this functor, 
positive rational functions 
are transferred to piecewise linear functions by 
the following simple correspondence:
\[
x\times y\mapsto x+y,\qq
\frac{x} {y}\mapsto x-y,\qq
x+y\mapsto \max(x,y).
\]

The purpose of this paper is to 
construct a
positive affine geometric crystal $\cV(\ge)$,
whose ultra-discretization is isomorphic to
a limit of perfect crystals 
$B_\ify(\ge^L)$, where $\ge^L$ is 
the Langlands dual of $\ge$.

Let $G$ (resp.~$\ge$) be the affine 
Kac-Moody group 
(resp.~algebra) associated with a 
generalized Cartan matrix $A=(a_{ij})_{i,j\in I}$. Let
 $B^\pm$ be the Borel subgroup and $T$ the 
maximal torus. Set $y_i(c)\seteq\exp(cf_i)$,
and let $\al_i^\vee(c)\in T$ be the image of $c\in\bbC^\times$
by the group morphism $\bbC^\times\to T$ induced by
the simple coroot $\alpha_i^\vee$ as in \ref{KM}.
We set $Y_i(c)\seteq y_i(c^{-1})\,\al_i^\vee(c)=\al_i^\vee(c)\,y_i(c)$.
Let $W$ (resp.~$\wtil W$) be the Weyl group 
(resp.~the extended Weyl group) associated with $\ge$.
The Schubert cell $X_w\seteq BwB/B$ $(w=s_{i_1}\cd s_{i_k}\in W)$ 
is isomorphic to the variety
\[
 B^-_w\seteq\set{Y_{i_1}(x_1)\cd Y_{i_k}(x_k)}%
{x_1,\cd,x_k\in \bbC^\times}\subset B^-,
\]
and $X_w$ has a natural geometric crystal structure
(\cite{BK}, \cite{N}).

We choose $0\in I$ as in \cite{K0},
and let $\{\varpi_i\}_{i\in I\setminus\{0\}}$
be the set of level $0$ fundamental weights.
Let $W(\varpi_i)$ be the fundamental representation
of $U_q(\ge)$ with $\varpi_i$ as an 
extremal weight (\cite{K0}).
Let us denote its reduction at $q=1$
by the same notation $W(\varpi_i)$. 
It is a finite-dimensional $\ge$-module.
The module $W(\varpi_i)$ is 
irreducible over $\uq$, but
its reduction at $q=1$ 
is not necessarily an irreducible $\ge$-module.
We set $\bbP(\varpi_i)\seteq
(W(\varpi_i)\setminus\{0\})/\bbC^\times$.

For any $i\in I$, define (see \ref{aff-wt} for the 
definition of the inner product).
\begin{eqnarray}
&& c_i^\vee\seteq
\mathrm{max }(1,\frac{2}{(\al_i,\al_i)}).
\label{eq:ci}
\end{eqnarray}
Then the translation $t(c^\vee_i\varpi_i)$ 
belongs to $\widetilde W$ (see \ref{shift}).
For a subset $J$ of $I$, let us denote by
$\ge_J$ the subalgebra of $\ge$ generated by
$\{e_i,f_i\}_{i\in J}$.
For an integral weight $\mu$, define 
$I(\mu)\seteq\set{j\in I}{\lan h_j,\mu\ran\geq0}$. 

Then we have the following conjecture (\cite{HKOTY,HKOTT,K0}).

\begin{conj}
\bnum
\item
For any $i\in I\setminus\{0\}$ and any positive integer
$n$, there exists
a finite-dimensional irreducible $\uq$-module $W(n\varpi_i)$
with extremal weight $n\varpi_i$, called ``Kirillov-Reshetikhin module''.
\item
The module $W(n\varpi_i)$ has a global crystal basis.
\item
Its crystal $B(n\varpi_i)$ is perfect if and only if $n\in c_i^\vee\bbZ$,
and its level is equal to $n/c_i^\vee$.
\item
$\{B(nc_i^\vee\varpi_i)\}_n$ is a coherent family
of perfect crystals. We denote by $B_\infty(\varpi_i)$
its limit crystal.
\ee
\end{conj}

Now, let us state our conjecture.

\begin{conj}
For any $i\in I\setminus\{0\}$, 
there exist a unique variety $X$ endowed with
a positive $\ge$-geometric crystal structure and a
rational mapping $\pi\cl X\longrightarrow
\bbP(\varpi_i)$ satisfying the following property:
\bnum
\item
for an arbitrary extremal vector $u\in W(\varpi_i)_\mu$,
writing the translation
$t(c_i^\vee\mu)$ as $\io w\in 
\wtil W$ with a Dynkin diagram automorphism $\io$
and $w=s_{i_1}\cd s_{i_k}$ \ro see {\rm \ref{shift}}\rf, 
there exists a birational mapping
$\xi\cl B^-_w\longrightarrow X$
such that $\xi$ is a morphism of $\ge_{I(\mu)}$-geometric crystals
and that
the composition
$\pi\circ\xi\cl B^-_w\to \bbP(\varpi_i)$
coincides with
$Y_{i_1}(x_1)\cdots Y_{i_k}(x_k)\mapsto 
Y_{i_1}(x_1)\cd Y_{i_k}(x_k)\ovl u$,
where $s_{i_1}\cd s_{i_k}$ is a reduced expression of
$w$ and $\ovl u$ is the line including $u$,
\item
the ultra-discretization of $X$
is isomorphic to the crystal $B_\infty(\varpi_i)$
of the Langlands dual $\ge^L$.
\ee
\end{conj}
In this paper, we construct
a positive geometric crystal for some $i=1$
and for $\ge=A_n^{(1)}$, $B_n^{(1)}$,
$C_n^{(1)}$, $D_n^{(1)}$, $A_{2n-1}^{(2)}$,
$D_{n+1}^{(2)}$, $A_{2n}^{(2)}$,
with this conjecture as a guide.


Let 
$W(\varpi_1)$ be
the fundamental $\ge$-module as above (\cite{K0}) and 
$u_1\seteq u_{\varpi_1}$ an extremal weight
vector of $W(\varpi_1)$ with the weight $\varpi_1$.

Then $t(c_1^\vee\varpi_1)$ belongs to $\wtil W$ and
(see \ref{shift}) there exist
$w_1\in W$ and a Dynkin diagram automorphism $\io$
such that $t(c_1^\vee\varpi_1)=\io\cdot w_1$. 
Associated with this $w_1$,
we define 
\[
X_{1}\seteq B^-_{w_1}.
\]
Then $I(\varpi_1)=I\setminus\{0\}$
and $X_{1}$ has a structure of
$\ge_{I\setminus\{0\}}$-geometric crystal.

Let us choose another extremal weight vector
$u_2$ with extremal weight $\eta$ such that
$I(\eta)=I\setminus\{i_2\}$ for some $i_2\not=0$.
Then  there exists
$w_2\in W$ such that $t(c_1^\vee\eta)=\io\cdot w_2$. 
We define similarly
\[
X_{2}\seteq B^-_{w_2} 
\]
which has a structure of
$\ge_{I\setminus\{i_2\}}$-geometric crystal.
We have the rational mapping $X_\nu\to \bbP(\varpi_1)$ ($\nu=1,2$).

Then we see (by a case-by-case calculation) that
there exists a unique positive birational mapping
$\xi\cl X_1\to X_2$ such that the diagram
$$\xymatrix{X_1\ar[r]^\xi\ar[dr]&X_2\ar[d]\\&\bbP(\varpi_1)}$$
commutes and
$\xi$ commutes with $\varepsilon_i$ for $i\not=0,i_2$.
Moreover, $\xi$ is an isomorphism of
$\ge_{I\setminus\{0,i_2\}}$-geometric crystals.

Now the $\ge$-geometric crystal $\cV(\ge)$
is obtained by patching $X_1$ and $X_2$ by $\xi$.
The relations (i)-(iv) in Definition~\ref{def-gc} is obvious
for $(i,j)\not=(0,i_2), (i_2,0)$ and
we check these relations for $(i,j)=(0,i_2), (i_2,0)$ by hand.
In this paper we choose $\eta=\sigma\varpi_1$ 
for some Dynkin diagram automorphism $\sigma$ 
except the case $\ge=A_{2n}^{(2)}$.

We then show that the ultra-discretization of
$\cV(\ge)$ is isomorphic to the crystal $B_\infty(\varpi_1)$
of the Langlands dual $\ge^L$.

\medskip
The organization of the paper is as follows.
In \S \ref{sec2}, we review basic facts about
geometric crystals following
\cite{BK}, \cite{N}, \cite{N2}. 
In \S \ref{limit}, we shall recall
the notion of the limit of 
perfect crystals and give examples of 
the limit $B_\infty(\ge^L)$ of perfect crystals following \cite{KKM}.
In \S \ref{sec4}, we present the explicit form of the fundamental 
representation $W(\varpi_1)$ for each affine Lie algebra $\ge$.
In \S \ref{sec5}, we shall construct the affine geometric crystal
$\cV(\ge)$
associated with the fundamental representation $W(\varpi_1)$.
In \S \ref{sec6}, we prove that the ultra-discretization
${\cal UD}(\cV(\ge))$ is isomorphic to $B_\ify(\ge^L)$.

\section{Geometric crystals}\label{sec2}

In this section, 
we review Kac-Moody groups and geometric crystals
following 
\cite{BK}, \cite{Kac}, \cite{KP}, \cite{Ku2}, \cite{N}, \cite{PK}.
\subsection{Kac-Moody algebras and Kac-Moody groups}
\label{KM}
Fix a symmetrizable generalized Cartan matrix
 $A=(a_{ij})_{i,j\in I}$ with a finite index set $I$.
Let $(\tt,\{\al_i\}_{i\in I},\{\al^\vee_i\}_{i\in I})$ 
be the associated
root data, where ${\tt}$ is a vector space 
over $\bbC$ and
$\{\al_i\}_{i\in I}\subset\tt^*$ and 
$\{\al^\vee_i\}_{i\in I}\subset\tt$
are linearly independent 
satisfying $\al_j(\al^\vee_i)=a_{ij}$.

The Kac-Moody Lie algebra $\ge=\ge(A)$ associated with $A$
is the Lie algebra over $\bbC$ generated by $\tt$, the 
Chevalley generators $e_i$ and $f_i$ $(i\in I)$
with the usual defining relations (\cite{KP},\cite{PK}).
There is the root space decomposition 
$\ge=\sbigoplus_{\al\in \tt^*}\ge_{\al}$.
Denote the set of roots by 
$\Delta\seteq\set{\al\in \tt^*}{\al\ne0,\,\,\ge_{\al}\ne(0)}$.
Set $Q=\sum_i\bbZ \al_i$, $Q_+=\sum_i\bbZ_{\geq0} \al_i$,
$Q^\vee=\sum_i\bbZ \al^\vee_i$
and $\Delta_+=\Delta\cap Q_+$.
An element of $\Delta_+$ is called 
a {\it positive root}.
Let $P\subset \tt^*$ be a weight lattice such that 
$\bbC\ot P=\tt^*$, whose element is called a
weight.

Define the simple reflections $s_i\in{\rm Aut}(\tt)$ $(i\in I)$ by
$s_i(h)\seteq h-\al_i(h)\al^\vee_i$, which generate the Weyl group $W$.
It induces the action of $W$ on $\tt^*$ by
$s_i(\lm)\seteq\lm-\lm(\al^\vee_i)\al_i$.
Set $\Delre\seteq\set{w(\al_i)}{w\in W,\,\,i\in I}$, whose element 
is called a real root.

Let $G$ be the Kac-Moody group associated 
with $\ge$ (\cite{PK}).
Let $U_{\al}\seteq\exp\ge_{\al}$ $(\al\in \Delre)$
be the one-parameter subgroup of $G$.
The group $G$ is generated by $U_{\al}$ $(\al\in \Delre)$.
Let $U^{\pm}$ be the subgroup generated by $U_{\pm\al}$
($\al\in \Delre_+=\Delre\cap Q_+$), {\it i.e.,}
$U^{\pm}\seteq\lan U_{\pm\al}|\al\in\Del^{\rm re}_+\ran$.

For any $i\in I$, there exists a unique group homomorphism
$\phi_i\cl SL_2(\bbC)\rightarrow G$ such that
\[
\phi_i\left(
\left(
\begin{array}{cc}
1&t\\
0&1
\end{array}
\right)\right)=\exp(t e_i),\,\,
 \phi_i\left(
\left(
\begin{array}{cc}
1&0\\
t&1
\end{array}
\right)\right)=\exp(t f_i)\qquad(t\in\bbC).
\]
Set $\al^\vee_i(c)\seteq
\phi_i\left(\left(
\begin{smallmatrix}
c&0\\
0&c^{-1}\end{smallmatrix}\right)\right)$,
$x_i(t)\seteq\exp{(t e_i)}$, $y_i(t)\seteq\exp{(t f_i)}$, 
$G_i\seteq\phi_i(SL_2(\bbC))$,
$T_i\seteq \alpha_i^\vee(\bbC^\times)$ 
and 
$N_i\seteq N_{G_i}(T_i)$. Let
$T$ be the subgroup of $G$ 
with the Lie algebra $\tt$
which is called a {\it maximal torus} in $G$, and let
$B^{\pm}=U^{\pm}T$ be the Borel subgroup of $G$.
Let $N$ be the subgroup of $G$ generated by
the $N_i$'s. Then we have the isomorphism
$\phi\cl W\mapright{\sim}N/T$ defined by $\phi(s_i)=N_iT/T$.
An element $\ovl s_i\seteq x_i(-1)y_i(1)x_i(-1)
=\phi_i\left(
\left(\begin{smallmatrix}
0&-1\\
1&0
\end{smallmatrix}
\right)\right)$ in 
$N_G(T)$ is a representative of 
$s_i\in W=N_G(T)/T$. 

\subsection{Geometric crystals}\label{sec:gc}
Let $W$ be the  Weyl group associated with $\ge$. 
Define $R(w)$ for $w\in W$ by
\[
 R(w)\seteq \set{(i_1,i_2,\cd,i_l)\in I^l}{w=s_{i_1}s_{i_2}\cd s_{i_l}},
\]
where $l$ is the length of $w$, {\em i.e.},
$R(w)$ is the set of reduced expressions of $w$.

Let $X$ be a variety, 
{$\gamma_i\cl X\rightarrow \bbC$} and 
$\vep_i\cl X\longrightarrow \bbC$ ($i\in I$) 
rational functions on $X$, and
{$e_i\cl\bbC^\times \times X\longrightarrow X$}
$((c,x)\mapsto e^c_i(x))$ a
rational $\bbC^\times$-action.
For $w\in W$ and ${\bf i}=(i_1,\cd,i_l)\in R(w)$, set 
$\al^{(j)}\seteq s_{i_l}\cd s_{i_{j+1}}(\al_{i_j})$ 
$(1\leq j\leq l)$ and 
\begin{eqnarray*}
e_{\bf i}\cl T\times X&\rightarrow &X\\
(t,x)&\mapsto &e_{\bf i}^t(x)\seteq e_{i_1}^{\al^{(1)}(t)}
e_{i_2}^{\al^{(2)}(t)}\cd e_{i_l}^{\al^{(l)}(t)}(x).
\label{tx}
\end{eqnarray*}
\begin{df}
\label{def-gc}
A quadruple $(X,\{e_i\}_{i\in I},\{\gamma_i,\}_{i\in I},
\{\vep_i\}_{i\in I})$ is a 
$G$ (or $\ge$)-{\it geometric crystal} 
if
\bnum
\item
$\{1\}\times X\subset dom(e_i)$ 
for any $i\in I$.
\item
$\gamma_j(e^c_i(x))=c^{a_{ij}}\gamma_j(x)$.
\item
{$e_{\bf i}=e_{\bf i'}$}
for any 
$w\in W$ and ${\bf i}$.
${\bf i'}\in R(w)$.
\item
$\vep_i(e_i^c(x))=c^{-1}\vep_i(x)$.
\ee
\end{df}
Note that the condition (iii) is 
equivalent to the following so-called 
{\it Verma relations}:
\[
 \begin{array}{lll}
&\hspace{-20pt}e^{c_1}_{i}e^{c_2}_{j}
=e^{c_2}_{j}e^{c_1}_{i}&
{\rm if }\,\,a_{ij}=a_{ji}=0,\\
&\hspace{-20pt} e^{c_1}_{i}e^{c_1c_2}_{j}e^{c_2}_{i}
=e^{c_2}_{j}e^{c_1c_2}_{i}e^{c_1}_{j}&
{\rm if }\,\,a_{ij}=a_{ji}=-1,\\
&\hspace{-20pt}
e^{c_1}_{i}e^{c^2_1c_2}_{j}e^{c_1c_2}_{i}e^{c_2}_{j}
=e^{c_2}_{j}e^{c_1c_2}_{i}e^{c^2_1c_2}_{j}e^{c_1}_{i}&
{\rm if }\,\,a_{ij}=-2,\,
a_{ji}=-1,\\
&\hspace{-20pt}
e^{c_1}_{i}e^{c^3_1c_2}_{j}e^{c^2_1c_2}_{i}
e^{c^3_1c^2_2}_{j}e^{c_1c_2}_{i}e^{c_2}_{j}
=e^{c_2}_{j}e^{c_1c_2}_{i}e^{c^3_1c^2_2}_{j}e^{c^2_1c_2}_{i}
e^{c^3_1c_2}_je^{c_1}_i&
{\rm if }\,\,a_{ij}=-3,\,
a_{ji}=-1,
\end{array}
\]
Note that the last formula is different from the one in 
\cite{BK}, \cite{N}, \cite{N2} which seems to be
incorrect. 

\subsection{Geometric crystal on Schubert cell}
\label{schubert}

Let $X\seteq G/B$ be the flag 
variety, which is the inductive limit 
of finite-dimensional projective varieties.
For $w\in W$, let $X_w\seteq BwB/B\subset X$ be the
Schubert cell associated with $w$, which has 
a natural geometric crystal structure
(\cite{BK}, \cite{N}).
For ${\bf i}=(i_1,\cd,i_k)\in R(w)$, set 
\begin{equation}
B_{\bf i}^-
\seteq \{Y_{\bf i}(c_1,\cd,c_k)
\seteq Y_{i_1}(c_1)\cd Y_{i_l}(c_k)
\,\vert\, c_1\cd,c_k\in\bbC^\times\}\subset B^-
\label{bw1}
\end{equation}
where $Y_i(c)\seteq y_i(c^{-1})\al_i^\vee(c)$.
Then $B_{\bf i}^-$ 
is birationally isomorphic to $X_w$
and endowed with the induced geometric crystal structure.
The explicit form of the action $e^c_i$ on 
$B_{\bf i}^-$ is given by
\[
e_i^c(Y_{i_1}(c_1)\cd Y_{i_l}(c_k))
=Y_{i_1}({\cal C}_1)\cd Y_{i_l}({\cal C}_k)
\]
where
\begin{equation}
{\cal C}_j\seteq 
c_j\cdot \frac{\displaystyle \sum_{1\leq m\leq j,\,i_m=i}
 \frac{c}
{c_1^{a_{i_1,i}}\cd c_{m-1}^{a_{i_{m-1},i}}c_m}
+\sum_{j< m\leq k,\,i_m=i} \frac{1}
{c_1^{a_{i_1,i}}\cd c_{m-1}^{a_{i_{m-1},i}}c_m}}
{\displaystyle\sum_{1\leq m<j,\,i_m=i} 
 \frac{c}
{c_1^{a_{i_1,i}}\cd c_{m-1}^{a_{i_{m-1},i}}c_m}+
\mathop\sum_{j\leq m\leq k,\,i_m=i}  \frac{1}
{c_1^{a_{i_1,i}}\cd c_{m-1}^{a_{i_{m-1},i}}c_m}}.
\label{eici}
\end{equation}
We also have the explicit forms of 
rational functions $\vep_i$ and $\gamma_i$:
\[
 \vep_i(Y_{i_1}(c_1)\cd Y_{i_l}(c_k))=
\hs{-3ex}\sum_{1\leq m\leq k,\,i_m=i} \frac{1}
{c_1^{a_{i_1,i}}\cd c_{m-1}^{a_{i_{m-1},i}}c_m},\q
\gamma_i(Y_{i_1}(c_1)\cd Y_{i_l}(c_k))
=c_1^{a_{i_1,i}}\cd c_k^{a_{i_k,i}}.
\]
\subsection{Positive structure,
Ultra-discretization and Tropicalization}
\label{positive-str}

Let us recall the notions of 
positive structure and ultra-discretization/tropicalization.

The setting below is same as in \cite{N2}.
Set $R\seteq \bbC(c)$ and define
$$
\begin{array}{cccc}
v\cl&R\setminus\{0\}&\longrightarrow &\ZZ\\
&f(c)&\mapsto
&{\rm deg}(f(c)).
\end{array}
$$
Here ${\rm deg}$ is the degree of poles at $c=\infty$.
Note that for $f_1,f_2\in R\setminus\{0\}$, we have
\begin{equation}
v(f_1 f_2)=v(f_1)+v(f_2),\q
v\left(\frac{f_1}{f_2}\right)=v(f_1)-v(f_2).
\label{ff=f+f}
\end{equation}
We say that a non-zero rational function $f(c)\in \bbC(c)$ is 
{\it positive} if $f$ can be expressed as a ratio
of polynomials with positive coefficients.
\nd
If $f_1,\,\,f_2\in R$ 
are positive, then we have 
\begin{equation}
v(f_1+f_2)={\rm max}(v(f_1),v(f_2)).
\label{+max}
\end{equation}

Let $T\simeq(\bbC^\times)^l$ be an algebraic torus over $\bbC$ and 
$X^*(T)\seteq {\rm Hom}(T,\bbC^\vee)$ 
(resp.~$X_*(T)\seteq {\rm Hom}(\bbC^\vee,T)$) 
be the lattice of characters
(resp.~co-characters)
of $T$. 
A non-zero rational function on
an algebraic torus $T$ is called {\em positive} if
it is written as $g/h$ where
$g$ and $h$ are a positive linear combination of
characters of $T$.
\begin{df}
Let 
$f\cl T\rightarrow T'$ be 
a rational mapping between
two algebraic tori $T$ and 
$T'$.
We say that $f$ is {\em positive},
if $\chi\circ f$ is positive
for any character $\chi\cl T'\to \C$.
\end{df}
Denote by ${\rm Mor}^+(T,T')$ the set of 
positive rational mappings from $T$ to $T'$.

\begin{lem}[\cite{BK}]
\label{TTT}
For any $f\in {\rm Mor}^+(T_1,T_2)$             
and $g\in {\rm Mor}^+(T_2,T_3)$, 
the composition $g\circ f$
is well-defined and belongs to ${\rm Mor}^+(T_1,T_3)$.
\end{lem}

By Lemma \ref{TTT}, we can define a category ${\cal T}_+$
whose objects are algebraic tori over $\bbC$ and arrows
are positive rational mappings.

Let $f\cl T\rightarrow T'$ be a 
positive rational mapping
of algebraic tori $T$ and 
$T'$.
We define a map $\what f\cl X_*(T)\rightarrow X_*(T')$ by 
\[
\langle\chi,\what f(\xi)\rangle
=v(\chi\circ f\circ \xi),
\]
where $\chi\in X^*(T')$ and $\xi\in X_*(T)$.


\begin{lem}[\cite{BK}]
For any algebraic tori $T_1$, $T_2$, $T_3$, 
and positive rational mappings
$f\in {\rm Mor}^+(T_1,T_2)$, 
$g\in {\rm Mor}^+(T_2,T_3)$, we have
$\what{g\circ f}=\what g\circ\what f.$
\end{lem}
By this lemma, we obtain a functor 
\[
\begin{array}{cccc}
{\cal UD}\cl&{\cal T}_+&\longrightarrow &{{\mathrm{Set}}}\\[2ex]
&T&\mapsto& X_*(T)\\
&(f\cl T\rightarrow T')&\mapsto& 
(\what f\cl X_*(T)\rightarrow X_*(T')).
\end{array}
\]

Let us come back to the situation in \S \ref{sec:gc}.
Hence $G$ is a Kac-Moody group and $T$ is its Cartan subgroup.
\begin{df}[\cite{BK}]
Let $\chi=(X,\{e_i\}_{i\in I},\{\gamma_i\}_{i\in I},
\{\vep_i\}_{i\in I})$ be a 
$G$-geometric crystal, $T'$ an algebraic torus
and $\theta\cl T'\rightarrow X$ 
a birational mapping.
The mapping $\theta$ is called 
a {\em positive structure} on
$\chi$ if it satisfies
\bnum
\item for any $i\in I$ the rational functions
$\gamma_i\circ \theta\cl T'\rightarrow \bbC$ and 
$\vep_i\circ \theta\cl T'\rightarrow \bbC$ 
are positive,
\item
for any $i\in I$, the rational mapping
$e_{i,\theta}\cl\bbC^\tm \tm T'\rightarrow T'$ defined by
$e_{i,\theta}(c,t)
\seteq \theta^{-1}\circ e_i^c\circ \theta(t)$
is positive.
\end{enumerate}
\end{df}
Let $\theta\cl T'\rightarrow X$ be a positive structure on 
a geometric crystal $\chi=(X,\{e_i\}_{i\in I},
\{\gamma_i\}_{i\in I},$\\
$\{\vep_i\}_{i\in I})$.
Applying the functor ${\cal UD}$ 
to the positive rational mappings
$e_{i,\theta}\cl\bbC^\tm \tm T'\rightarrow T'$ and
$\gamma_i\circ \theta\cl T'\ra \bbC^\times$,
we obtain
\begin{eqnarray*}
\til e_i\seteq {\cal UD}(e_{i,\theta})&\cl&
\ZZ\tm X_*(T') \rightarrow X_*(T')\\
{\rm wt}_i\seteq {\cal UD}(\gamma_i\circ\theta)&\cl&
X_*(T')\rightarrow \bbZ.
\end{eqnarray*}
Hence
the quadruple $(X_*(T'),\{\til e_i\}_{i\in I},
\{{\rm wt}_i\}_{i\in I},\{\vep_i\}_{i\in I})$
is a free pre-crystal structure (see \cite[2.2]{BK}) 
and we denote it by ${\cal UD}_{\theta,T'}(\chi)$.
We have thus the following theorem:

\begin{thm}[\cite{BK}\cite{N}]
For any geometric crystal 
$\chi=(X,\{e_i\}_{i\in I},\{\gamma_i\}_{i\in I},
\{\vep_i\}_{i\in I})$ and a positive structure
$\theta\cl T'\rightarrow X$, the associated pre-crystal 
${\cal UD}_{\theta,T'}(\chi)$
is a crystal {\rm (see \cite[2.2]{BK})}.
\end{thm}

Now, let ${\cal GC}^+$ be the category whose 
object is a triplet
$(\chi,T',\theta)$ where 
$\chi=(X,\{e_i\},$\\
$\{\gamma_i\},\{\vep_i\})$ 
is a geometric crystal and $\theta\cl T'\rightarrow X$ 
is a positive structure on $\chi$, and morphism
$f\cl(\chi_1,T'_1,\theta_1)\longrightarrow 
(\chi_2,T'_2,\theta_2)$ is given by a morphism 
$\vp\cl X_1\longrightarrow X_2$  of geometric crystals
such that 
\[
f\seteq \theta_2^{-1}\circ\vp\circ\theta_1\cl T'_1\longrightarrow T'_2,
\]
is a positive rational mapping. Let ${\cal CR}$
be the category of crystals. 
Then by the theorem above, we have
\begin{cor}
$\cal UD$ defines a functor
\begin{eqnarray*}
&&\ba{ccc}
 {\cal UD}:{\cal GC}^+&\longrightarrow &{\cal CR}\\[2ex]
(\chi,T',\theta)&\mapsto& X_*(T')\\
(f\cl(\chi_1,T'_1,\theta_1)\rightarrow 
(\chi_2,T'_2,\theta_2))&\mapsto&
(\what f\cl X_*(T'_1)\rightarrow X_*(T'_2)).
\ea
\end{eqnarray*}

\end{cor}
We call the functor $\cal UD$
{\it ``ultra-discretization''} as in \cite{N},\cite{N2},
while it is called ``tropicalization'' in \cite{BK}.
For a crystal $B$, if there
exists a geometric crystal $\chi$ and a positive 
structure $\theta\cl T'\rightarrow X$ on $\chi$ such that 
${\cal UD}(\chi,T',\theta)\simeq B$ as crystals, 
we call an object $(\chi,T',\theta)$ in ${\cal GC}^+$
a {\it tropicalization} of $B$.

\section{Limit of perfect crystals}
\label{limit}
\renewcommand{\theequation}{\thesection.\arabic{equation}}
We review the limits of perfect crystals following \cite{KKM}.
(See also \cite{KMN1},\cite{KMN2}.)

\subsection{Crystals}

First we recall the notion of crystals,
which is obtained by
abstracting the combinatorial 
properties of crystal bases.
\begin{df}
A {\it crystal} $B$ is a set endowed with the following maps:
\begin{eqnarray*}
&& {\rm wt}\cl B\lar P,\\
&&\vep_i\cl B\lar\ZZ\sqcup\{-\infty\},\q
  \vp_i\cl B\lar\ZZ\sqcup\{-\infty\} \q{\hbox{for}}\q i\in I,\\
&&\eit\cl B\sqcup\{0\}\lar B\sqcup\{0\},
\q\fit\cl B\sqcup\{0\}\lar B\sqcup\{0\}\q{\hbox{for}}\q i\in I,\\
&&\eit(0)=\fit(0)=0.
\end{eqnarray*}
Those maps satisfy the following axioms: for
 all $b,b_1,b_2 \in B$, we have
\begin{eqnarray*}
&&\vp_i(b)=\vep_i(b)+\lan \al^\vee_i,{\rm wt}
(b)\ran,\\
&&\wt(\eit b)=\wt(b)+\al_i{\hbox{ if  }}\eit b\in B,\\
&&\wt(\fit b)=\wt(b)-\al_i{\hbox{ if  }}\fit b\in B,\\
&&\eit b_2=b_1 \Longleftrightarrow \fit b_1=b_2,\\
&&\vep_i(b)=-\ify
   \Longrightarrow \eit b=\fit b=0.
\end{eqnarray*}
\end{df}
The following tensor product structure 
is one of the most crucial properties of crystals.
\begin{thm}
\label{tensor}
Let $B_1$ and $B_2$ be crystals.
Set
$B_1\ot B_2\seteq 
\set{b_1\otimes b_2}{b_j\in B_j\;(j=1,2)}$. Then we have 
\bnum
\item $B_1\ot B_2$ is a crystal,
\item
for $b_1\in B_1$ and $b_2\in B_2$, we have
$$
\tilde f_i(b_1\otimes b_2)=
\left\{\begin{array}{ll}\tilde f_ib_1\otimes b_2&
{\rm if}\;\varphi_i(b_1)>\vep_i(b_2),\\
b_1\otimes\tilde f_ib_2&{\rm if}\;
\varphi_i(b_1)\leq\vep_i(b_2).
\end{array}\right.
$$
$$
\tilde e_i(b_1\otimes b_2)=\left\{\begin{array}{ll}
b_1\otimes \tilde e_ib_2&
{\rm if}\;\varphi_i(b_1)<\vep_i(b_2),\\
\tilde e_ib_1\otimes b_2
&{\rm if}\;\varphi_i(b_1)\geq\vep_i(b_2).
\end{array}\right.
$$
\end{enumerate}
\end{thm}

\begin{df}
Let $B_1$ and $B_2$ be crystals. A {\it morphism} of crystals
$\psi\cl B_1\lar B_2$ is a map
$\psi\cl B_1\sqcup\{0\} \lar B_2\sqcup\{0\}$
satisfying
\bnum
\item
$\psi(0)=0$, $\psi(B_1)\subset B_2$,
\item
for $b,b'\in B_1$, $\fit b=b'$ implies
$\fit\psi(b)=\psi(b')$,
\item
$
\wt(\psi(b))=\wt(b),\q \vep_i(\psi(b))=\vep_i(b),\q
  \vp_i(\psi(b))=\vp_i(b)
\quad\text{for any $b\in B_1$.}
$
\enum
\end{df}
In particular, 
a bijective morphism is called an 
{\it isomorphism of crystals}. 

\begin{ex}
\label{ex-tlm}
\bnum
\item
If $(L,B)$ is a crystal base, then $B$ is a crystal.
\item For the crystal base $(L(\ify),B(\ify))$
of the subalgebra $\uqm$ of 
the quantum algebra $\uq$, 
$B(\ify)$ is a crystal. \label{exam2}
\item
\label{tlm}
For $\lm\in P$, set $T_\lm\seteq \{t_\lm\}$. We define a crystal
structure on $T_\lm$ by 
\[
 \eit(t_\lm)=\fit(t_\lm)=0,\q\vep_i(t_\lm)=
\vp_i(t_\lm)=-\ify,\q \wt(t_\lm)=\lm.
\]
\ee
\end{ex}
\begin{df}
To a crystal $B$, a colored oriented graph
is associated by 
\[
 b_1\mapright{i}b_2\Longleftrightarrow 
\fit b_1=b_2.
\]
We call this graph the {\em crystal graph}
of $B$.
\end{df}

\subsection{Affine weights}
\label{aff-wt}

Let $\ge$ be an affine Lie algebra. 
The sets $\mathfrak t$, 
$\{\al_i\}_{i\in I}$ 
and $\{\al^\vee_i\}_{i\in I}$ be as in \ref{KM}. 
We take $\mathfrak t$ so that $\dim\mathfrak t=\sharp I+1$.
Let $\del\in Q_+$ be a unique element 
satisfying $$\set{\lm\in Q}{\lan \al^\vee_i,\lm\ran=0
\text{ for any $i\in I$}}=\bbZ\del,$$
and let ${\bf c}\in \sum_i\bbZ_{\ge0}\alpha_i^\vee\subset \ge$ 
be a unique central element
satisfying $$\set{h\in Q^\vee}{\lan h,\al_i\ran=0
\text{ for any }i\in I}=\bbZ c.$$
We write (\cite[6.1]{Kac})
\begin{eqnarray}
{\bf c}=\sum_i a_i^\vee \al^\vee_i,\qq
\del=\sum_i a_i\al_i.\label{eq:cdel}
\end{eqnarray}
Let $(\ ,\ )$ be the non-degenerate
$W$-invariant symmetric bilinear form on $\mathfrak t^*$
normalized by $(\del,\lm)=\lan {\bf c},\lm\ran$
for $\lm\in\mathfrak{t}^*$.
Let us set $\tt^*_{\rm cl}\seteq \tt^*/\bbC\del$ and let
${\rm cl}\cl\tt^*\longrightarrow \tt^*_{\rm cl}$
be the canonical projection. 
Then we have 
$\tt^*_{\rm cl}\cong \sbigoplus_i(\bbC \al^\vee_i)^*$.
Set $\tt^*_0\seteq \set{\lm\in\tt^*}{\lan {\bf c},\lm\ran=0}$,
$(\tt^*_{\rm cl})_0\seteq {\rm cl}(\tt^*_0)$. 
Then we have a positive-definite
symmetric form on $(\tt^*_{\rm cl})_0$ 
induced by the one on 
$\tt^*$. 
Let $\Lm_i\in \tt^*_{\rm cl}$ $(i\in I)$ be a 
weight such that $\lan \al^\vee_i,\Lm_j\ran=\del_{i,j}$, which 
is called a {\em fundamental weight}.
We choose 
$P$ so that $P_{\rm cl}\seteq {\rm cl}(P)$ coincides with 
$\oplus_{i\in I}\bbZ\Lm_i$ and 
we call $P_{\rm cl}$ the
{\it classical weight lattice}.

\subsection{Definitions of perfect crystal and its limit}
\label{def-perfect}

Let $\ge$ be an affine Lie algebra, $P_{\mathrm{cl}}$ 
the classical weight lattice as above and set 
$(P_{\mathrm{cl}})^+_l\seteq \set{\lm\in P_{\mathrm{cl}}}{
\lan c,\lm\ran=l,\,\,\lan \al^\vee_i,\lm\ran\geq0}$ 
for $l\in\ZZ_{>0}$.
\begin{df}
\label{perfect-def}
We say that a crystal $B$ is {\it perfect} of level $l$ if 
\bnum
\item
$B\ot B$ is connected as a crystal graph.
\item
There exists $\lm_0\in P_{\rm cl}$ such that 
\[
 \wt(B)\subset \lm_0+\sum_{i\ne0}\ZZ_{\leq0}
{\rm cl}(\al_i),\qq
\sharp B_{\lm_0}=1
\]
\item There exists a finite-dimensional 
$U'_q(\ge)$-module $V$ with a
crystal pseudo-base $B_{ps}$ 
such that $B\cong B_{ps}/{\pm1}$
\item
The maps 
$\vep,\vp\cl B^{min}\seteq \set{b\in B}{\lan c,\vep(b)\ran=l}
\mapright{}(P_{\rm cl}^+)_l$ are bijective, where 
$\vep(b)\seteq \sum_i\vep_i(b)\Lm_i$ and 
$\vp(b)\seteq \sum_i\vp_i(b)\Lm_i$.
\enum
\end{df}

Let $\{B_l\}_{l\geq1}$ be a family of 
perfect crystals of level $l$ and set 
$J\seteq \set{(l,b)}{l>0,\,b\in B^{min}_l}$.
\begin{df}
\label{def-limit}
A crystal $B_\ify$ with an element $b_\ify$ is called the
{\it limit of $\{B_l\}_{l\geq1}$}
if 
\bnum
\item
$\wt(b_\ify)=\vep(b_\ify)=\vp(b_\ify)=0$.
\item
For any $(l,b)\in J$, there exists an
embedding of 
crystals:
\begin{eqnarray*}
 f_{(l,b)}\cl&
T_{\vep(b)}\ot B_l\ot T_{-\vp(b)}\hookrightarrow
B_\ify\\
&t_{\vep(b)}\ot b\ot t_{-\vp(b)}\mapsto b_\ify
\end{eqnarray*}
\item
$B_\ify=\bigcup_{(l,b)\in J} {\rm Im}f_{(l,b)}$.
\end{enumerate}
\end{df}
\noindent
As for the crystal $T_\lm$, see Example \ref{ex-tlm} \eqref{tlm}.
If the limit of a family $\{B_l\}$ exists, 
we say that $\{B_l\}$
is a {\it coherent family} of perfect crystals.

The following is one of the most 
important properties of the limits of perfect crystals.
\begin{pro}
Let $B(\ify)$ be the crystal as in 
{\rm Example \ref{ex-tlm} \eqref{exam2}}. Then we have
the following isomorphism of crystals:
\[
B(\ify)\ot B_\ify\mapright{\sim}B(\ify).
\]
\end{pro}

In the rest of this section, for each 
affine Lie algebra 
$\ge=A_n^{(1)}$, $B_n^{(1)}$,
$C_n^{(1)}$, $D_n^{(1)}$, $A_{2n-1}^{(2)}$,
$D_{n+1}^{(2)}$, $A_{2n}^{(2)}$,
we explicitly describe an example of $B_\ify=B_\ify(\ge)$ 
 following \cite{KKM}.

\subsection{$\TY(A,1,n)$ $(n\geq2)$}
\label{a-inf} 
The Cartan matrix $(a_{ij})_{i,j\in I}$ 
$(I\seteq \{0,1,\cd,n\})$ of type $\TY(A,1,n)$ is
\[
 a_{ij}=\begin{cases}
2&\text{if $i=j$,}\\
-1&\text{if $i\equiv j\pm1\mod n+1$,}\\
0&\text{otherwise.}
\end{cases}
\]
We have ${\bf c}=\sum_{i\in I}\al^\vee_i$ and 
$\del=\sum_{i\in I}\al_i$.
A limit of perfect crystals is given as follows:
\[
 B_\ify(\TY(A,1,n))=
\bigl\{(b_1,\cd,b_{n+1})\in\bbZ^{n+1}\,\vert\,\sum_{i=1}^{n+1}b_i=0\bigr\},
\]
for $b=(b_1,\cd,b_{n+1})\in B_\ify(\TY(A,1,n))$, we
have
\[
 \begin{cases}
\til e_0b=(b_1-1,b_2,\cd,b_n,b_{n+1}+1),\\
\til e_ib=(b_1,\cd,b_{i}+1,b_{i+1}-1,\cd,b_{n+1}),
\,\,(i=1,\cd,n)\\
\fit b=\eit^{-1},
\end{cases}
\]
and
\[
 \begin{cases}
\wt(b)=(b_{n+1}-b_1)\Lm_0+\sum_{i=1}^{n}(b_i-b_{i+1})\Lm_i,\\
\vep_0(b)=b_1,\q
\vep_i(b)=b_{i+1}\,\,(i=1,\cd,n),\\
\vp_0(b)=b_{n+1},\q
\vp_i(b)=b_{i}\,\,(i=1,\cd,n).
\end{cases}
\]
\subsection{$\TY(B,1,n)$ $(n\geq3)$}
\label{b-b1n}
The Cartan matrix $(a_{ij})_{i,j\in I}$ 
$(I\seteq \{0,1,\cd,n\})$ of type $\TY(B,1,n)$ is
\[
 a_{ij}=\begin{cases}
2&i=j,\\
-1&|i-j|=1\text{ and }(i,j)\ne(0,1),(1,0),
(n,n-1)\text{ or }(i,j)=(0,2),(2,0)\\
-2&(i,j)=(n,n-1),\\
0&\text{otherwise.}
\end{cases}
\]
The Dynkin diagram is 
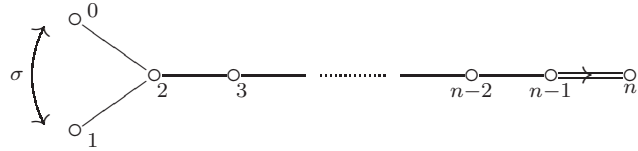
\[\SelectTips{cm}{}
\xymatrix@R=3ex{
*{\circ}<3pt> \ar@{-}[dr]^<{0} \ar@{<->}@/_/@<-2ex>[dd]_{\sigma}\\
&*{\circ}<3pt> \ar@{-}[r]_<{2} & *{\circ}<3pt> \ar@{-}[r]_<{3}
& {} \ar@{.}[r]&{} \ar@{-}[r]_>{\,\,\,\,n-2}
& *{\circ}<3pt> \ar@{-}[r]_>{\,\,\,\,n-1} &
*{\circ}<3pt> \ar@{=}[r] |-{\object@{>}}& *{\circ}<3pt>\ar@{}_<{n}\\
*{\circ}<3pt> \ar@{-}[ur]_<{1}
}
\]
where $\sigma$ is a Dynkin diagram 
automorphism which we use later. We have
\[
 {\bf c}=\al_0^\vee+\al_1^\vee+2\sum_{i=2}^{n-1}
{\al^\vee_i}+\al_n^\vee,\qq\q
\del=\al_0+\al_1+2\sum_{i=2}^n\al_i.
\]
A limit of perfect crystals is given as follows:
\[
\hspace{-5pt} B_\ify(\TY(B,1,n))
= \{(b_1,\cd,b_n,\ovl b_n,\cd,b_1)\in\bbZ^{n-1}\times
\left(\frac{1}{2}\bbZ\right)^2\times\bbZ^{n-1}
\big{\vert}\,
\sum_{i=1}^{n}(b_i+\ovl b_i)=0,\,\,
b_n+\ovl b_n\in\bbZ\},
\]
for $b=(b_1,\cd,b_n,\ovl b_n,\cd,b_1)
\in B_\ify(\TY(B,1,n))$, we have
\begin{eqnarray*}
&&\til e_0b=\begin{cases}
(b_1,b_2-1,\cd,\ovl b_2,\ovl b_1+1)&\text{if }
b_2>\ovl b_2,\\
(b_1-1,b_2,\cd,\ovl b_2+1,\ovl b_1)&\text{if }
b_2\leq\ovl b_2,
\end{cases}\\
&&\eit b=\begin{cases}
(b_1\cd,b_i+1,b_{i+1}-1,\cd,\ovl b_1)&
\text{if }b_{i+1}>\ovl b_{i+1},\\
(b_1\cd,\ovl b_{i+1}+1,\ovl b_{i}-1,\cd,\ovl b_1)
&\text{if }b_{i+1}\leq\ovl b_{i+1},
\end{cases}
\quad(i=1,\cd,n-1),\\
&&\til e_nb=(b_1,\cd,b_n+\frac{1}{2},\ovl b_{n}-\frac{1}{2},
\cd,\ovl b_1),\\
&&\fit=\eit^{-1},
\end{eqnarray*}
and
\[
 \begin{cases}
\wt(b)=(\ovl b_1-b_1+\ovl b_2-b_2)\Lm_0+
\sum_{i=1}^{n-1}(b_i-\ovl b_i+\ovl b_{i+1}-b_{i+1})\Lm_i
+2(b_n-\ovl b_n)\Lm_n,\\
\vep_0(b)=b_1+(b_2-\ovl b_2)_+,\,\,
\vep_i(b)=\ovl b_i+(b_{i+1}-\ovl b_{i+1})_+\,\,(i=1,\cd,n-1),\,
\vep_n(b)=2\ovl b_n,\\
\vp_0(b)=\ovl b_1+(\ovl b_2-b_2)_+,\,\,
\vp_i(b)=b_i+(\ovl b_{i+1}-b_{i+1})_+\,\,(i=1,\cd,n-1),\,
\vp_n(b)=2b_n.
\end{cases}
\]
Note that the presentation above is slightly 
different from the one
in \cite{KKM}. But we see that they are equivalent 
by the correspondence
$\nu_n+\frac{1}{2}\nu_0\leftrightarrow b_n$ and
$\ovl\nu_n+\frac{1}{2}\nu_0\leftrightarrow \ovl b_n$.

\subsection{$\TY(C,1,n)$ $(n\geq 2)$}
\label{c1n-ify}
The Cartan matrix $(a_{ij})_{i,j\in I}$ 
$(I\seteq \{0,1,\cd,n\})$ of type $\TY(C,1,n)$ is
\[
 a_{ij}=\begin{cases}
2&i=j,\\
-1&|i-j|=1\text{ and }(i,j)\ne(1,0),
(n-1,n),\\
-2&(i,j)=(1,0),(n-1,n),\\
0&\text{otherwise.}
\end{cases}
\]
Then the Dynkin diagram is 
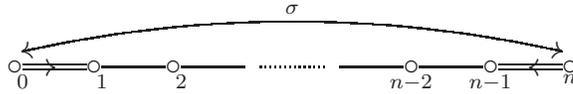
\[\SelectTips{cm}{}
\xymatrix{
*{\circ}<3pt> \ar@{=}[r] |-{\object@{>}}_<{0} 
\ar@{<->}@/^1pc/@<1ex>[rrrrrrr]^{\sigma}
&*{\circ}<3pt> \ar@{-}[r]_<{1} & *{\circ}<3pt> \ar@{-}[r]_<{2}
& {} \ar@{.}[r]&{} \ar@{-}[r]_>{\,\,\,\,n-2}
& *{\circ}<3pt> \ar@{-}[r]_>{\,\,\,\,n-1} &
*{\circ}<3pt> \ar@{=}[r] |-{\object@{<}}
& *{\circ}<3pt>\ar@{}_<{n}
}
\]
where $\sigma$ is the Dynkin 
diagram automorphism 
$\sigma\cl \al_i\mapsto\al_{n-i}$. We have 
\[
 {\bf c}=\sum_{i\in I}
{\al^\vee_i},\qq\q
\del=\al_0+2\sum_{i=1}^{n-1}\al_i+\al_n.
\]
A limit of perfect crystals is given as follows:
\[
 B_\ify(\TY(C,1,n))
= \{(b_1,\cd,b_n,\ovl b_n,\cd,\ovl b_1)\in\bbZ^{2n}|
\sum_{i=1}^{n}(b_i+\ovl b_i)\in2\bbZ\},
\]
for $b=(b_1,\cd,b_n,\ovl b_n,\cd,\ovl b_1)\in 
B_\ify(\TY(C,1,n))$, we have
\begin{eqnarray*}
&&\til e_0b=\begin{cases}
(b_1-2,b_2,\cd,\ovl b_2,\ovl b_1)&\text{if }
b_1>\ovl b_1+1,\\
(b_1-1,b_2,\cd,\ovl b_2,\ovl b_1+1)&\text{if }
b_1=\ovl b_1+1,\\
(b_1,b_2,\cd,\ovl b_2,\ovl b_1+2)&\text{if }
b_1\leq \ovl b_1,
\end{cases}\\
&&\eit b=\begin{cases}
(b_1\cd,b_i+1,b_{i+1}-1,\cd,\ovl b_1)&\text{if }
b_{i+1}>\ovl b_{i+1},\\
(b_1\cd,\ovl b_{i+1}+1,\ovl b_{i}-1,\cd,\ovl b_1)
&\text{if }b_{i+1}\leq\ovl b_{i+1},
\end{cases}
\,(1\leq i<n),\\
&&\til e_nb=(b_1,\cd,b_n+1,\ovl b_{n}-1,
\cd,\ovl b_1),\\
&&\fit=\eit^{-1},
\end{eqnarray*}
and
\[\hspace{-10pt}
 \begin{cases}
\wt(b)=(\ovl b_1-b_1)\Lm_0+
\sum_{i=1}^{n-1}(b_i-\ovl b_i+\ovl b_{i+1}-b_{i+1})\Lm_i
+(b_n-\ovl b_n)\Lm_n,\\
\vep_0(b)=-\frac{1}{2}l(b)+(b_1-\ovl b_1)_+,\,\,
\vep_i(b)=\ovl b_i+(b_{i+1}-\ovl b_{i+1})_+\,\,(1\leq i<n),\,
\vep_n(b)=\ovl b_n,\\
\vp_0(b)=-\frac{1}{2}l(b)+(\ovl b_1-b_1)_+,\,\,
\vp_i(b)=b_i+(\ovl b_{i+1}-b_{i+1})_+\,\,(1\leq i<n),\,
\vp_n(b)=b_n,
\end{cases}
\]
where $l(b)\seteq \sum_{i=1}^{n}(b_i+\ovl b_i)$.
\subsection{$\TY(D,1,n)$ $(n\geq 4)$}
\label{d1n-ify}
The Cartan matrix $(a_{ij})_{i,j\in I}$ 
$(I\seteq \{0,1,\cd,n\})$ of type $\TY(D,1,n)$ is
\[
 a_{ij}=\begin{cases}
2&i=j,\\
-1&|i-j|=1 \text{ and }1\leq i,j\leq n-1
\text{ or }(i,j)=(0,2),(2,0),
(n-2,n), (n,n-2),\\
0&\text{otherwise.}
\end{cases}
\]
The Dynkin diagram is 
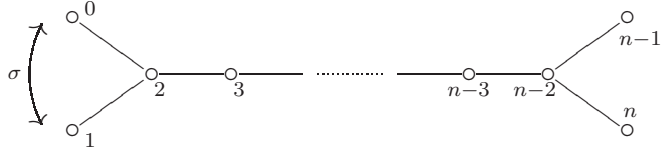
\[\SelectTips{cm}{}
\xymatrix@R=3ex{
*{\circ}<3pt> \ar@{-}[dr]^<{0} \ar@{<->}@/_/@<-2ex>[dd]_{\sigma}&&&&&&&*{\circ}<3pt> \ar@{-}[dl]^<{n-1}\\
&*{\circ}<3pt> \ar@{-}[r]_<{2} & *{\circ}<3pt> \ar@{-}[r]_<{3}
& {} \ar@{.}[r]&{} \ar@{-}[r]_>{\,\,\,\,n-3}
& *{\circ}<3pt> \ar@{-}[r]_>{n-2\,\,\,} &*{\circ}<3pt>\\
*{\circ}<3pt> \ar@{-}[ur]_<{1}&&&&&&&
*{\circ}<3pt> \ar@{-}[ul]_<{n}
}
\]
where $\sigma$ is the Dynkin diagram automorphism
$\sigma\cl\al_0\leftrightarrow \al_1$ and 
$\sigma\alpha_i=\alpha_i$ for $i\not=0,1$. We have
\[
 {\bf c}=\al_0^\vee+\al_1^\vee+2\sum_{i=2}^{n-2}
{\al^\vee_i}+\al^\vee_{n-1}+\al_n^\vee,\qq\q
\del=\al_0+\al_1+2\sum_{i=2}^{n-2}\al_i
+\al_{n-1}+\al_n.
\]
A limit of perfect crystals is given as follows:
\[
 B_\ify(\TY(D,1,n))
= \{(b_1,\cd,b_n,\ovl b_{n-1},\cd,\ovl b_1)
\in\bbZ^{2n-1}\big\vert
\sum_{i=1}^{n}b_i+\sum_{1}^{n-1}\ovl b_i=0\},
\]
for $b=(b_1,\cd,b_n,\ovl b_{n-1},\cd,\ovl b_1)
\in B_\ify(\TY(D,1,n))$, we have
\begin{eqnarray*}
&&\til e_0b=\begin{cases}
(b_1,b_2-1,\cd,\ovl b_2,\ovl b_1+1)&\text{if }
b_2>\ovl b_2,\\
(b_1-1,b_2,\cd,\ovl b_2+1,\ovl b_1)&\text{if }
b_2\leq\ovl b_2,
\end{cases}\\
&&\eit b=\begin{cases}
(b_1\cd,b_i+1,b_{i+1}-1,\cd,\ovl b_1)&
\text{if }b_{i+1}>\ovl b_{i+1},\\
(b_1\cd,\ovl b_{i+1}+1,\ovl b_{i}-1,\cd,\ovl b_1)
&\text{if }b_{i+1}\leq\ovl b_{i+1},
\end{cases}
\,(i=1,\cd,n-2),\\
&&\til e_{n-1}b=
(b_1,\cd,b_{n-1}+1,b_{n}-1,\ovl b_{n-1}, 
\cd,\ovl b_1),\\
&&\til e_nb=(b_1,\cd,b_{n-1},b_{n}+1,\ovl b_{n-1}-1
\cd,\ovl b_1),\\
&&\fit=\eit^{-1},
\end{eqnarray*}
and
\[
 \begin{cases}
\wt(b)=(\ovl b_1-b_1+\ovl b_2-b_2)\Lm_0+
\sum_{i=1}^{n-2}(b_i-\ovl b_i+\ovl b_{i+1}-b_{i+1})
\Lm_i\\
\qq\qq\qq+(b_{n-1}-\ovl b_{n-1}-b_n)\Lm_{n-1}
+(b_{n-1}-\ovl b_{n-1}+b_n)\Lm_n,\\
\vep_0(b)=b_1+(b_2-\ovl b_2)_+,\,\,
\vep_i(b)=\ovl b_i+(b_{i+1}-\ovl b_{i+1})_+\,\,
(i=1,\cd,n-2),\\
\vep_{n-1}(b)=b_n+\ovl b_{n-1},\,\,
\vep_n(b)=\ovl b_{n-1},\\
\vp_0(b)=\ovl b_1+(\ovl b_2-b_2)_+,\,\,
\vp_i(b)=b_i+(\ovl b_{i+1}-b_{i+1})_+\,\,
(i=1,\cd,n-2),\\
\vp_{n-1}(b)=b_{n-1},\,\,
\vp_n(b)=b_{n-1}+b_n.
\end{cases}
\]

\subsection{$\TY(A,2,2n-1)$ $(n\geq3)$}
\label{aon-ify}

The Cartan matrix $(a_{ij})_{i,j\in I}$ 
$(I\seteq \{0,1,\cd,n\})$ of type $\TY(A,2,2n-1)$ is
\[
 a_{ij}=\begin{cases}
2&i=j,\\
-1&|i-j|=1\text{ and }1\leq i,j\leq n-1
\text{ or }(i,j)=(0,2),(2,0),(n,n-1)\\
-2&(i,j)=(n-1,n),\\
0&\text{otherwise.}
\end{cases}
\]
The Dynkin diagram is 
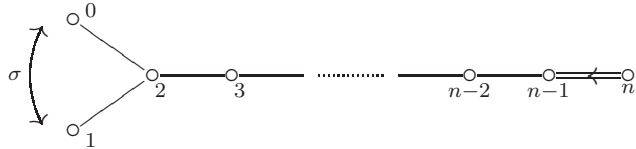
\[\SelectTips{cm}{}
\xymatrix@R=3ex{
*{\circ}<3pt> \ar@{-}[dr]^<{0} \ar@{<->}@/_/@<-2ex>[dd]_{\sigma}\\
&*{\circ}<3pt> \ar@{-}[r]_<{2} & *{\circ}<3pt> \ar@{-}[r]_<{3}
& {} \ar@{.}[r]&{} \ar@{-}[r]_>{\,\,\,\,n-2}
& *{\circ}<3pt> \ar@{-}[r]_>{\,\,\,\,n-1} &
*{\circ}<3pt> \ar@{=}[r] |-{\object@{<}}& *{\circ}<3pt>\ar@{}_<{n}\\
*{\circ}<3pt> \ar@{-}[ur]_<{1}
}
\]
where $\sigma$ is the Dynkin diagram 
automorphism $\sigma\cl\al_0\leftrightarrow \al_1$. 
We have
\[
 {\bf c}=\al_0^\vee+\al_1^\vee+2\sum_{i=2}^{n}
{\al^\vee_i},\qq\q
\del=\al_0+\al_1+2\sum_{i=2}^{n-1}\al_i
+\al_n.
\]
A limit of perfect crystals is given as follows:
\[
B_\ify(\TY(A,2,2n-1))
= \{(b_1,\cd,b_n,\ovl b_n,\cd,\ovl b_1)\in\bbZ^{2n}
\big\vert
\sum_{i=1}^{n}(b_i+\ovl b_i)=0\},
\]
for $b=(b_1,\cd,b_n,\ovl b_n,\cd,\ovl b_1)
\in B_\ify(\TY(A,2,2n-1))$, we have
\begin{eqnarray*}
&&\til e_0b=\begin{cases}
(b_1,b_2-1,\cd,\ovl b_2,\ovl b_1+1)&\text{if }
b_2>\ovl b_2,\\
(b_1-1,b_2,\cd,\ovl b_2+1,\ovl b_1)&\text{if }
b_2\leq\ovl b_2,
\end{cases}\\
&&\eit b=\begin{cases}
(b_1\cd,b_i+1,b_{i+1}-1,\cd,\ovl b_1)&
\text{if }b_{i+1}>\ovl b_{i+1},\\
(b_1\cd,\ovl b_{i+1}+1,\ovl b_{i}-1,\cd,\ovl b_1)
&\text{if }b_{i+1}\leq\ovl b_{i+1},
\end{cases}
\,(i=1,\cd,n-1),\\
&&\til e_nb=(b_1,\cd,b_{n-1},b_{n}+1,\ovl b_{n}-1,
\cd,\ovl b_1),\\
&&\fit=\eit^{-1},
\end{eqnarray*}
and
\[
 \begin{cases}
\wt(b)=(\ovl b_1-b_1+\ovl b_2-b_2)\Lm_0+
\sum_{i=1}^{n-1}(b_i-\ovl b_i+\ovl b_{i+1}-b_{i+1})
\Lm_i
+(b_{n}-\ovl b_{n})\Lm_n,\\
\vep_0(b)=b_1+(b_2-\ovl b_2)_+,\,\,
\vep_i(b)=\ovl b_i+(b_{i+1}-\ovl b_{i+1})_+\,\,
(i=1,\cd,n-1),\,\,
\vep_n(b)=\ovl b_n,\\
\vp_0(b)=\ovl b_1+(\ovl b_2-b_2)_+,\,\,
\vp_i(b)=b_i+(\ovl b_{i+1}-b_{i+1})_+\,\,
(i=1,\cd,n-1),\,
\vp_n(b)=b_n.
\end{cases}
\]

\subsection{$\TY(D,2,n+1)$ $(n\geq 2)$}
\label{d2n-ify}

The Cartan matrix $(a_{ij})_{i,j\in I}$ 
$(I\seteq \{0,1,\cd,n\})$ of type $\TY(D,2,n+1)$ is
\[
 a_{ij}=\begin{cases}
2&i=j,\\
-1&|i-j|=1\text{ and }(i,j)\ne(0,1),
(n,n-1),\\
-2&(i,j)=(0,1),(n,n-1),\\
0&\text{otherwise.}
\end{cases}
\]
The Dynkin diagram is
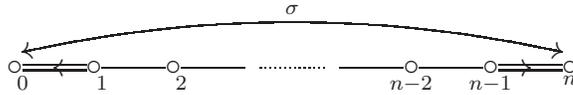
\[\SelectTips{cm}{}
\xymatrix{
*{\circ}<3pt> \ar@{=}[r] |-{\object@{<}}_<{0} 
\ar@{<->}@/^1pc/@<1ex>[rrrrrrr]^{\sigma}
&*{\circ}<3pt> \ar@{-}[r]_<{1} & *{\circ}<3pt> \ar@{-}[r]_<{2}
& {} \ar@{.}[r]&{} \ar@{-}[r]_>{\,\,\,\,n-2}
& *{\circ}<3pt> \ar@{-}[r]_>{\,\,\,\,n-1} &
*{\circ}<3pt> \ar@{=}[r] |-{\object@{>}}
& *{\circ}<3pt>\ar@{}_<{n}
}
\]
where $\sigma$ is the Dynkin diagram 
automorphism $\sigma\cl\al_i
\leftrightarrow \al_{n-i}$ $(i=0,1\cd,n)$. We have
\[
 {\bf c}=\al_0^\vee+\al_1^\vee+2\sum_{i=2}^{n-1}
{\al^\vee_i}+\al_n^\vee,\qq\q
\del=\sum_{i=0}^n\al_i.
\]
A limit of perfect crystals is given as follows:
\[
 B_\ify(\TY(D,2,n+1))
=\{(b_1,\cd,b_n,\ovl b_n,\cd,b_1)\in\bbZ^{n-1}\times
\left(\frac{1}{2}\bbZ\right)^2
\times\bbZ^{n-1}\big\vert b_n+\ovl b_n\in\ZZ\},
\]
for $b=(b_1,\cd,b_n,\ovl b_n,\cd,b_1)
\in B_\ify(\TY(D,2,n+1))$, we have
\begin{eqnarray*}
&&\til e_0b=\begin{cases}
(b_1-1,b_2,\cd,\ovl b_2,\ovl b_1)&\text{if }
b_1>\ovl b_1,\\
(b_1,b_2,\cd,\ovl b_2,\ovl b_1+1)&\text{if }
b_1\leq\ovl b_1,
\end{cases}\\
&&\eit b=\begin{cases}
(b_1\cd,b_i+1,b_{i+1}-1,\cd,\ovl b_1)&
\text{if }b_{i+1}>\ovl b_{i+1},\\
(b_1\cd,\ovl b_{i+1}+1,\ovl b_{i}-1,\cd,\ovl b_1)
&\text{if }b_{i+1}\leq\ovl b_{i+1},
\end{cases}
\,(1\leq i<n),\\
&&\til e_nb=(b_1,\cd,b_n+\frac{1}{2},\ovl b_{n}-\frac{1}{2},
\cd,\ovl b_1),\\
&&\fit=\eit^{-1},
\end{eqnarray*}
and
\[\hspace{-7pt}
 \begin{cases}
\wt(b)=2(\ovl b_1-b_1)\Lm_0+
\sum_{i=1}^{n-1}(b_i-\ovl b_i+\ovl b_{i+1}-b_{i+1})
\Lm_i
+2(b_{n}-\ovl b_{n})\Lm_n,\\
\vep_0(b)=-l(b)+2(b_1-\ovl b_1)_+,\,\,
\vep_i(b)=\ovl b_i+(b_{i+1}-\ovl b_{i+1})_+\,\,
(1\leq i<n),\,\,
\vep_n(b)=2\ovl b_n,\\
\vp_0(b)=-l(b)+2(\ovl b_1-b_1)_+,\,\,
\vp_i(b)=b_i+(\ovl b_{i+1}-b_{i+1})_+\,\,
(1\leq i<n),\,
\vp_n(b)=2b_n,
\end{cases}
\]
where $l(b)\seteq \sum_{i=1}^n(b_i+\ovl b_i)$.
Note that the presentation above is slightly 
different from the one
in \cite{KKM}. As in \ref{b-b1n}, 
they are equivalent 
by the correspondence
$\nu_n+\frac{1}{2}\nu_0\leftrightarrow b_n$ and 
$\ovl\nu_n+\frac{1}{2}\nu_0\leftrightarrow \ovl b_n$.

\subsection{$\TY(A,2,2n)$ $(n\geq 2)$}
\label{aen-ify}

The Cartan matrix $(a_{ij})_{i,j\in I}$ 
$(I\seteq \{0,1,\cd,n\})$ of type $\TY(A,2,2n)$ is
\[
 a_{ij}=\begin{cases}
2&i=j,\\
-1&|i-j|=1\text{ and }(i,j)\ne(0,1),
(n-1,n),\\
-2&(i,j)=(0,1),(n-1,n),\\
0&\text{otherwise.}
\end{cases}
\]
The Dynkin diagram is
\[\SelectTips{cm}{}
\xymatrix{
*{\circ}<3pt> \ar@{=}[r] |-{\object@{<}}_<{0} 
&*{\circ}<3pt> \ar@{-}[r]_<{1} & *{\circ}<3pt> \ar@{-}[r]_<{2}
& {} \ar@{.}[r]&{} \ar@{-}[r]_>{\,\,\,\,n-2}
& *{\circ}<3pt> \ar@{-}[r]_>{\,\,\,\,n-1} &
*{\circ}<3pt> \ar@{=}[r] |-{\object@{<}}
& *{\circ}<3pt>\ar@{}_<{n}
}
\]
Note that there exists no Dynkin diagram 
automorphism in this case. We have
\[
 {\bf c}=\al_0^\vee+2\sum_{i=1}^{n}
{\al^\vee_i},\qq\q
\del=2\sum_{i=0}^{n-1}\al_i+\al_n.
\]
A limit of perfect crystals is given as follows:
\[
B_\ify(\TY(A,2,2n))
= \{(b_1,\cd,b_n,\ovl b_n,\cd,\ovl b_1)\in\bbZ^{2n}\},
\]
for $(b_1,\cd,b_n,\ovl b_n,\cd,\ovl b_1)
\in B_\ify(\TY(A,2,2n))$, we have
\begin{eqnarray*}
&&\til e_0b=\begin{cases}
(b_1-1,b_2,\cd,\ovl b_2,\ovl b_1)&\text{if }
b_1>\ovl b_1,\\
(b_1,b_2,\cd,\ovl b_2,\ovl b_1+1)&\text{if }
b_1\leq\ovl b_1,
\end{cases}\\
&&\eit b=\begin{cases}
(b_1\cd,b_i+1,b_{i+1}-1,\cd,\ovl b_1)&
\text{if }b_{i+1}>\ovl b_{i+1},\\
(b_1\cd,\ovl b_{i+1}+1,\ovl b_{i}-1,\cd,\ovl b_1)
&\text{if }b_{i+1}\leq\ovl b_{i+1},
\end{cases}
\,(1\leq i<n),\\
&&\til e_nb=(b_1,\cd,b_{n-1},b_{n}+1,\ovl b_{n}-1,
\cd,\ovl b_1),\\
&&\fit=\eit^{-1},
\end{eqnarray*}
and
\[\hspace{-7pt}
 \begin{cases}
\wt(b)=2(\ovl b_1-b_1)\Lm_0+
\sum_{i=1}^{n-1}(b_i-\ovl b_i+\ovl b_{i+1}-b_{i+1})
\Lm_i
+(b_{n}-\ovl b_{n})\Lm_n,\\
\vep_0(b)=-l(b)+2(b_1-\ovl b_1)_+,\,\,
\vep_i(b)=\ovl b_i+(b_{i+1}-\ovl b_{i+1})_+\,\,
(1\leq i<n),\,\,
\vep_n(b)=\ovl b_n,\\
\vp_0(b)=-l(b)+2(\ovl b_1-b_1)_+,\,\,
\vp_i(b)=b_i+(\ovl b_{i+1}-b_{i+1})_+\,\,
(1\leq i<n),\,
\vp_n(b)=b_n,
\end{cases}
\]
where $l(b)\seteq \sum_{i=1}^n(b_i+\ovl b_i)$.

\section{Fundamental Representations}\label{sec4}

\subsection{Fundamental representation 
$W(\varpi_1)$}
\label{fundamental}

Let $\ge=A_n^{(1)}$, $B_n^{(1)}$,
$C_n^{(1)}$, $D_n^{(1)}$, $A_{2n-1}^{(2)}$,
$D_{n+1}^{(2)}$, $A_{2n}^{(2)}$,
and let $\set{\Lm_i}{i\in I}$ be the set of fundamental 
weight as in \S \ref{limit}.
Let $\varpi_1\seteq \Lm_1-a^\vee_1\Lm_0$ be the
(level 0) fundamental weight, where
$i=1$ is the node of the Dynkin diagram as in \S \ref{limit}
and $a^\vee_i$ is given in \eqref{eq:cdel}.

Let $V(\varpi_1)$ be the extremal weight module
over $\uq$
associated with $\varpi_1$ (\cite{K0}) and 
let $W(\varpi_1)\cong 
V(\varpi_1)/(z_1-1)V(\varpi_1)$ be the 
fundamental representation of $\uqp$
where $z_1$ is a 
$\uqp$-linear automorphism on $V(\varpi_1)$
 (see \cite[Sect 5.]{K0}). 

By \cite[Theorem 5.17]{K0}, $W(\varpi_1)$ is a
finite-dimensional irreducible integrable 
$\uqp$-module and has a global basis
with a simple crystal. Thus, we can consider 
its specialization $q=1$ and obtain a
finite-dimensional $\ge$-module 
which will be denoted by 
the same notation $W(\varpi_1)$.

Now we present the explicit form of 
$W(\varpi_1)$ for 
$\ge=A_n^{(1)}$, $B_n^{(1)}$,
$C_n^{(1)}$, $D_n^{(1)}$, $A_{2n-1}^{(2)}$,
$D_{n+1}^{(2)}$, $A_{2n}^{(2)}$.

\subsection{$\TY(A,1,n)$ $(n\geq2)$}
The global basis of $W(\varpi_1)$ is
\[
 \{\fsquare(0.5cm,1),\,\,\fsquare(0.5cm,2),\cd,
\,\,\fsquare(0.5cm,n+1)\},
\]
and we have
\[
 \wt(\fsquare(0.5cm,i))=\Lm_i-\Lm_{i-1} \q(i=1\cd,n+1),
\]
where we understand $\Lm_{n+1}=\Lm_0$.
The explicit actions of $f_i$'s are  
\begin{eqnarray*}
&&f_i\fsquare(0.5cm,i)
=\fsquare(0.5cm,i+1)\q(1\leq i\leq n),\q
f_0\fsquare(0.5cm,n+1)=\fsquare(0.5cm,1),\\
&&f_i\fsquare(0.5cm,j)=0 \qq \text{otherwise.}
\end{eqnarray*}
Its crystal graph is:
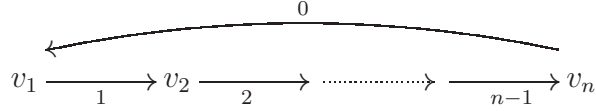
\[\SelectTips{cm}{}
\xymatrix@C=8ex{
\fsquare(5mm,1) \ar@{->}[r]_{1}&\fsquare(5mm,2)
\ar@{->}[r]_2&\ar@{.>}[r]&\ar@{->}[r]_{n-1}
&\fsquare(5mm,n)
\ar@{->}@/_1pc/@<-2ex>[llll]_0}
\]

\subsection{$\TY(B,1,n)$ $(n\geq 3)$}
\label{bn-w1}

The global basis of $\W1$ is 
\[
 \{\fsquare(0.5cm,1),\fsquare(0.5cm,2),\,\,\cd,
\,\,\fsquare(0.5cm,n),\,\,\fsquare(0.5cm,0),
\,\,\fsquare(0.5cm,\ovl n),\cd,
\,\,\fsquare(0.5cm,\ovl 2),
\,\,\fsquare(0.5cm,\ovl 1)\},
\]
and we have
\begin{eqnarray*}
&&\wt(\fsquare(5mm,i))=\Lm_i-\Lm_{i-1},\qq
\wt(\fsquare(5mm,\ovl i))=\Lm_{i-1}-\Lm_i
\qq(i\ne 0,2,n),\\
&&\wt(\fsquare(5mm,2))=-\Lm_0-\Lm_1+\Lm_2,\qq
\wt(\fsquare(5mm,\ovl 2))=\Lm_0+\Lm_1-\Lm_2,\\
&&\wt(\fsquare(5mm,n))=2\Lm_n-\Lm_{n-1},\qq
\wt(\fsquare(5mm,\ovl n))=\Lm_{n-1}-2\Lm_n,\qq
\wt(\fsquare(5mm,0))=0.
\end{eqnarray*}
The explicit forms of the actions by $f_i$'s are
\begin{eqnarray*}
&&f_i\fsquare(5mm,i)=\fsquare(5mm,i+1),\qq
f_i\fsquare(5mm,\ovl {i+1})=\fsquare(5mm,\ovl i)
\qq(i=1,\cd,n-1),\\
&&f_{n}\fsquare(5mm,n)=\fsquare(5mm,0),\qq
f_n\fsquare(5mm,0)=2\fsquare(5mm,\ovl n),\\
&&f_0\fsquare(5mm,\ovl 2)=\fsquare(5mm,1),\qq
f_0\fsquare(5mm,\ovl 1)=\fsquare(5mm,2),\\
&&f_i\fsquare(5mm,j)=0\qq\text{otherwise}.
\end{eqnarray*}
Its crystal graph is

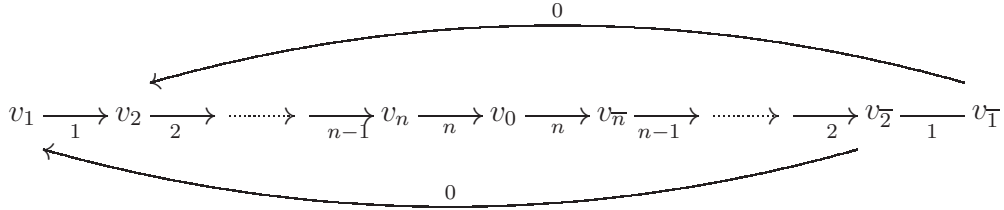
\[\SelectTips{cm}{}
\xymatrix{
\fsquare(5mm,1) \ar@{->}[r]_{1}&\fsquare(5mm,2)
\ar@{->}[r]_2&\ar@{.>}[r]&\ar@{->}[r]_{n-1}
&\fsquare(5mm,n)
\ar@{->}[r]_n&\fsquare(5mm,0) \ar@{->}[r]_n
&\fsquare(5mm,\ovl n) \ar@{->}[r]_{n-1}&
\ar@{.>}[r]&\ar@{->}[r]_2
&\fsquare(5mm,\ovl 2) \ar@{-}[r]_1
\ar@{->}@/^2pc/@<2ex>[lllllllll]_0
&\fsquare(5mm,\ovl 1)
\ar@{->}@/_2pc/@<-2ex>[lllllllll]_0}\]
\subsection{$\TY(C,1,n)$ $(n\geq 2)$}
\label{cn-w1}

The global basis of $\W1$ is 
\[
 \{\fsquare(0.5cm,1),\,\,\fsquare(0.5cm,2),\cd,
\,\,\fsquare(0.5cm,n),
\,\,\fsquare(0.5cm,\ovl n),\cd,
\,\,\fsquare(0.5cm,\ovl 2),
\,\,\fsquare(0.5cm,\ovl 1)\},
\]
and we have
\begin{eqnarray*}
&&\wt(\fsquare(5mm,i))=\Lm_i-\Lm_{i-1},\qq
\wt(\fsquare(5mm,\ovl i))=\Lm_{i-1}-\Lm_i
\qq(i=1\cd,n).
\end{eqnarray*}
The explicit forms of the actions by $f_i$'s are
\begin{eqnarray*}
&&f_i\fsquare(5mm,i)=\fsquare(5mm,i+1),\qq
f_i\fsquare(5mm,\ovl {i+1})=\fsquare(5mm,\ovl i)
\qq(i=1,\cd,n-1),\\
&&f_{n}\fsquare(5mm,n)=\fsquare(5mm,\ovl n),\qq
f_0\fsquare(5mm,\ovl 1)=\fsquare(5mm,1),\\
&&f_i\fsquare(5mm,j)=0\qq\text{otherwise}.
\end{eqnarray*}
Its crystal graph is
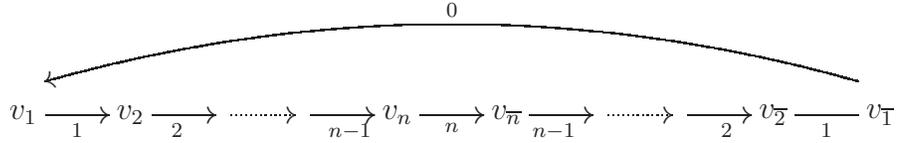
\[\SelectTips{cm}{}
\xymatrix{
\fsquare(5mm,1) \ar@{->}[r]_{1}&\fsquare(5mm,2)
\ar@{->}[r]_2&\ar@{.>}[r]&\ar@{->}[r]_{n-1}
&\fsquare(5mm,n)\ar@{->}[r]_n
&\fsquare(5mm,\ovl n) \ar@{->}[r]_{n-1}&
\ar@{.>}[r]&\ar@{->}[r]_2
&\fsquare(5mm,\ovl 2) \ar@{-}[r]_1
&\fsquare(5mm,\ovl 1) 
\ar@{->}@/_2pc/@<-2ex>[lllllllll]_0
}
\]
\subsection{$\TY(D,1,n)$ $(n\geq 4)$}
\label{dn-w1}

The global basis of $\W1$ is 
\[
 \{\fsquare(0.5cm,1),\,\,
\fsquare(0.5cm,2),\cd,
\,\,\fsquare(0.5cm,n),
\,\,\fsquare(0.5cm,\ovl n),\cd,
\,\,\fsquare(0.5cm,\ovl 2),
\,\,\fsquare(0.5cm,\ovl 1)\},
\]
and we have
\begin{eqnarray*}
&&\wt(\fsquare(5mm,i))=\Lm_i-\Lm_{i-1},\qq
\wt(\fsquare(5mm,\ovl i))=\Lm_{i-1}-\Lm_i
\qq(i\ne 2,n-1),\\
&&\wt(\fsquare(5mm,2))=-\Lm_0-\Lm_1+\Lm_2,\qq
\wt(\fsquare(5mm,\ovl 2))=\Lm_0+\Lm_1-\Lm_2,\\
&&\wt(\fsquare(5mm,n-1))=\Lm_{n-1}+
\Lm_n-\Lm_{n-1},\qq
\wt(\fsquare(6mm,\ovl{n-1}))
=\Lm_{n-2}-\Lm_{n-1}-\Lm_n.
\end{eqnarray*}
The explicit forms of the actions by $f_i$'s are
\begin{eqnarray*}
&&f_i\fsquare(5mm,i)=\fsquare(5mm,i+1),\qq
f_i\fsquare(5mm,\ovl {i+1})=\fsquare(5mm,\ovl i)
\qq(i=1,\cd,n-1),\\
&&f_{n}\fsquare(5mm,n)=\fsquare(6mm,\ovl{n-1}),
\qq
f_n\fsquare(5mm,n-1)=\fsquare(5mm,\ovl n),\\
&&f_0\fsquare(5mm,\ovl 2)=\fsquare(5mm,1),\qq
f_0\fsquare(5mm,\ovl 1)=\fsquare(5mm,2),\\
&&f_i\fsquare(5mm,j)=0\qq\text{otherwise}.
\end{eqnarray*}
Its crystal graph is 

\[
\xymatrix@R=1.5ex{
&&&&&\fsquare(5mm,n)\ar@{->}[dr]_n&&&&\\
\fsquare(5mm,1) \ar@{->}[r]_{1}&\fsquare(5mm,2)
\ar@{->}[r]_2&\ar@{.>}[r]&\ar@{->}[r]_{n-2}
&\fsquare(5mm,n-1)\ar@{->}[ur]_{n-1}
\ar@{->}[dr]_n
&&\fsquare(6mm,\ovl{n-1}) \ar@{->}[r]_{n-2}&
\ar@{.>}[r]&\ar@{->}[r]_2
&\fsquare(5mm,\ovl 2) \ar@{->}[r]_1
\ar@{->}@/^4pc/@<2ex>[lllllllll]_0
&\fsquare(5mm,\ovl 1) 
\ar@{->}@/_4pc/@<-2ex>[lllllllll]^0\\
&&&&&\fsquare(5mm,\ovl n)\ar@{->}[ur]_{n-1}&&&&
}
\]

\subsection{$\TY(A,2,2n-1)$ $(n\geq 3)$}
\label{aon-w1}

The global basis of $\W1$ is 
\[
 \{\fsquare(0.5cm,1),
\,\,\fsquare(0.5cm,2),\cd,
\,\,\fsquare(0.5cm,n),
\,\,\fsquare(0.5cm,\ovl n),\cd,
\,\,\fsquare(0.5cm,\ovl 2),
\,\,\fsquare(0.5cm,\ovl 1)\},
\]
and we have
\begin{eqnarray*}
&&\wt(\fsquare(5mm,i))=\Lm_i-\Lm_{i-1},\qq
\wt(\fsquare(5mm,\ovl i))=\Lm_{i-1}-\Lm_i
\qq(i\ne 2)\\
&&\wt(\fsquare(5mm,2))=-\Lm_0-\Lm_1+\Lm_2,\qq
\wt(\fsquare(5mm,\ovl 2))=\Lm_0+\Lm_1-\Lm_2.
\end{eqnarray*}
The explicit forms of the actions by $f_i$'s are
\begin{eqnarray*}
&&f_i\fsquare(5mm,i)=\fsquare(5mm,i+1),\qq
f_i\fsquare(5mm,\ovl {i+1})=\fsquare(5mm,\ovl i)
\qq(i=1,\cd,n-1),\\
&&f_{n}\fsquare(5mm,n)=\fsquare(5mm,\ovl{n}),\\
&&f_0\fsquare(5mm,\ovl 2)=\fsquare(5mm,1),\qq
f_0\fsquare(5mm,\ovl 1)=\fsquare(5mm,2),\\
&&f_i\fsquare(5mm,j)=0\qq\text{otherwise}.
\end{eqnarray*}
Its crystal graph is 

\[\xymatrix{
\fsquare(5mm,1) \ar@{->}[r]_{1}&\fsquare(5mm,2)
\ar@{->}[r]_2&\ar@{.>}[r]&\ar@{->}[r]_{n-1}
&\fsquare(5mm,n)\ar@{->}[r]_n
&\fsquare(5mm,\ovl n) \ar@{->}[r]_{n-1}&
\ar@{.>}[r]&\ar@{->}[r]_2
&\fsquare(5mm,\ovl 2) \ar@{->}[r]_1
\ar@{->}@/^2pc/@<2ex>[llllllll]_0
&\fsquare(5mm,\ovl 1) 
\ar@{->}@/_2pc/@<-2ex>[llllllll]_0
}\]
\subsection{$\TY(D,2,n+1)$ $(n\geq 2)$}
\label{d2n-w1}

The global basis of $\W1$ is 
\[
 \{\fsquare(0.5cm,1),
\,\,\fsquare(0.5cm,2),\cd,
\,\,\fsquare(0.5cm,n),
\,\,\fsquare(0.5cm,0),
\,\,\fsquare(0.5cm,\ovl n),\cd,
\,\,\fsquare(0.5cm,\ovl 2),
\,\,\fsquare(0.5cm,\ovl 1),
\,\,\phi\},
\]
and we have
\begin{eqnarray*}
&&\wt(\fsquare(5mm,i))=\Lm_i-\Lm_{i-1},\qq
\wt(\fsquare(5mm,\ovl i))=\Lm_{i-1}-\Lm_i
\qq(i\ne 0,1,n),\\
&&\wt(\fsquare(5mm,1))=\Lm_1-2\Lm_0,\qq
\wt(\fsquare(5mm,\ovl 1))=2\Lm_0-\Lm_1,\\
&&\wt(\fsquare(5mm,n))=2\Lm_n-\Lm_{n-1},\qq
\wt(\fsquare(5mm,\ovl n))=\Lm_{n-1}-2\Lm_n,\\
&&\wt(\fsquare(5mm,0))=0,\qq
\wt(\phi)=0
\end{eqnarray*}
The explicit forms of the actions by $f_i$'s are
\begin{eqnarray*}
&&f_i\fsquare(5mm,i)=\fsquare(5mm,i+1),\qq
f_i\fsquare(5mm,\ovl {i+1})=\fsquare(5mm,\ovl i)
\qq(i=1,\cd,n-1),\\
&&f_{n}\fsquare(5mm,n)=\fsquare(5mm,0),\qq
f_n\fsquare(5mm,0)=2\fsquare(5mm,\ovl n),\\
&&f_0\fsquare(5mm,\ovl 1)=\phi,\qq
f_0\phi=2\fsquare(5mm,1),\\
&&f_i\fsquare(5mm,j)=0,\q
f_i\phi=0\qq\text{otherwise}.
\end{eqnarray*}
Its crystal graph is
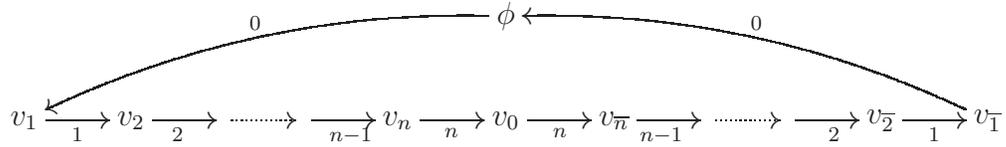
\[\SelectTips{cm}{}
\xymatrix{
&&&&&\phi
\ar@{->}@/_1pc/[dlllll]_0&&&&\\
\fsquare(5mm,1) \ar@{->}[r]_{1}&\fsquare(5mm,2)
\ar@{->}[r]_2&\ar@{.>}[r]&\ar@{->}[r]_{n-1}
&\fsquare(5mm,n)
\ar@{->}[r]_n&\fsquare(5mm,0) \ar@{->}[r]_n
&\fsquare(5mm,\ovl n) \ar@{->}[r]_{n-1}&
\ar@{.>}[r]&\ar@{->}[r]_2
&\fsquare(5mm,\ovl 2) \ar@{->}[r]_1
&\fsquare(5mm,\ovl 1)
\ar@{->}@/_1pc/[ulllll]_0}\]


\subsection{$A_{2n}^{(2)\,\dagger}$ $(n\geq 2)$}
\label{aen-v1}
We take the 
Cartan data transposed the one in \ref{aen-ify}.
Then the Dynkin diagram is
\[\SelectTips{cm}{}
\xymatrix{
*{\circ}<3pt> \ar@{=}[r] |-{\object@{>}}_<{0} 
&*{\circ}<3pt> \ar@{-}[r]_<{1} & *{\circ}<3pt> \ar@{-}[r]_<{2}
& {} \ar@{.}[r]&{} \ar@{-}[r]_>{\,\,\,\,n-2}
& *{\circ}<3pt> \ar@{-}[r]_>{\,\,\,\,n-1} &
*{\circ}<3pt> \ar@{=}[r] |-{\object@{>}}
& *{\circ}<3pt>\ar@{}_<{n}
}
\]
In this case, we denote this type by 
$A_{2n}^{(2)\,\dagger}$ in order to distinguish it 
with the one in \ref{aen-ify}.
Then the global basis of $\W1$ is 
\[
 \{\fsquare(0.5cm,1),
\,\,\fsquare(0.5cm,2),\cd,
\,\,\fsquare(0.5cm,n),
\,\,\fsquare(0.5cm,0),
\,\,\fsquare(0.5cm,\ovl n),\cd,
\,\,\fsquare(0.5cm,\ovl 2),
\,\,\fsquare(0.5cm,\ovl 1)
\},
\]
and we have
\begin{eqnarray*}
&&\wt(\fsquare(5mm,i))=\Lm_i-\Lm_{i-1},\qq
\wt(\fsquare(5mm,\ovl i))=\Lm_{i-1}-\Lm_i
\qq(i=1,\cd,n-1),\\
&&\wt(\fsquare(5mm,n))=2\Lm_n-\Lm_{n-1},\qq
\wt(\fsquare(5mm,0))=0,\qq
\wt(\fsquare(5mm,\ovl 1))=\Lm_{n-1}-2\Lm_n.
\end{eqnarray*}
The explicit forms of the actions by $f_i$'s are
\begin{eqnarray*}
&&f_i\fsquare(5mm,i)=\fsquare(5mm,i+1),\qq
f_i\fsquare(5mm,\ovl {i+1})=\fsquare(5mm,\ovl i)
\qq(i=1,\cd,n-1),\\
&&f_{n}\fsquare(5mm,n)=\fsquare(5mm,\ovl 0),
\qq f_{n}\fsquare(5mm,0)=
2\fsquare(5mm,\ovl n),\\
&&f_0\fsquare(5mm,\ovl 1)=\fsquare(5mm,1),\\
&&f_i\fsquare(5mm,j)=0,\q
\text{otherwise}.
\end{eqnarray*}
Its crystal graph is
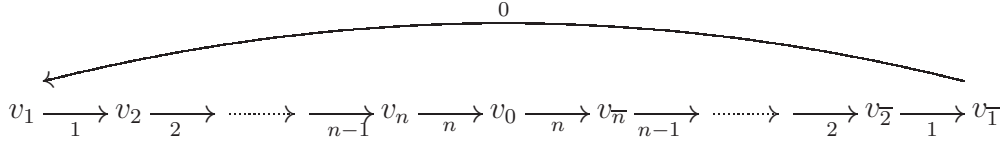
\[\SelectTips{cm}{}
\xymatrix{
\fsquare(5mm,1) \ar@{->}[r]_{1}&\fsquare(5mm,2)
\ar@{->}[r]_2&\ar@{.>}[r]&\ar@{->}[r]_{n-1}
&\fsquare(5mm,n)
\ar@{->}[r]_n&\fsquare(5mm,0) \ar@{->}[r]_n
&\fsquare(5mm,\ovl n) \ar@{->}[r]_{n-1}&
\ar@{.>}[r]&\ar@{->}[r]_2
&\fsquare(5mm,\ovl 2) \ar@{->}[r]_1
&\fsquare(5mm,\ovl 1)
\ar@{->}@/_2pc/@<-2ex>[llllllllll]_0}\]


\section{Affine Geometric Crystals}\label{sec5}

In this section, we shall construct the  
affine geometric crystal $\cV(\ge)$, 
which is realized in 
the fundamental representation $W(\varpi_1)$. 

\subsection{Translation $t(\til\varpi_1)$}
\label{shift}
For $\xi_0\in (\frt^*_{\rm cl})_0$, let $t(\xi_0)$ be 
as in \cite[Sect.4]{K0}:
\[
 t(\xi_0)(\lm)\seteq \lm+(\del,\lm)
\xi-
(\xi,\lm)\del
-\frac{(\xi,\xi)}{2}
(\del,\lm)\del
\]
for $\xi\in \mathfrak{t}^*$ such that $\mathrm{cl}(\xi)=\xi_0$.
Then $t(\xi_0)$ does not depend on the choice of $\xi$,
and it is well-defined.

Let  $c_i^\vee$ be as in \eqref{eq:ci}.
Then $t(m\varpi_i)$ belongs to $\widetilde W$ if and only if 
$m\in c_i^\vee\bbZ$.
Setting $\wtil\varpi_i\seteq c_i^\vee\varpi_i$ $(i\in I)$
(\cite{K1}),
$t(\wtil\varpi_1)$ 
is expressed as follows (see {\it e.g.} \cite{KMOTU}):
\[
 t(\wtil\varpi_1)=
\begin{cases}
\io(s_ns_{n-1}\cd s_2s_1)&\TY(A,1,n)\text{ case},\\
\io(s_1\cd s_n)(s_{n-1}\cd s_2s_1)
&\TY(B,1,n),\,\TY(A,2,2n-1)\text{ cases},\\
(s_0s_1\cd s_n)(s_{n-1}\cd s_2s_1)
\quad&\TY(C,1,n),\,\TY(D,2,n+1)\text{ cases},\\
\io(s_1\cd s_n)(s_{n-2}\cd s_2s_1)
&\TY(D,1,n)\text{ case},\\
(s_0s_1\cd s_n)(s_{n-1}\cd s_2s_1)
&A_{2n}^{(2)\,\dagger}\text{ cases},
\end{cases}
\]
where $\io$ is the Dynkin diagram automorphism
\[
\io=
 \begin{cases}\sigma&\ge=\TY(A,1,n), \TY(B,1,n),
\TY(A,2,2n-1),\\
\al_0\leftrightarrow \al_1\text{ and }
\al_{n-1}\leftrightarrow \al_{n}&\ge=\TY(D,1,n).
\end{cases}
\]
Now, we know that each $t(\til\varpi_1)$ is 
in the form $w_1$ or $\io\cdot w_1$ for $w_1\in W$,
{\it e.g.,} $w_1=s_n\cd s_1$ for $\TY(A,1,n)$,
$w_1=(s_1\cd s_n)(s_{n-1}\cd s_1)$ for 
$\TY(B,1,n)$, {\it etc.}, \ldots.

In the case $\ge=A_{2n}^{(2)\,\dagger}$,
$\eta\seteq 
\on{wt}(\fsquare(5mm,\ovl n))
=-\Lm_{n-1}+\Lm_n$ 
is a unique weight of $W(\varpi_1)$
which satisfies $\lan \al^\vee_i,\eta\ran\geq0$
for $i\ne n$.
For this $\eta$ we have
\begin{equation}
t(\eta)=(s_ns_{n-1}\cd s_1)(s_0s_1\cd s_{n-1})
=:w_2,
\label{t-pi-n}
\end{equation}
which will be  used later.
\subsection{Affine geometric crystals}

Let $\sigma$ be the Dynkin diagram automorphism
as in \S \ref{a-inf}--\ref{d2n-ify} and 
$w_1=s_{i_1}\cd s_{i_k}$ be as 
in the previous subsection.
Set
\begin{equation}
{\cV}(\ge)\seteq 
=\{v(x_1,\cd,x_k)\seteq 
Y_{i_1}(x_1)\cd Y_{i_k}(x_k)\fsquare(5mm,1)\,
 \big\vert \,
x_1\cd,x_k\in \bbC^\times\}\subset W(\varpi_1)
\label{Vxrn}
\end{equation}
Since the vector $\fsquare(5mm,1)$ 
is a highest weight
vector in $W(\varpi_1)$ as a $\ge_0$-module,
$\cV(\ge)$ has a 
$G_0$-geometric crystal structure.
Moreover $\bbC^\times\to {\cV}(\ge)$ is a birational morphism.
We shall define a $G$-geometric crystal structure 
on $\cV(\ge)$ by using the Dynkin diagram 
automorphism $\sigma$ except for $\TY(A,2,2n)$. 
This $\sigma$ induces an automorphism 
of $W(\varpi_1)$, which is also denoted by 
$\sigma\cl W(\varpi_1)\longrightarrow W(\varpi_1)$.
In the subsequent subsections, we shall show the 
following theorems by case-by-case arguments.

\begin{thm}
\label{birat}
\bnum
\item
Case $\ge\ne A_{2n}^{(2)\,\dagger}$.
For $x=(x_1\cd,x_k)\in (\bbC^\times)^k$, there 
exist a unique $y=(y_1,\cd,y_k)\in (\bbC^\times)^k$
and a positive rational function $a(x)$ such that
\begin{equation}
v(y)=a(x)\sigma(v(x)),\q
\vep_{\sigma(i)}(v(y))=\vep_i(v(x))\quad
\text{if $i,\sigma(i)\not=0$.}
\end{equation}
\item
Case $\ge=A_{2n}^{(2)\,\dagger}$. Associated with 
$w_1$ and $w_2$ as in the previous section, we define
\begin{eqnarray*}
 \cV(\ge)&\seteq &\{v_1(x)
=Y_0(x_0)Y_1(x_1)\cd Y_n(x_n)
Y_{n-1}(\bar x_{n-1})\cd Y_1(\bar x_1)
\fsquare(5mm,1)\,\,\big\vert\,\,x_i,\bar x_i\in\bbC^\times\},\\
\cV_2(\ge)&\seteq &\{v_2(y)
=Y_n(y_n)\cd Y_1(y_1)Y_0(y_0)Y_{1}(\ovl y_1)\cd
Y_{n-1}(\ovl y_{n-1})\fsquare(5mm,\ovl n)
\,\big\vert\,y_i,\ovl y_i\in\bbC^\times\}.
\end{eqnarray*}
For any $x\in (\bbC^\times)^{2n}$ 
there exist a unique $y\in (\bbC^\times)^{2n}$ and 
a rational function $a(x)$ 
such that $v_2(y)=a(x)v_1(x)$.
\end{enumerate}
\end{thm}

Now, using this theorem, we define the rational mapping
\begin{equation}
\begin{array}{cccccccccc}
\ovl\sigma\cl \cV(\ge)&\longrightarrow &\cV(\ge),&&&
\ovl\sigma\cl\cV(\ge)
&\longrightarrow &\cV_2(\ge),&\\
v(x)&\mapsto &v(y)&(\ge\ne A_{2n}^{(2)\,\dagger}),&\quad&
v_1(x)&\mapsto &v_2(y)&
(\ge=A_{2n}^{(2)\,\dagger}),
\end{array}
\end{equation}

\begin{thm}
\label{aff-geo}
The rational mapping $\ovl\sigma$ is birational. If we define 
\begin{equation}
\begin{cases}
e_0^c\seteq \ovl\sigma^{-1}\circ 
e_{\sigma(0)}^c\circ\ovl\sigma
,\q \vep_0\seteq \vep_{\sigma(0)}\circ\ovl\sigma,\q
\gamma_0\seteq \gamma_{\sigma(0)}\circ\ovl\sigma,&
\text{for $\ge\ne A_{2n}^{(2)\,\dagger}$,}\\
e_0^c\seteq \ovl\sigma^{-1}\circ 
e_0^c\circ\ovl\sigma
,\q \vep_0\seteq \vep_0\circ\ovl\sigma,\q
\gamma_0\seteq \gamma_0\circ\ovl\sigma,&
\text{for $\ge=A_{2n}^{(2)\,\dagger}$}.
\end{cases}
\label{ese}
\end{equation}
then $(\cV(\ge),\{e_i\}_{i\in I},
\{\gamma_i\}_{i\in I}, \{\vep_i\}_{i\in I})$
is an affine $\ge$-geometric crystal.
\end{thm}
{\sl Remark.}
In the case $\ge=A_{2n}^{(2)\,\dagger}$,
$\cV_2(\ge)$ 
has a $\ge_{I\setminus\{n\}}$-geometric crystal structure. 
Thus, $e_0,\gamma_0,\vep_0$ 
are well-defined on $\cV_2(\ge)$.

\medskip
The following lemma is obvious and it shows 
Theorem \ref{aff-geo} partially.
\begin{lem}
\label{io}
Suppose that $\ge\ne
A_{2n}^{(2)\,\dagger}$. 
If there exists 
$\osigma$ as above and 
\begin{equation}
e_{\sigma(i)}^c=\osigma\circ e_i^c\circ
\osigma^{-1}, \q\gamma_{i}
=\gamma_{\sigma(i)}\circ\osigma,\q
\vep_{i}=\vep_{\sigma(i)}\circ
\osigma, 
\label{cond-i}
\end{equation}
for $i\ne \sigma^{-1}(0),0$,  then we obtain
\bnum
\item
\begin{eqnarray*}
&&e_0^{c_1}e_i^{c_2}=e_i^{c_2}e_0^{c_1}\q
\text{if }a_{0i}=a_{i0}=0,\\
&&e^{c_1}_{0}e^{c_1c_2}_{i}e^{c_2}_{0}
=e^{c_2}_{1}e^{c_1c_2}_{0}e^{c_2}_{1}
\q\text{if }a_{0i}=a_{i0}=-1,\\
&&e^{c_1}_{i}e^{c^2_1c_2}_{j}e^{c_1c_2}_{i}
e^{c_2}_{j}
=e^{c_2}_{j}e^{c_1c_2}_{i}e^{c^2_1c_2}_{j}
e^{c_1}_{i}\q
{\rm if }\,\,a_{0i}=-2, \,\,a_{i0}=-1
\end{eqnarray*}
\item
$\gamma_0(e_i^c(v(x)))=c^{a_{i0}}
\gamma_0(v(x))$ and 
$\gamma_i(e_0^c(v(x)))=c^{a_{0i}}
\gamma_i(v(x))$.
\item
$\vep_0(e_0^c(v(x)))=c^{-1}\vep_0(v(x))$.
\enum
\end{lem}
{\sl Proof.} 
For example, we have
\begin{eqnarray*}
\gamma_0(e_i^c(v(x)))
&=&\gamma_{\sigma(0)}
(\osigma e_i^c\osigma^{-1}
(\osigma(v(x))))\\
&=&\gamma_{\sigma(0)}
(e_{\sigma(i)}^c(\osigma(v(x))))
=c^{a_{\sigma(i),\sigma(0)}}
\gamma_{\sigma(0)}(\osigma(v(x)))\\
&=&c^{a_{i,0}}\gamma_0(v(x)),
\end{eqnarray*}
where we use
$a_{\sigma(0),\sigma(i)}=a_{0i}$ in the last equality. 
The other assertions are obtained similarly.\qed

In the rest of this section, we shall prove
Theorems~\ref{birat} and \ref{aff-geo}
in each case.

\subsection{$A^{(1)}_n$-case}

We have $w_1\seteq s_ns_{n-1}\cd s_2s_1$, and
\[
{\cal V}(\TY(A,1,n))
\seteq \{Y_n(x_n)\cd Y_2(x_2)Y_1(x_1)\fsquare(5mm,1)
\,\big\vert\,x_i\in \bbC^\times\}\subset W(\varpi_1).
\]

Since $\exp(c^{-1}f_i)=1+c^{-1}f_i$ on 
$W(\varpi_1)$,  
$v(x)=Y_n(x_n)\cd Y_2(x_2)Y_1(x_1)\fsquare(5mm,1)$
is explicitly written as
\[
 v(x)=v(x_1,\cd,x_n)=
\left(\sum_{i=1}^nx_i\fsquare(5mm,i)\right)
+\fsquare(6mm,n+1).
\]
Let  $\sigma\cl\al_k\mapsto \al_{k+1}$ 
$(k\in I)$ be 
the  Dynkin diagram automorphism for $\TY(A,1,n)$, 
which gives rise to the automorphism 
$\sigma\cl W(\varpi_1)\rightarrow W(\varpi_1)$.
We have
$\sigma\fsquare(5mm,i)=\fsquare(5mm,i+1).$
Then, we obtain
\[
 \sigma(v(x))=\fsquare(5mm,1)+
\sum_{i=1}^nx_i\fsquare(5mm,i+1).
\]
Then the equation $v(y)=a(x)\sigma(v(x))$, {\em i.e.}
\[
 \sum_{i=1}^ny_i\fsquare(5mm,i)
+\fsquare(6mm,n+1)
=a(x)(\fsquare(5mm,1)+\sum_{i=1}^nx_i\fsquare(5mm,i+1))
\]
is solved by
\begin{equation}
a(x)=\frac{1}{x_n}, \q y_1=\frac{1}{x_n},\q
y_i=\frac{x_{i-1}}{x_n} \,\,(i=2,\cd,n), 
\end{equation}
that is
\begin{equation}
\ovl\sigma(v(x_1\cd,x_n))
=v(\frac{1}{x_n},\frac{x_1}{x_n},\cd,
\frac{x_{n-1}}{x_n}).
\end{equation}
The $A_n$-geometric crystal structure on $\cV(\TY(A,1,n))$ 
induced from the one on $B^-_{w_1}$
is given by:
\begin{eqnarray}
&&e_i^c(v(x_1,\cd,x_n))=v(x_1,\cd, cx_i,\cd,x_n)
\q(i=1,\cd,n), \label{ei-A}\\
&&\gamma_1(v(x))=\frac{x_1^2}{x_2},\q
\gamma_i(v(x))=\frac{x_i^2}{x_{i-1}x_{i+1}}\q
(i=2,\cd,n-1),\q
\gamma_n(v(x))=\frac{x_n^2}{x_{n-1}},\\
&&\vep_i(v(x))=\frac{x_{i+1}}{x_i}\q
(i=1,\cd, n-1),\q
\vep_n(v(x))=\frac{1}{x_n}.
\end{eqnarray}
Then we have 
\[
 \vep_{i+1}(\osigma(v(x)))=
\begin{cases}
\frac{x_{i+1}}{x_{i}}&\text{if }i=1\cd,n-2,\\
\frac{x_n}{x_{n-1}}&\text{if }i=n-1,
\end{cases}
\]
which implies 
$\vep_{\sigma(i)}(\osigma(v(x)))=\vep_i(v(x))$,
and then we completed the 
proof of Theorem \ref{birat}.

\medskip
Now, we define $e_0^c,\gamma_0$ and $\vep_0$ by 
\begin{equation}
e_0^c\seteq \ovl\sigma^{-1}\circ e_1^c
\circ \ovl\sigma,\q
\gamma_0\seteq \gamma_1\circ\ovl\sigma,\q
\vep_0\seteq \vep_1\circ\ovl\sigma.
\end{equation}
Their explicit forms are
\begin{eqnarray}
&&e_0^c(v(x))=v(\frac{x_1}{c},\frac{x_2}{c},\cd,
\frac{x_n}{c}),\\
&&\gamma_0(v(x))=\frac{1}{x_1x_n},\q
\vep_0(v(x))=x_1.
\label{ge0}
\end{eqnarray}
Thus, we can check \eqref{cond-i} easily,
and then Lemma \ref{io} 
reduces the proof of Theorem \ref{aff-geo} to
the statements:
\begin{eqnarray}
&&e^{c_1}_{0}e^{c_1c_2}_{n}e^{c_2}_{0}
=e^{c_2}_{n}e^{c_1c_2}_{0}e^{c_2}_{n},\\
&&\gamma_0(e_n^c(v(x)))=c^{-1}
\gamma_0(v(x)),\q
\gamma_n(e_0^c(v(x)))=c^{-1}
\gamma_n(v(x)).
\end{eqnarray}
These are immediate from (\ref{ei-A})--
(\ref{ge0}). Thus, we obtain 
Theorem \ref{aff-geo}.

\subsection{$B^{(1)}_n$-case}
\label{b1n}
We have $w_1=s_1\cd s_{n-1}s_ns_{n-1}\cd s_1$, and
\[
 {\cV}(\TY(B,1,n))
\seteq \{v(x)=Y_1(x_1)\cd Y_n(x_n)
Y_{n-1}(\bar x_{n-1})\cd Y_1(\bar x_1)
\fsquare(5mm,1)\,\,\big\vert\,\,x_i,\bar x_i\in\bbC^\times\}
\]
It follows from the explicit description of $W(\varpi_1)$
as in \ref{bn-w1}
that $y_i(c^{-1})=
\exp(c^{-1}f_i)$ on $W(\varpi_1)$ can be written as:
\[
 \exp(c^{-1}f_i)=\begin{cases}
1+c^{-1}f_i&i\ne n,\\
1+c^{-1}f_n+\dfrac{1}{2c^2}f_n^2&i= n.
\end{cases}
\]
Therefore, we have 
\begin{eqnarray*}
&&v(x_1,\cd,x_n,\ovl x_{n-1},\cd,\ovl x_1)
=\left(\sum_{i=1}^n\xi_i(x)\fsquare(0.5cm,i)\right)
+x_n\fsquare(0.5cm,0)
+\left(\sum_{i=2}^{n}x_{i-1}\fsquare(0.5cm,\ovl i)
\right)
+\fsquare(0.5cm,\ovl 1),\\
&&\qq\qq\q{\rm where}\q \xi_i(x)\seteq \begin{cases}
x_1\ovl x_1&i=1\\
\dfrac{x_{i-1}\ovl x_{i-1}+x_i\ovl x_i}{x_{i-1}}&i\ne1,n\\
\dfrac{x_{n-1}\ovl x_{n-1}+x^2_n}{x_{n-1}}&i=n
\end{cases}
\end{eqnarray*}
Since $\sigma\fsquare(0.5cm,1)
=\fsquare(0.5cm,\ovl 1)$,
$\sigma\fsquare(0.5cm,\ovl 1)
=\fsquare(0.5cm,1)$ and 
$\sigma\fsquare(0.5cm,k)
=\fsquare(0.5cm,k)$ otherwise,
we have 
\[
\sigma(v(x))
=\fsquare(0.5cm,1)
+\left(\sum_{i=2}^n \xi_i(x)\fsquare(0.5cm,i)
\right)
+x_n\fsquare(0.5cm,0)
+\left(\sum_{i=2}^{n}x_{i-1}\fsquare(0.5cm,\ovl i)
\right)
+x_1\ovl x_1\fsquare(0.5cm,\ovl 1),
\]
The equation $v(y)=a(x)\sigma(v(x))$
($x,y\in (\bbC^\times)^{2n-1}$)
has a unique solution:
\begin{equation}
a(x)=\frac{1}{x_1\ovl x_1},\q
y_i=a(x)x_i=
\frac{x_i}{x_1\ovl x_1}\q(1\leq i\leq n),\q
\ovl y_i=a(x)\ovl x_i
=\frac{\ovl x_i}{x_1\ovl x_i}\q
(1\leq i<n).
\end{equation}
Hence we have the rational mapping:
\begin{equation}
\osigma(v(x))\seteq v(y)=
v\left(\frac{x_1}{x_1\ovl x_1},
\frac{x_2}{x_1\ovl x_1},\cd,
\frac{x_n}{x_1\ovl x_1},
\frac{\ovl x_{n-1}}{x_1\ovl x_i},
\cd,
\frac{\ovl x_1}{x_1\ovl x_i}\right).
\label{bn-osigma}
\end{equation}
By the explicit form of 
$\osigma$ in (\ref{bn-osigma}), 
we have $\osigma^2={\rm id}$, which means 
that the morphism $\osigma$ is birational.
In this case, the second condition in 
Theorem \ref{birat} is trivial since 
$\sigma(i)=i$ if $i,\,\sigma(i)\ne 0$.
Thus, the proof of 
Theorem \ref{birat} in this case
is completed.

Now, we set 
$e_0^c\seteq \osigma\circ e_1^c\circ\osigma$, 
$\gamma_0\seteq \gamma_1\circ\osigma$ and 
$\vep_0\seteq \vep_1\circ\osigma$.

The explicit forms of $e_i$'s, 
$\vep_i$'s and $\gamma_i$'s are:
\begin{eqnarray*}
e_0^c\cl
&x_1\mapsto&x_1\frac{cx_1\bar x_1+x_2\bar x_2}
{c(x_1\bar x_1+x_2\bar x_2)}\qq
x_i\mapsto\frac{x_i}{c}\,\,(2\leq i\leq n)\\
&\bar x_1\mapsto&\bar x_1\frac{x_1\bar x_1+x_2\bar x_2}
{cx_1\bar x_1+x_2\bar x_2},\qq
\bar x_i\mapsto\frac{\bar x_i}{c}\,\,(2\leq i\leq n-1),\\
\hspace{-20pt}e_i^c\cl&x_i\mapsto&
x_i\frac{cx_i\ovl x_i+x_{i+1}\bar x_{i+1}}
{x_i\ovl x_i+x_{i+1}\ovl x_{i+1}},\q
\ovl x_i\mapsto
\ovl x_i\frac{c(x_i\ovl x_i+x_{i+1}\bar x_{i+1})}
{cx_i\ovl x_i+x_{i+1}\ovl x_{i+1}},\\
&x_j\mapsto& x_j,\q \ovl x_j\mapsto \ovl x_j\,\;
(j\ne i)
\qq\qq(1\leq i<n-1),\\
e_{n-1}^c\cl&x_{n-1}\mapsto&
x_{n-1}\frac{cx_{n-1}\ovl x_{n-1}+x_{n}^2}
{x_{n-1}\ovl x_{n-1}+x_{n}^2},\q
\ovl x_{n-1}\mapsto
\ovl x_{n-1}
\frac{c(x_{n-1}\ovl x_{n-1}+x_{n}^2)}
{cx_{n-1}\ovl x_{n-1}+x_{n}^2},\\
&x_j\mapsto& x_j,\q \ovl x_j\mapsto \ovl x_j\,\,
(j\ne n-1),\\
e^c_n\cl&x_n\mapsto& cx_n,\qq
x_j\mapsto x_j\q \ovl x_j\mapsto \ovl x_j\,\,
(j\ne n),
\end{eqnarray*}
\begin{eqnarray*}
&&\vep_0(v(x))=
\frac{x_1\ovl x_1+x_2\ovl x_2}{x_1},\q
\vep_1(v(x))=\frac{1}{x_1}
\left(1+\frac{x_{2}\ovl x_{2}}{x_{1}\ovl x_1}\right),\\
&&\vep_i(v(x))=\frac{x_{i-1}}{x_i}
\left(1+\frac{x_{i+1}\ovl x_{i+1}}{x_{i}\ovl x_i}\right)\,\,(2\leq i\leq n-2),
\\
&&\vep_{n-1}(v(x))=\frac{x_{n-2}}{x_{n-1}}
\left(1+\frac{x_n^2}{x_{n-1}\ovl x_{n-1}}\right),\,\,
\vep_n(v(x))=\frac{x_{n-1}}{x_n},
\end{eqnarray*}
\begin{eqnarray*}
&&\hspace{-10pt}\gamma_0(v(x))=\frac{1}{x_2\ovl x_2},\,\,
\gamma_1(v(x))=\frac{(x_1\ovl x_1)^2}{x_2\ovl x_2},\,\,
\gamma_i(v(x))=
\frac{(x_i\ovl x_i)^2}{x_{i-1}\ovl x_{i-1}x_{i+1}\ovl x_{i+1}}
\,(2\leq i\leq n-2),\\
&&\gamma_{n-1}(v(x))=
\frac{(x_{n-1}\ovl x_{n-1})^2}{x_{n-2}\ovl x_{n-2}x_n^2},\q
\gamma_n(v(x))=\frac{x_n^2}{x_{n-1}\ovl x_{n-1}}.
\end{eqnarray*}
Since $\sigma(i)=i$ for $i\ne0,1$,
the condition (\ref{cond-i}) in Lemma \ref{io}
can be easily seen by (\ref{bn-osigma}) and by the 
explicit form of $e_i$, $\gamma_i$ and $\vep_i$
($i\in I$).
Thus, in order to prove Theorem \ref{aff-geo}, 
it suffices to show that 
\begin{eqnarray}
&&e_0^{c_1}e_1^{c_2}=e_1^{c_2}e_0^{c_1},
\label{b1}\\
&&\gamma_0(e_1^c(v(x)))=\gamma_0(v(x)),\q
\gamma_1(e_0^c(v(x)))=\gamma_1(v(x)).
\label{b2}
\end{eqnarray}
It follows from the explicit formula above that
\[
 e_0^{c_1}e_1^{c_2}(v(x))=e_1^{c_2}e_0^{c_1}(v(x))
=v\left(
x_1\frac{c_1c_2x_1\ovl x_1+x_2\ovl x_2}
{c_1(x_1\ovl x_1+x_2\ovl x_2)},\frac{x_2}{c_1},\cd,
\frac{\ovl x_2}{c_1},
\ovl x_1\frac{c_2(x_1\ovl x_1+x_2\ovl x_2)}
{c_1c_2x_1\ovl x_1+x_2\ovl x_2}
\right),
\]
which implies (\ref{b1}). 
We get  (\ref{b2}) 
immediately from the formula above and we 
complete the proof of Theorem \ref{aff-geo}
for $\TY(B,1,n)$.

\subsection{$C^{(1)}_n$-case}

We have $w_1=s_0s_1\cd s_ns_{n-1}\cd s_1$, and
\[
 {\cV}(\TY(C,1,n))\seteq \{v(x)=Y_0(x_0)Y_1(x_1)\cd Y_n(x_n)
Y_{n-1}(\bar x_{n-1})\cd Y_1(\bar x_1)
\fsquare(5mm,1)\,\,\big\vert\,\,x_i,\bar x_i\in\bbC^\times\}.
\]
Due to the explicit description of $\W1$ in \ref{cn-w1},
we obtain $y_i(c^{-1})=\exp(c^{-1}f_i)=
1+c^{-1}f_i$ on $\W1$. Hence 
we have 
\begin{eqnarray*}
&&v(x_0,x_1,\cd,x_n,\ovl x_{n-1},\cd,\ovl x_1)\\
&&\qq =\left(\sum_{i=1}^n\xi_i\fsquare(0.5cm,i)\right)
+\left(\sum_{i=1}^{n}x_{i-1}\fsquare(0.5cm,\ovl i)
\right)
\text{ where } \xi_i\seteq \begin{cases}
\frac{x_{0}+x_1\ovl x_1}{x_{0}}&i=1\\
\frac{x_{i-1}\ovl x_{i-1}+x_i\ovl x_i}{x_{i-1}}&i\ne1,n\\
\frac{x_{n-1}\ovl x_{n-1}+x_n}{x_{n-1}}&i=n
\end{cases}
\end{eqnarray*}
The explicit forms of $\vep_i(x)$ $(1\leq i\leq n)$ are:
\eq
&&\vep_i(v(x))=\frac{x_{i-1}}{x_i}
\left(1+\frac{x_{i+1}\ovl x_{i+1}}{x_{i}\ovl x_i}\right)\,
(1\leq i\leq n-2),\nn\\
&&\vep_{n-1}(v(x))=\frac{x_{n-2}}{x_{n-1}}\left(
1+\frac{x_n}{x_{n-1}\ovl x_{n-1}}\right),\,
\vep_n(v(x))=\frac{x_{n-1}^2}{x_n}.\nn
\eneq
Since $\sigma\fsquare(5mm,i+1)
=\fsquare(6mm,\ovl{n-i})$ and 
$\sigma\fsquare(5mm,\ovl{i+1})=\fsquare(6mm,{n-i})$
($0\leq i<n$), 
we have 
\[
\sigma(v(x))
=\left(\sum_{i=1}^n x_{n-i}\fsquare(0.5cm,i)\right)
+\left(\sum_{i=1}^{n}\xi_{n-i+1}\fsquare(0.5cm,\ovl i)
\right).
\]
We obtain a unique solution of 
the equations $v(y)=a(x)\sigma(v(x))$, $\vep_{n-1}(v(y))=\vep_1(v(x))$
($x,y\in (\bbC^\times)^{2n}$):
\begin{eqnarray*}
&&a(x)=\frac{\ovl x_{n-1}}{x_n}+\frac{1}{x_{n-1}},\\
&&y_0=x_n
\left(\frac{\ovl x_{n-1}}{x_n}+\frac{1}{x_{n-1}}\right)^2\\
&&y_i=
\left(\ovl x_{n-i-1}+\frac{x_{n-i}\ovl x_{n-i}}{x_{n-i-1}}\right)
\left(\frac{\ovl x_{n-1}}{x_n}+\frac{1}{x_{n-1}}\right)\q
(1\leq i<n)\\
&&y_{n-1}=
\left(1+\frac{x_1\ovl x_1}{x_0}\right)
\left(\frac{\ovl x_{n-1}}{x_n}+\frac{1}{x_{n-1}}\right)\\
&&y_n=x_0
\left(\frac{\ovl x_{n-1}}{x_n}+\frac{1}{x_{n-1}}\right)^2
\\
&&\ovl y_{n-1}
=\frac{(x_{n-1}\ovl x_{n-1}+x_n)x_0x_1\ovl x_1}
{(x_0+x_1\ovl x_1)x_{n-1}x_n}
\\
&&\ovl y_i
=\frac{(x_{n-1}\ovl x_{n-1}+x_n)x_{n-i-1}x_{n-i}\ovl x_{n-i}}
{(x_{n-i-1}\ovl x_{n-i-1}+x_{n-i}\ovl x_{n-i})x_{n-1}x_n}
\q
(1\leq i<n-1).
\end{eqnarray*}
Now we have the rational mapping
$\osigma\cl\cV(\TY(C,1,n))\longrightarrow \cV(\TY(C,1,n))$
defined by $v(x)\mapsto v(y)$. By the 
above explicit from of $y$, we have 
$\osigma^2=$id, which means that 
$\osigma$ is birational. 

By direct calculations, we have 
$\vep_{n-i}(v(y))=\vep_i(v(x))$ ($1\leq i\le n-1$),
and we complete
the proof of Theorem \ref{birat} for $\TY(C,1,n)$.

Let us define $e^c_0\seteq \osigma\circ e^c_n\circ\osigma$
($\osigma^2={\rm id}$), $\gamma_0\seteq \gamma_n\circ\osigma$
and $\vep_0\seteq \vep_n\circ\osigma$. 
The explicit forms of $e_i$, $\gamma_i$ and $\vep_0$ are 
\begin{eqnarray*}
e_0^c\cl&x_0\mapsto&x_0\frac{(cx_0+x_1\bar x_1)^2}
{c(x_0+x_1\bar x_1)^2}\qq
x_i\mapsto{x_i}\frac{cx_0+x_1\bar x_1}
{c(x_0+x_1\bar x_1)}\,\,(1\leq i\leq n-1)\\
&x_n\mapsto&{x_n}\frac{(cx_0+x_1\bar x_1)^2}
{c^2(x_0+x_1\bar x_1)^2},\qq
\bar x_i\mapsto{\bar x_i}\frac{cx_0+x_1\bar x_1}
{c(x_0+x_1\bar x_1)}\,\,(1\leq i\leq n-1),\\
e_i^c\cl&x_i\mapsto&
x_i\frac{cx_i\ovl x_i+x_{i+1}\bar x_{i+1}}
{x_i\ovl x_i+x_{i+1}\ovl x_{i+1}}\qq
\ovl x_i\mapsto
\ovl x_i\frac{c(x_i\ovl x_i+x_{i+1}\bar x_{i+1})}
{cx_i\ovl x_i+x_{i+1}\ovl x_{i+1}},\\
&x_j\mapsto& x_j,\q \ovl x_j\mapsto \ovl x_j\q(j\ne i),
\qq\qq\qq \q(1\leq i<n-1),\\
e_{n-1}^c\cl&x_{n-1}\mapsto&
x_{n-1}\frac{cx_{n-1}\ovl x_{n-1}+x_{n}}
{x_{n-1}\ovl x_{n-1}+x_{n}}\qq
\ovl x_{n-1}\mapsto
\ovl x_{n-1}\frac{c(x_{n-1}\ovl x_{n-1}+x_{n})}
{cx_{n-1}\ovl x_{n-1}+x_{n}},\\
&x_j\mapsto& x_j,\q \ovl x_j\mapsto \ovl x_j\q(j\ne n-1),\\
e^c_n\cl&x_n\mapsto& cx_n,\q
x_j\mapsto x_j,\q \ovl x_j\mapsto \ovl x_j\q(j\ne n),
\end{eqnarray*}
\begin{eqnarray*}
&&\gamma_0(v(x))=\frac{x_0^2}{(x_1\ovl x_1)^2},\q
\gamma_1(v(x))=\frac{(x_1\ovl x_1)^2}{x_0x_2\ovl x_2},\q
\gamma_i=
\frac{(x_i\ovl x_i)^2}{x_{i-1}\ovl x_{i-1}x_{i+1}\ovl x_{i+1}}\q
(2\leq i\leq n-2),\\
&&\gamma_{n-1}(v(x))=
\frac{(x_{n-1}\ovl x_{n-1})^2}{x_{n-2}\ovl x_{n-2}x_n},\q
\gamma_n(v(x))=\frac{x_n^2}{(x_{n-1}\ovl x_{n-1})^2},\\
&&\vep_0(v(x))=\frac{1}{x_0}
\left(1+\frac{x_1\ovl x_1}{x_0}\right)^2.
\end{eqnarray*}
Let us check the condition (\ref{cond-i}) 
in Lemma \ref{io}. 
The following are useful for this purpose:
\begin{eqnarray}
&&y_i\ovl y_i=a(x)^2x_{n-i}\ovl x_{n-i},\q
y_0=a(x)^2x_n,\label{u1}\\
&&a(v(y))=a(\osigma(v(x)))=\frac{1}{a(v(x))}.
\label{u2}
\end{eqnarray}
Using these we can easily check 
the two conditions 
$\gamma_{i}
=\gamma_{\sigma(i)}\circ\osigma$ and 
$\vep_{i}=\vep_{\sigma(i)}\circ
\osigma$.
The condition 
$e_{\sigma(i)}^c=\osigma\circ e_i^c\circ
\osigma^{-1}$ for $i=2,\cd,n-2$ is also 
immediate from (\ref{u1}) and (\ref{u2}).
Next let us see the case $i=1,n-1$.
We have 
\[
a(e_{n-1}^c(v(y)))=
\frac{y_n+cy_{n-1}\ovl y_{n-1}}{y_{n-1}y_n}\cdot
\frac{x_1\ovl x_1+x_0}{cx_1\ovl x_1+x_0}=\frac{1}{a(v(x))}.
\]
Using this, we can get 
$e_{n-i}^c=\osigma\circ e_1^c\circ
\osigma^{-1}$ and then 
$e_{1}^c=\osigma\circ e_{n-1}^c\circ
\osigma^{-1}$ since $\osigma^2={\rm id}$.
Now, it remains to show that 
\begin{eqnarray*}
&&e_0^{c_1}e_n^{c_2}=e_n^{c_2}e_0^{c_1},\q
\vep_0(e_n^c(v(x)))=\vep_0(v(x)),\q
\vep_n(e_0^c(v(x)))=\vep_n(v(x)).
\end{eqnarray*}
They easily follow from the explicit form of 
$e_0^c$. Thus, the proof of Theorem \ref{aff-geo}
in this case is completed. 
\subsection{$D^{(1)}_n$-case}
\label{d1n}
We have $w_1=s_1s_2\cd s_{n-1}s_ns_{n-2}s_{n-3}\cd s_2s_1$, and
\[
 {\cV}(\TY(D,1,n))\seteq 
\{v(x)=Y_1(x_1)\cd Y_{n-1}(x_{n-1})Y_n(x_n)
Y_{n-2}(\bar x_{n-2})\cd Y_1(\bar x_1)\fsquare(5mm,1)
\,\,\big\vert\,
\,x_i,\bar x_i\in\bbC^\times\}.
\]
It follows from the explicit form of $\W1$ in \ref{dn-w1}
that $y_i(c^{-1})=\exp(c^{-1}f_i)=
1+c^{-1}f_i$ on $\W1$. Thus, 
we have 
\begin{eqnarray*}
&&v(x)
=\left(\sum_{i=1}^{n-1}\xi_i(x)\fsquare(0.5cm,i)\right)
+x_n\fsquare(0.5cm,n)
+\left(\sum_{i=2}^{n}x_{i-1}
\fsquare(0.5cm,\ovl i)\right)
+\fsquare(0.5cm,\ovl 1),\\
&&\qq\qq\qq\q{\rm where}\q \xi_i(x)\seteq \begin{cases}
x_1\ovl x_1&i=1\\
\frac{x_{i-1}\ovl x_{i-1}+x_i\ovl x_i}{x_{i-1}}&i\ne1,n-1\\
\frac{x_{n-2}\ovl x_{n-2}+x_{n-1}x_n}{x_{n-2}}&i=n-1
\end{cases}
\end{eqnarray*}
Since $\sigma\fsquare(0.5cm,1)=\fsquare(0.5cm,\ovl 1)$,
$\sigma\fsquare(0.5cm,\ovl 1)=\fsquare(0.5cm,1)$ and 
$\sigma\fsquare(0.5cm,k)=\fsquare(0.5cm,k)$ otherwise,
we have 
\[
\sigma(v(x))
=\fsquare(0.5cm,1)+
\left(\sum_{i=2}^n \xi_i(x)\fsquare(0.5cm,i)\right)
+x_n\fsquare(0.5cm,n)
+\left(\sum_{i=1}^{n-1}x_{i-1}\fsquare(0.5cm,\ovl i)
\right)
+\xi_1\fsquare(0.5cm,\ovl 1).
\]
Then the equation $v(y)=a(x)\sigma(v(x))$
($x,y\in (\bbC^\times)^{2n-2}$)
has the following unique
solution:
\begin{equation}
a(x)=\frac{1}{x_1\ovl x_1},\q
y_i=a(x)x_i=
\frac{x_i}{x_1\ovl x_1}\q(1\leq i\leq n),\q
\ovl y_i=a(x)\ovl x_i
=\frac{\ovl x_i}{x_1\ovl x_i}\q
(1\leq i\leq n-2).
\end{equation}
We define the rational mapping $\osigma\cl
\cV(\TY(D,1,n))\longrightarrow 
\cV(\TY(D,1,n))$ by 
\begin{equation}
\osigma(v(x))=
v\left(\frac{x_1}{x_1\ovl x_1},
\frac{x_2}{x_1\ovl x_1},\cd,
\frac{x_n}{x_1\ovl x_1},
\frac{\ovl x_{n-1}}{x_1\ovl x_i},
\cd,
\frac{\ovl x_1}{x_1\ovl x_i}\right).
\label{dn-osigma}
\end{equation}
It is immediate from 
(\ref{dn-osigma}) that 
$\osigma^2={\rm id}$, which implies
the birationality of 
the morphism $\osigma$.
In this case, the second condition in 
Theorem \ref{birat} is trivial since 
$\sigma(i)=i$ if $i,\,\sigma(i)\ne 0$.
Thus, the proof of 
Theorem \ref{birat} for $\TY(D,1,n)$
is completed.

Now, we set 
$e_0^c\seteq \osigma\circ e_1^c\circ\osigma$, 
$\gamma_0\seteq \gamma_1\circ\osigma$ and 
$\vep_0\seteq \vep_1\circ\osigma$.
The explicit forms of $e_i$, $\vep_i$ and $\gamma_i$ are:
\begin{eqnarray*}
e_0^c\cl
&x_1\mapsto&x_1\frac{cx_1\bar x_1+x_2\bar x_2}
{c(x_1\bar x_1+x_2\bar x_2)}\qq
x_i\mapsto\frac{x_i}{c}\,\,(2\leq i\leq n)\\
&\bar x_1\mapsto&\bar x_1\frac{x_1\bar x_1+x_2\bar x_2}
{cx_1\bar x_1+x_2\bar x_2},\qq
\bar x_i\mapsto\frac{\bar x_i}{c}\,\,(2\leq i\leq n-2),\\
e_i^c\cl&x_i\mapsto&
x_i\frac{cx_i\ovl x_i+x_{i+1}\bar x_{i+1}}
{x_i\ovl x_i+x_{i+1}\ovl x_{i+1}}\q
\ovl x_i\mapsto
\ovl x_i\frac{c(x_i\ovl x_i+x_{i+1}\bar x_{i+1})}
{cx_i\ovl x_i+x_{i+1}\ovl x_{i+1}},\\
&x_j\mapsto& x_j,\q \ovl x_j\mapsto\ovl x_j\q(j\ne i),
\qq\q(1\leq i\leq n-3),\\
e_{n-2}^c\cl&x_{n-2}\mapsto&
x_{n-2}\frac{cx_{n-2}\ovl x_{n-2}+x_{n-1}x_{n}}
{x_{n-2}\ovl x_{n-2}+x_{n-1}x_{n}},\q
\ovl x_{n-2}\mapsto
\ovl x_{n-2}\frac{c(x_{n-2}\ovl x_{n-2}+x_{n-1}x_{n})}
{cx_{n-2}\ovl x_{n-2}+x_{n-1}x_{n}},\\
&x_j\mapsto& x_j,\q \ovl x_j\mapsto\ovl x_j\q(j\ne n-2),\\
e^c_{n-1}\cl&x_{n-1}\mapsto& cx_{n-1},\q
x_j\mapsto x_j,\q \ovl x_j\mapsto\ovl x_j\q(j\ne n-1),\\
e^c_n\cl&x_n\mapsto& cx_n,\q
x_j\mapsto x_j,\q \ovl x_j\mapsto\ovl x_j\q(j\ne n),
\end{eqnarray*}
\begin{eqnarray*}
&&\vep_0(v(x))=
\frac{x_1\ovl x_1+x_2\ovl x_2}{x_1},\q
\vep_1(v(x))=\frac{1}{x_1}
\left(1+\frac{x_{2}\ovl x_{2}}{x_{1}\ovl x_1}\right),\q\\
&&\vep_i(v(x))=\frac{x_{i-1}}{x_i}
\left(1+\frac{x_{i+1}\ovl x_{i+1}}{x_{i}\ovl x_i}\right)\q
(2\leq i\leq n-3),\\
&&\vep_{n-2}(v(x))=\frac{x_{n-3}}{x_{n-2}}
\left(1+\frac{x_{n-1}x_n}{x_{n-2}\ovl x_{n-2}}\right),\q
\vep_{n-1}(v(x))=\frac{x_{n-2}}{x_{n-1}},\q
\vep_n(v(x))=\frac{x_{n-2}}{x_n},
\end{eqnarray*}
\begin{eqnarray*}
&&\hspace{-10pt}
\gamma_0(v(x))=\frac{1}{x_2\ovl x_2},\,\,
\gamma_1(v(x))=\frac{(x_1\ovl x_1)^2}{x_2\ovl x_2},\,\,
\gamma_i(v(x))=
\frac{(x_i\ovl x_i)^2}{x_{i-1}\ovl x_{i-1}x_{i+1}\ovl x_{i+1}}
\,(2\leq i\leq n-3),\\
&&\gamma_{n-2}(v(x))=\frac{(x_{n-2}\ovl x_{n-2})^2}
{x_{n-3}\ovl x_{n-3}x_{n-1}x_n},\q
\gamma_{n-1}(v(x))=
\frac{x_{n-1}^2}{x_{n-2}\ovl x_{n-2}},\q
\gamma_n(v(x))=\frac{x_n^2}{x_{n-2}\ovl x_{n-2}}.
\end{eqnarray*}
By these formulas, we can show Theorem \ref{aff-geo} 
for $\TY(D,1,n)$ similarly to the one for $\TY(B,1,n)$.

\subsection{$A^{(2)}_{2n-1}$-case}
\label{aon}

We have $w_1=s_1s_2\cd s_ns_{n-1}\cd s_2s_1$, and
\[
 {\cV}(\TY(A,2,2n-1))\seteq \{v(x)\seteq Y_1(x_1)\cd Y_n(x_n)
Y_{n-1}(\bar x_{n-1})\cd Y_1(\bar x_1)\fsquare(5mm,1)\,
\,\big\vert\,\,x_i,\bar x_i\in\bbC^\times\}.
\]
In this case, $y_i(c^{-1})=\exp(c^{-1}f_i)=
1+c^{-1}f_i$ on $\W1$, and
we have 
\begin{eqnarray*}
&&v(x_1,\cd,x_n,\ovl x_{n-1},\cd,\ovl x_1)\\
&&\qq =\left(\sum_{i=1}^n\xi_i\fsquare(0.5cm,i)\right)
+\left(\sum_{i=2}^{n}x_{i-1}\fsquare(0.5cm,\ovl i)
\right)
+\fsquare(0.5cm,\ovl 1),
\q{\rm where}\q \xi_i\seteq \begin{cases}
x_1\ovl x_1&i=1\\
\frac{x_{i-1}\ovl x_{i-1}+x_i\ovl x_i}{x_{i-1}}&i\ne1,n\\
\frac{x_{n-1}\ovl x_{n-1}+x_n}{x_{n-1}}&i=n
\end{cases}
\end{eqnarray*}
Since $\sigma\fsquare(0.5cm,1)=\fsquare(0.5cm,\ovl 1)$,
$\sigma\fsquare(0.5cm,\ovl 1)=\fsquare(0.5cm,1)$ and 
$\sigma\fsquare(0.5cm,k)=\fsquare(0.5cm,k)$ otherwise,
we have 
\[
\sigma(v(x))
=\fsquare(0.5cm,1)+\left(
\sum_{i=2}^n \xi_i\fsquare(0.5cm,i)\right)
+\left(\sum_{i=2}^{n}x_{i-1}\fsquare(0.5cm,\ovl i)
\right)
+x_1\ovl x_1\fsquare(0.5cm,\ovl 1),
\]
Solving $v(y)=a(x)\sigma(v(x))$
($x, y\in (\bbC^\times)^{2n-1}$), we obtain 
a unique solution:
\begin{equation}
a(x)=\frac{1}{x_1\ovl x_1},\q
y_i=a(x)x_i=
\frac{x_i}{x_1\ovl x_1},\q 
\ovl y_i=a(x)\ovl x_i
=\frac{\ovl x_i}{x_1\ovl x_i}\q
(1\leq i\leq n-1),\ 
y_n=\frac{x_n}{(x_1\ovl x_1)^2}.
\end{equation}
The we have 
\begin{equation}
\osigma(v(x))=
v\left(\frac{x_1}{x_1\ovl x_1},
\frac{x_2}{x_1\ovl x_1},\cd,
\frac{x_n}{(x_1\ovl x_1)^2},
\frac{\ovl x_{n-1}}{x_1\ovl x_i},
\cd,
\frac{\ovl x_1}{x_1\ovl x_i}\right).
\label{a2o-osigma}
\end{equation}
By the explicit form of 
$\osigma$ in (\ref{a2o-osigma}), 
we have $\osigma^2={\rm id}$, which means 
that the morphism $\osigma$ is birational.
In this case, the second condition in 
Theorem \ref{birat} is trivial since 
$\sigma(i)=i$ if $i,\,\sigma(i)\ne 0$.
Thus, the proof of 
Theorem \ref{birat} for $\TY(A,2,2n-1)$
is completed.

Now, we set 
$e_0^c\seteq \osigma\circ e_1^c\circ\osigma$, 
$\gamma_0\seteq \gamma_1\circ\osigma$ and 
$\vep_0\seteq \vep_1\circ\osigma$.
The explicit forms of $e_i$, $\vep_i$ and $\gamma_i$ are:
\begin{eqnarray*}
e_0^c\cl
&x_1\mapsto&x_1
\frac{cx_1\bar x_1+x_2\bar x_2}
{c(x_1\bar x_1+x_2\bar x_2)}\qq
x_i\mapsto\frac{x_i}{c}\,\,
(2\leq i\leq n-1),\q
x_n\mapsto \frac{x_n}{c^2}\\
&\bar x_1\mapsto&\bar x_1\frac{x_1\bar x_1+x_2\bar x_2}
{cx_1\bar x_1+x_2\bar x_2},\qq
\bar x_i\mapsto\frac{\bar x_i}{c}\,\,(2\leq i\leq n-1),\\
e_i^c\cl&x_i\mapsto&
x_i\frac{cx_i\ovl x_i+x_{i+1}\bar x_{i+1}}
{x_i\ovl x_i+x_{i+1}\ovl x_{i+1}}\qq
\ovl x_i\mapsto
\ovl x_i\frac{c(x_i\ovl x_i+x_{i+1}\bar x_{i+1})}
{cx_i\ovl x_i+x_{i+1}\ovl x_{i+1}},\\
&x_j\mapsto& x_j,\q \ovl x_j\mapsto \ovl x_j\q(j\ne i)
\qq(1\leq i<n-1),\\
e_{n-1}^c\cl&x_{n-1}\mapsto&
x_{n-1}\frac{cx_{n-1}\ovl x_{n-1}+x_{n}}
{x_{n-1}\ovl x_{n-1}+x_{n}}\qq
\ovl x_{n-1}\mapsto
\ovl x_{n-1}\frac{c(x_{n-1}\ovl x_{n-1}+x_{n})}
{cx_{n-1}\ovl x_{n-1}+x_{n}},\\
&x_j\mapsto& x_j,\q \ovl x_j\mapsto \ovl x_j\q(j\ne n-1),\\
e^c_n\cl&x_n\mapsto& cx_n,\q
x_j\mapsto x_j,\q \ovl x_j\mapsto \ovl x_j\q(j\ne n).
\end{eqnarray*}
\begin{eqnarray*}
&&\vep_0(v(x))=
\frac{x_1\ovl x_1+x_2\ovl x_2}{x_1},\q
\vep_1(v(x))=\frac{1}{x_1}
\left(1+\frac{x_{2}\ovl x_{2}}{x_{1}\ovl x_1}\right),\\
&&\vep_i(v(x))=\frac{x_{i-1}}{x_i}
\left(1+\frac{x_{i+1}\ovl x_{i+1}}{x_{i}\ovl x_i}\right)\q
(2\leq i\leq n-2),\\
&& \vep_{n-1}(v(x))=\frac{x_{n-2}}{x_{n-1}}
\left(1+\frac{x_n}{x_{n-1}\ovl x_{n-1}}\right),\q
\vep_n(v(x))=\frac{x_{n-1}^2}{x_n},
\end{eqnarray*}
\begin{eqnarray*}
&&\gamma_0(v(x))=\frac{1}{x_2\ovl x_2},\q
\gamma_1(v(x))=\frac{(x_1\ovl x_1)^2}{x_2\ovl x_2},\q
\gamma_i=
\frac{(x_i\ovl x_i)^2}{x_{i-1}\ovl x_{i-1}x_{i+1}\ovl x_{i+1}}
\q(2\leq i\leq n-2),\\
&&\gamma_{n-1}(v(x))=
\frac{(x_{n-1}\ovl x_{n-1})^2}{x_{n-2}\ovl x_{n-2}x_n},\q
\gamma_n(v(x))=\frac{x_n^2}{(x_{n-1}\ovl x_{n-1})^2}.
\end{eqnarray*}
We can show Theorem \ref{aff-geo} for $\TY(A,2,2n-1)$ 
similarly to the one for $\TY(B,1,n)$.

\subsection{$D^{(2)}_{n+1}$-case}
\label{d2n}

We have $w_1=s_0s_1\cd s_ns_{n-1}\cd s_2s_1$, and
\[
\hspace{-20pt}
 {\cV}(\TY(D,2,n+1))\seteq \{v(x)\seteq Y_0(x_0)Y_1(x_1)\cd Y_n(x_n)
Y_{n-1}(\bar x_{n-1})\cd Y_1(\bar x_1)
\fsquare(5mm,1)\,\,\big\vert\,\,x_i,\bar x_i\in\bbC^\times\}.
\]
It follows from the explicit description of $W(\varpi_1)$
as in \ref{d2n-w1}
that 
on $W(\varpi_1)$:
\[
y_i(c^{-1})= \exp(c^{-1}f_i)=\begin{cases}
1+c^{-1}f_i&i\ne 0,n,\\
1+c^{-1}f_i+\frac{1}{2c^2}f_i^2&i=0, n.
\end{cases}
\]
Then we have 
\[
v(x)
=\left(\sum_{i=1}^n\xi_i(x)\fsquare(0.5cm,i)\right)
+x_n\fsquare(0.5cm,0)+x_0\phi
+\left(\sum_{i=2}^{n}x_{i-1}\fsquare(0.5cm,\ovl i)
\right)
+x_0^2\fsquare(0.5cm,\ovl 1)\]
where
\[\xi_i(x)\seteq \begin{cases}
\frac{x_{0}^2+x_1\ovl x_1}{x_{0}^2}&i=1\\
\frac{x_{i-1}\ovl x_{i-1}+x_i\ovl x_i}{x_{i-1}}&i\ne1,n\\
\frac{x_{n-1}\ovl x_{n-1}+x_n^2}{x_{n-1}}&i=n
\end{cases}
\]
Since $\sigma\fsquare(5mm,i+1)
=\fsquare(5mm,\ovl{n-i})$ and 
$\sigma\fsquare(5mm,\ovl{i+1})=\fsquare(5mm,{n-i})$
($0\leq i<n$) and 
$\sigma\cl\fsquare(5mm,0)\leftrightarrow \phi$,
we also have 
\[
\sigma(v(x))
=\left(\sum_{i=1}^{n-1} x_{n-i}\fsquare(0.5cm,i)
\right)+x_0^2\fsquare(0.5cm,n)
+\left(\sum_{i=1}^{n}\xi_{n-i+1}\fsquare(0.5cm,\ovl i)
\right)
+x_n\phi+x_0\fsquare(0.5cm,0),
\]
Solving $v(y)=a(x)\sigma(v(x))$
($x, y\in (\bbC^\times)^{2n}$), we get a unique solution:
\begin{eqnarray*}
&&a(x)=\frac{x_{n-1}\ovl x_{n-1}+x_n^2}{x_{n-1}x_n^2},\\
&&y_0=
\frac{x_{n-1}\ovl x_{n-1}+x_n^2}{x_{n-1}x_n},\\
&&y_i=
\frac{(x_{n-i-1}\ovl x_{n-i-1}+x_{n-i}\ovl x_{n-i})
(x_{n-1}\ovl x_{n-1}+x_n^2)}{x_{n-i-1}x_{n-1}x_n^2}
\q(1\leq i<n),\\
&&y_{n-1}=
\frac{(x_0^2+x_1\ovl x_1)(x_{n-1}\ovl x_{n-1}+x_n^2)}
{x_0^2x_{n-1}x_n^2}
,\\
&&y_n=
\frac{x_0(x_{n-1}\ovl x_{n-1}+x_n^2)}{x_{n-1}x_n^2},\\
&&\ovl y_i
=\frac{(x_{n-1}\ovl x_{n-1}+x_n^2)x_{n-i+1}x_{n-i}\ovl x_{n-i}}
{(x_{n-i-1}\ovl x_{n-i-1}+x_{n-i}\ovl x_{n-i})x_{n-1}x_n^2}\q
(1\leq i\leq n-2),\\
&&\ovl y_{n-1}
=\frac{(x_{n-1}\ovl x_{n-1}+x_n^2)x_0^2x_{1}\ovl x_{1}}
{(x_0^2+x_{1}\ovl x_{1})x_{n-1}x_n^2}.
\end{eqnarray*}
Then we have the rational mapping 
$\osigma\cl\cV(\TY(D,2,n+1))\longrightarrow \cV(\TY(D,2,n+1))$
defined by $v(x)\mapsto v(y)$.
The explicit forms of $\vep_i$ ($1\leq i\leq n$) are 
as follows:
\begin{eqnarray*}
&&\hspace{-15pt}\vep_1(v(x))=\frac{x_0^2}{x_1}
\left(1+\frac{x_2\ovl x_2}{x_{1}\ovl x_1}\right),\q
\vep_n(v(x))=\frac{x_{n-1}}{x_n},\\
&&\hspace{-15pt}
\vep_i(v(x))=\frac{x_{i-1}}{x_i}
\left(1+\frac{x_{i+1}\ovl x_{i+1}}{x_{i}\ovl x_i}\right)\,\,
(2\leq i\leq n-2),\,\,
\vep_{n-1}(v(x))=\frac{x_{n-2}}{x_{n-1}}
\left(1+\frac{x_n^2}{x_{n-1}\ovl x_{n-1}}\right).
\end{eqnarray*}
Then we get easily that 
$\vep_{n-i}(v(y))=\vep_i(v(x))$ ($1\leq i\leq n-1$), which 
finishes the proof of Theorem \ref{birat} for $\TY(D,2,n+1)$.

Let us define $e^c_0\seteq \osigma\circ e^c_n\circ\osigma$
($\osigma^2={\rm id}$), $\gamma_0\seteq \gamma_n\circ\osigma$
and $\vep_0\seteq \vep_n\circ\osigma$. 
The explicit forms of $e_i$, $\gamma_i$ and $\vep_0$ are 
\begin{eqnarray*}
e_0^c\cl&x_0\mapsto&x_0\frac{c^2x_0^2+x_1\bar x_1}
{c(x_0^2+x_1\bar x_1)},\qq
x_i\mapsto{x_i}\frac{c^2x_0^2+x_1\bar x_1}
{c^2(x_0^2+x_1\bar x_1)}\,\,(1\leq i\leq n),\\
&
\bar x_i\mapsto&{\bar x_i}\frac{c^2x_0^2+x_1\bar x_1}
{c^2(x_0^2+x_1\bar x_1)}\,\,(1\leq i\leq n-1),\\
e_i^c\cl&x_i\mapsto&
x_i\frac{cx_i\ovl x_i+x_{i+1}\bar x_{i+1}}
{x_i\ovl x_i+x_{i+1}\ovl x_{i+1}},\qq
\ovl x_i\mapsto
\ovl x_i\frac{c(x_i\ovl x_i+x_{i+1}\bar x_{i+1})}
{cx_i\ovl x_i+x_{i+1}\ovl x_{i+1}},\\
&x_j\mapsto& x_j,\q \ovl x_j\mapsto \ovl x_j\q(j\ne i),
\qq\qq\qq\q(1\leq i<n-1),\\
e_{n-1}^c\cl&x_{n-1}\mapsto&
x_{n-1}\frac{cx_{n-1}\ovl x_{n-1}+x_{n}^2}
{x_{n-1}\ovl x_{n-1}+x_{n}^2},\q
\ovl x_{n-1}\mapsto
\ovl x_{n-1}\frac{c(x_{n-1}\ovl x_{n-1}+x_{n}^2)}
{cx_{n-1}\ovl x_{n-1}+x_{n}^2},\\
&x_j\mapsto& x_j,\q \ovl x_j\mapsto \ovl x_j\q(j\ne n-1),\\
e^c_n\cl&x_n\mapsto& cx_n,\q
x_j\mapsto x_j,\q \ovl x_j\mapsto \ovl x_j\q(j\ne n),
\end{eqnarray*}
\begin{eqnarray*}
&&\hspace{-15pt}
\gamma_0(v(x))=\frac{x_0^2}{x_1\ovl x_1},\q
\gamma_1(v(x))=\frac{(x_1\ovl x_1)^2}{x_0^2x_2\ovl x_2},\,\,
\gamma_i(v(x))=
\frac{(x_i\ovl x_i)^2}{x_{i-1}\ovl x_{i-1}x_{i+1}\ovl x_{i+1}}\,\,
(2\leq i\leq n-2),\\
&&\hspace{-15pt}\gamma_{n-1}(v(x))=
\frac{(x_{n-1}\ovl x_{n-1})^2}{x_{n-2}\ovl x_{n-2}x_n^2},\,\,
\gamma_n(v(x))=\frac{x_n^2}{x_{n-1}\ovl x_{n-1}},\\
&&\hspace{-15pt}\vep_0(v(x))=
\frac{x_0^2+x_1\ovl x_1}{x_0^3}.
\end{eqnarray*}
Then similarly to the case $\TY(C,1,n)$, we can show
Theorem \ref{aff-geo} for $\TY(D,2,n+1)$.

\subsection{$A^{(2)\,\dagger}_{2n}$-case}
\label{aen-1}

As in the beginning of this section, 
we have $w_1=s_0s_1\cd s_ns_{n-1}\cd s_2s_1$
and
\[
 {\cal V}(A_{2n}^{(2)\,\dagger})\seteq \{v_1(x)
=Y_0(x_0)Y_1(x_1)\cd Y_n(x_n)
Y_{n-1}(\bar x_{n-1})\cd Y_1(\bar x_1)\fsquare(5mm,1)
\,\,\big\vert\,\,x_i,\bar x_i\in\bbC^\times\}.
\]
By the explicit description of $W(\varpi_1)$
as in \ref{aen-v1}
on $W(\varpi_1)$ we have:
\[
y_i(c^{-1})=\exp(c^{-1}f_i)=\begin{cases}
1+c^{-1}f_i&i\ne n,\\
1+c^{-1}f_i+\frac{1}{2c^2}f_n^2&i= n.
\end{cases}
\]
Then we have 
\[
v_1(x)
=\left(\sum_{i=1}^{n}\xi_i(x)\fsquare(0.5cm,i)\right)
+x_n\fsquare(0.5cm,0)
+\left(\sum_{i=1}^{n}x_{i-1}\fsquare(0.5cm,\ovl i)
\right)
\q{\rm where}\q \xi_i(x)\seteq \begin{cases}
\frac{x_{0}+x_1\ovl x_1}{x_{0}}&i=1\\
\frac{x_{i-1}\ovl x_{i-1}+x_i\ovl x_i}{x_{i-1}}&i\ne1,n\\
\frac{x_{n-1}\ovl x_{n-1}+x_n^2}{x_{n-1}}&i=n
\end{cases}
\]
Next, for $w_2=s_ns_{n-1}\cd s_1
s_0s_1\cd s_{n-1}$ we set 
\[
\cV_2(A_{2n}^{(2)\,\dagger})\seteq \{v_2(y)
=Y_n(y_n)\cd Y_1(y_1)Y_0(y_0)Y_{1}(\ovl y_1)\cd
Y_{n-1}(\ovl y_{n-1})\fsquare(5mm,\ovl n)
\,\big\vert\,y_i,\ovl y_i\in\bbC^\times\}.
\]
Then we have 
\[
v_2(y)
=\left(\sum_{i=1}^{n-1} y_i\fsquare(0.5cm,i)\right)
+y_n^2\fsquare(0.5cm,n)
+y_n\fsquare(0.5cm,0)
+\left(\sum_{i=1}^{n}\eta_i(y)\fsquare(0.5cm,\ovl i)
\right)
\ {\rm where}\ \eta_i(y)\seteq \begin{cases}
\frac{y_{0}+y_1\ovl y_1}{y_{1}}&i=1\\
\frac{y_{i-1}\ovl y_{i-1}+y_i\ovl y_i}{y_{i}}&i\ne1,n\\
\frac{y_{n-1}\ovl y_{n-1}+y_n^2}{y_{n}^2}&i=n.
\end{cases}
\]
For $x\in(\bbC^\times)^{2n}$ 
there exist a unique
$y=(y_0,\cd,\ovl y_1)\in (\bbC^\times)^{2n}$ 
and $a(x)$ such 
that $v_2(y)=a(x)v_1(x)$.
They are given by
\begin{eqnarray*}
&&a(x)=\frac{x_{n-1}\ovl x_{n-1}+x_n^2}
{x_{n-1}x_n^2},\\
&&y_0=a(x)^2x_0=
\frac{x_0(x_{n-1}\ovl x_{n-1}
+x_n^2)^2}{(x_{n-1}x_n^2)^2},\\
&&y_1=a(x)\xi_1(x)=
\frac{(x_{n-1}\ovl x_{n-1}
+x_{n}^2)(x_0+x_1\ovl x_1)}
{x_0x_{n-1}x_n^2},\\
&&y_i=a(x)\xi_i(x)=
\frac{(x_{i-1}\ovl x_{i-1}
+x_{i}\ovl x_{i})
(x_{n-1}\ovl x_{n-1}+x_n^2)}
{x_{i-1}x_{n-1}x_n^2}
\q(1\leq i<n),\\
&&y_n=a(x)x_n=
\frac{x_{n-1}\ovl x_{n-1}
+x_n^2}{x_{n-1}x_n},\\
&&\ovl y_{1}=a(x)\frac{x_0x_1\ovl x_1}
{x_0+x_1\ovl x_1}
=\frac{(x_{n-1}\ovl x_{n-1}
+x_n^2)x_0x_{1}\ovl x_{1}}
{(x_0+x_{1}\ovl x_{1})x_{n-1}x_n^2},\\
&&\ovl y_i
=\frac{(x_{n-1}\ovl x_{n-1}+x_n^2)
x_{i-1}x_{i}\ovl x_{i}}
{(x_{i-1}\ovl x_{i-1}+x_{i}\ovl x_{i})x_{n-1}x_n^2}\q
(1\leq i\leq n-2).
\end{eqnarray*}
It defines a rational mapping
$\osigma\cl\cV(A_{2n}^{(2)\,\dagger})
\longrightarrow 
\cV_2(A_{2n}^{(2)\,\dagger})$
($v_1(x)\mapsto v_2(y)$).\\
The inverse 
$\osigma^{-1}\cl\cV_2(A_{2n}^{(2)\,\dagger})
\longrightarrow \cV(A_{2n}^{(2)\,\dagger})$
$(v_2(y)\mapsto v_1(x))$ is given by 
\begin{eqnarray*}
&&a(y)\seteq \frac{y_0y_1}{y_0+y_1\ovl y_1}(=a(x)),\\
&&x_0=a(y)^{-1}\frac{y_0+y_1\ovl y_1}{y_1},\\
&&x_i=a(y)^{-1}\frac{y_i\ovl y_i
+y_{i+1}\ovl y_{i+1}}{y_{i+1}}\q(1\leq i\leq n-2),
\\
&& x_{n-1}=a(y)^{-1}
\frac{y_{n-1}\ovl y_{n-1}+y_n^2}{y_n^2},\\
&&x_n=a(y)^{-1}y_n,\\
&&\ovl x_i=a(y)^{-1}\frac{y_i\ovl y_iy_{i+1}}
{y_i\ovl y_i+y_{i+1}\ovl y_{i+1}}\q
(1\leq i\leq n-2),\\
&&\ovl x_{n-1}=a(y)^{-1}\frac{y_{n-1}\ovl y_{n-1}
y_n^2}{y_{n-1}\ovl y_{n-1}+y_n^2},
\end{eqnarray*}
which means that the morphism $\osigma\cl
\cV(A_{2n}^{(2)\,\dagger})\longrightarrow 
\cV_2(A_{2n}^{(2)\,\dagger})$ is birational.
Thus, we obtain Theorem \ref{birat}(ii). 

The actions
of $e_i$ $(i=0,1,\cd,n-1)$ on $v_2(y)$
are induced from the ones on 
$Y_{\bf i_2}(y)$
(${\bf i_2}=(n,\cd,1,0,1,\cd,n-1)$)
since $e_i\fsquare(0.5cm,\ovl n)=0$ for $i=0,1,\cd,n-1$.
We also get $\gamma_i(v_2(y))$ and 
$\vep_i(v_2(y))$ from the ones for 
$Y_{\bf i_2}(y)$ where $v_2(y)=\osigma(v_1(x))$:
\begin{eqnarray*}
&&e_0^c\cl y_0\mapsto cy_0,\\
&&e_1^c\cl y_1\mapsto y_1\frac{cy_1\ovl y_1+y_0}
{y_1\ovl y_1+y_0}, \q
\ovl y_1\mapsto \ovl y_1
\frac{c(y_1\ovl y_1+y_0)}{cy_1\ovl y_1+y_0},
\\
&&e_{i}^c\cl y_i\mapsto y_i
\frac{cy_i\ovl y_i+y_{i-1}\ovl y_{i-1}}
{y_i\ovl y_i+y_{i-1}\ovl y_{i-1}},\q
\ovl y_i\mapsto \ovl y_i
\frac{c(y_i\ovl y_i+y_{i-1}\ovl y_{i-1})}
{cy_i\ovl y_i+y_{i-1}\ovl y_{i-1}},\q
\q(i=2,\cd,n-1),\\
&&\gamma_0(v_2(y))
=\frac{y_0^2}{(y_1\ovl y_1)^2},\q
\gamma_1(v_2(y))=
\frac{(y_1\ovl y_1)^2}{y_0y_2\ovl y_2},\\
&&\gamma_i(v_2(y))=
\frac{(y_i\ovl y_i)^2}
{y_{i-1}\ovl y_{i-1}y_{i+1}\ovl y_{i+1}}
\q(i=2,\cd, n-2),\q
\gamma_{n-1}(v_2(y))
=\frac{(y_{n-1}\ovl y_{n-1})^2}
{y_{n-2}\ovl y_{n-2}y_{n}^2},\\
&&\vep_0(v_2(y))=\frac{y_1^2}{y_0},\,
\vep_1(v_2(y))=\frac{y_2}{y_1}
\left(1+\frac{y_0}{y_1\ovl y_1}\right),\,
\vep_{n-1}(v_2(y))=
\frac{y_n^2}{y_{n-1}}\left(
1+\frac{y_{n-2}\ovl y_{n-2}}
{y_{n-1}\ovl y_{n-1}}\right),\\
&&\vep_i(v_2(y))=\frac{y_{i+1}}{y_i}
\left(1+\frac{y_{i-1}\ovl y_{i-1}}
{y_{i}\ovl y_i}\right)\q
(i=2\cd,n-1).
\end{eqnarray*}

The explicit forms of $\vep_i(v_1(x))$and 
$\gamma_i(v_1(x))$
$(1\leq i\leq n)$ are also induced from 
the ones
for $Y_{\bf i_1}(x)\seteq 
Y_0(x_0)\cd Y_1(\ovl x_1)$ and, we define
$\vep_0(v_1(x))\seteq \vep_0(v_2(y))$ and 
$\gamma_0(v_1(x))\seteq \gamma_0(v_2(y))$
$(v_2(y)\seteq \osigma(v_1(x)))$:
\begin{eqnarray*}
&&\vep_0(v_1(x))=
\frac{1}{x_0}\left(1+\frac{x_1\ovl x_1}{x_0}
\right)^2,\q
\vep_i(v_1(x))=\frac{x_{i-1}}{x_i}
\left(1+\frac{x_{i+1}\ovl x_{i+1}}
{x_{i}\ovl x_i}\right)\q
(1\leq i\leq n-2),\\
&&\vep_{n-1}(v_1(x))=\frac{x_{n-2}}{x_{n-1}}
\left(1+\frac{x_n^2}{x_{n-1}\ovl x_{n-1}}
\right),\q
\vep_n(v_1(x))=\frac{x_{n-1}}{x_n}.
\end{eqnarray*}
\begin{eqnarray*}
&&\hspace{-23pt}\gamma_0(v_1(x))=\frac{x_0^2}{(x_1\ovl x_1)^2},\,\,
\gamma_1(v_1(x))=\frac{(x_1\ovl x_1)^2}
{x_0x_2\ovl x_2},\,\,
\gamma_i(v_1(x))=
\frac{(x_i\ovl x_i)^2}{x_{i-1}\ovl x_{i-1}x_{i+1}\ovl x_{i+1}}\,\,
(2\leq i\leq n-2),\\
&&\hspace{-23pt}\gamma_{n-1}(v_1(x))=
\frac{(x_{n-1}\ovl x_{n-1})^2}
{x_{n-2}\ovl x_{n-2}x_n^2},\q
\gamma_n(v_1(x))=\frac{x_n^2}{x_{n-1}\ovl x_{n-1}}
\end{eqnarray*}
For $i=0$, we define 
$e^c_0(v_1(x))=
\osigma^{-1}\circ e^c_0\circ\osigma(v_1(x))
=\osigma^{-1}\circ e^c_0(v_2(y))$.
Then we get
\begin{eqnarray*}
e_0^c\cl&x_0\mapsto&x_0
\frac{(cx_0+x_1\bar x_1)^2}
{c(x_0+x_1\bar x_1)^2},\qq
x_i\mapsto{x_i}\frac{cx_0+x_1\bar x_1}
{c(x_0+x_1\bar x_1)}\,\,(1\leq i\leq n),\\
&\bar x_i\mapsto&{\bar x_i}
\frac{cx_0+x_1\bar x_1}
{c(x_0+x_1\bar x_1)}\,\,(1\leq i\leq n-1),\\
e_i^c\cl&x_i\mapsto&
x_i\frac{cx_i\ovl x_i+x_{i+1}\bar x_{i+1}}
{x_i\ovl x_i+x_{i+1}\ovl x_{i+1}},\q
\ovl x_i\mapsto
\ovl x_i\frac{c(x_i\ovl x_i+x_{i+1}\bar x_{i+1})}
{cx_i\ovl x_i+x_{i+1}\ovl x_{i+1}},\ (1\leq i<n-1),\\
e_{n-1}^c\cl&x_{n-1}\mapsto&
x_{n-1}\frac{cx_{n-1}\ovl x_{n-1}+x_{n}^2}
{x_{n-1}\ovl x_{n-1}+x_{n}^2},\qq
\ovl x_{n-1}\mapsto
\ovl x_{n-1}
\frac{c(x_{n-1}\ovl x_{n-1}+x_{n}^2)}
{cx_{n-1}\ovl x_{n-1}+x_{n}^2}\q,\\
e^c_n\cl&x_n\mapsto& cx_n.
\end{eqnarray*}
In order to prove Theorem \ref{aff-geo}, it 
suffices to show the following:
\begin{eqnarray}
&&
e_i^c=\osigma^{-1}\circ e_{i}^c\circ \osigma,\,\,
\gamma_i=\gamma_{i}\circ \osigma,\,\,
\vep_i=\vep_{i}\circ \osigma
\,\,(i\ne0,n),
\label{0i}\\
&&e_0^{c_1}e_n^{c_2}=e_n^{c_2}e_0^{c_1},
\label{e000}\\
&&\gamma_n(e_0^c(v_1(x)))=\gamma_n(v_1(x)),\q
\gamma_0(e_n^c(v_1(x)))=\gamma_0(v_1(x)),
\label{0nn0}\\
&&\vep_0(e_0^c(v_1(x)))=c^{-1}\vep_0(v_1(x)),
\label{vep00}
\end{eqnarray}
which are immediate from the above formulae.
Let us show \eqref{0i}.
Set $v_2(y)\seteq \osigma(v_1(x))$ and $v_1(x')
\seteq \osigma^{-1}(e_i^c(v_2(y)))$
for $i=2,\cd,n-2$. Then
$x_j'=x_j$ and $\ovl x_j'=\ovl x_j$ for $j\ne i-1,i$, and we have 
\begin{eqnarray*}
&&a(v_2(y))=a(v_1(x)),\\
&&x_i'=\frac{1}{a(v_2(y))}\left(
\ovl y_{i+1}+\frac{cy_i\ovl y_i}{y_{i+1}}\right)
=x_i\frac{cx_i\ovl x_i+x_{i+1}\ovl x_{i+1}}
{x_i\ovl x_i+x_{i+1}\ovl x_{i+1}},\\
&&\ovl x_i'=\frac{1}{a(v_2(y))}
\frac{cy_i\ovl y_iy_{i+1}}
{cy_i\ovl y_i+y_{i+1}\ovl y_{i+1}}
=\ovl x_i\frac{c(x_i\ovl x_i
+x_{i+1}\ovl x_{i+1})}
{cx_i\ovl x_i+x_{i+1}\ovl x_{i+1}},\\
&&x_{i-1}'=\frac{1}{a(v_2(y))}\left(
\ovl y_{i}+\frac{cy_{i-1}\ovl y_{i-1}}
{y_{i}}\right)
=x_{i-1},\\
&&\ovl x_{i-1}'=\frac{1}{a(v_2(y))}
\frac{cy_{i-1}\ovl y_{i-1}y_{i}}
{cy_i\ovl y_i+y_{i-1}\ovl y_{i-1}}
=\ovl x_{i-1},
\end{eqnarray*}
where the formula
$y_i\ovl y_i=a(v_1(x))x_i\ovl x_i$
is useful to obtain these results.
Therefore we have $e_i^c=\osigma^{-1}\circ 
e_i^c\circ\osigma$ for $i=2,\cd,n-2$. 
Others are obtained similarly.

\section{Ultra-discretization of affine geometric crystals}\label{sec6}

In this section, we shall see that
the limit $B_\ify$ of perfect crystals 
as in \S \ref{limit} are isomorphic to the
ultra-discretization of the geometric crystal
obtained in the previous section. 

Let $\cV(\ge)$ be the affine geometric crystal
for $\ge$ as in the previous section.
Let $m={\rm dim}\cV(\ge)$. 
Then due to their explicit forms in the previous section,
we have the following simple positive structure
on $\cV(\ge)$:
\begin{eqnarray}
\Theta\cl(\bbC^\times)^m&\longrightarrow& \cV(\ge),
\\
x&\mapsto& v(x).\nn
\end{eqnarray}
Therefore, we obtain an affine crystal
${\cal UD}(\cV(\ge))
=(\cB(\ge),\{\eit\}_{i\in I},\{\wt_i\}_{i\in I},
\{\vep_i\}_{i \in I})$.

\begin{thm}
\label{ud}
For 
$\ge=A_n^{(1)},B_n^{(1)},
C_n^{(1)},D_n^{(1)},A_{2n-1}^{(2)},
D_{n+1}^{(2)}$ and $A_{2n}^{(2)\,\dagger}$,
let us denote  its Langlands dual by $\ge^{L}$.
Then we have the following 
isomorphism of crystals:
\begin{equation}
{\cal UD}(\cV(\ge))\cong B_{\ify}(\ge^{L}).
\end{equation}
\end{thm}
Here
$L\cl\TY(A,1,n)\leftrightarrow \TY(A,1,n)$, 
$L\cl\TY(B,1,n)\leftrightarrow \TY(A,2,2n-1)$,
$L\cl\TY(C,1,n)\leftrightarrow \TY(D,2,n+1)$, 
$L\cl\TY(D,1,n)\leftrightarrow \TY(D,1,n)$ and 
$L\cl{\TY(A,2,2n)}\leftrightarrow A_{2n}^{(2)\,\dagger}$.

In the following subsections, we shall show 
Theorem \ref{ud} in each case.

The following formula is useful for this purpose:
\begin{eqnarray}
&&
{\cal UD}\bigl(\dfrac{cx+y}{x+y}\bigr)
={\rm max}(c+x,y)-{\rm max}(x,y),
\end{eqnarray}
Note that, if $c=1$,
$${\rm max}(c+x,y)-{\rm max}(x,y)
=\begin{cases}1&x\geq y,\\0&x<y.\end{cases}$$

\subsection{$\TY(A,1,n)$ $(n\geq2)$}

The crystal structure of 
${\cal UD}(\cV(\TY(A,1,n)))$ is given as follows:
$\cB(\TY(A,1,n))=\bbZ^n$ and for 
$x=(x_1,\cd,x_n)\in\bbZ^n$,
we have
\begin{eqnarray*}
&&\til e_0(x)=(x_1-1,x_2-1,\cd,x_n-1),\\
&&\til e_i(x)=(\cd,x_i+1,\cd),\q
(1\leq i\leq n),\\
&&\wt_0(x)=-x_1-x_n,\q \wt_1(x)=2x_1-x_2,\\
&&\wt_i(x)=-x_{i-1}+2x_i-x_{i+1}
\q(i=2,\cd, n-1),\q \wt_n(x)=2x_n-x_{n-1},\\
&&\vep_0(x)=x_1,\q\vep_i(x)=x_{i+1}-x_i\q
(i=1,\cd,n-1),\q \vep_n(x)=-x_n.
\end{eqnarray*}
Comparing with the result in \ref{a-inf} 
and the above formulae,
we easily see that the following 
map gives an isomorphism of 
$\TY(A,1,n)$-crystals:
\begin{eqnarray*}
\mu\cl\cB(\TY(A,1,n))=
\bbZ^n&\longrightarrow & B_\ify(\TY(A,1,n)),\\
(x_1\cd,x_n)&\mapsto&
(b_1,\cd,b_n)=(x_1,{x_2}-{x_1},\cd,
{x_n}-{x_{n-1}}).\nn
\end{eqnarray*}
Thus we have proved Theorem \ref{ud} for $\TY(A,1,n)$.
\subsection{$\TY(B,1,n)$ ($n\geq 3$)}

Let us see the crystal structure of 
${\cal UD}(\cV(\TY(B,1,n)))
=(\cB(\TY(B,1,n)),\{\eit\}, 
\{\wt_i\},\{\vep _i\})$. 
Due to the formula in \ref{b1n}, 
we have ${\cB}(\TY(B,1,n))=\bbZ^{2n-1}$ 
and for 
$x=(x_1,\cd,x_n,\ovl x_{n-1},\cd,\ovl x_1)\in
{\cB}(\TY(B,1,n))$, we have
\begin{eqnarray*}
&&\hspace{-15pt}\til e_0(x)=
\begin{cases}
(x_1,x_2-1,\cd,x_n-1,\cd,\ovl x_2-1,\ovl x_1-1)
&\text{if }
x_1+\ovl x_1\geq x_2+\ovl x_2,\\
(x_1-1,x_2-1,\cd,x_n-1,\cd,\ovl x_2-1,\ovl x_1)
&\text{if }
x_1+\ovl x_1<x_2+\ovl x_2,
\end{cases}\\
&&\hspace{-15pt}\til e_i(x)=
\begin{cases}
(x_1,\cd,x_i+1,\cd,\ovl x_i,\cd,\ovl x_1)
&\text{if }
x_i+\ovl x_i\geq x_{i+1}+\ovl x_{i+1},\\
(x_1,\cd,x_i,\cd,\ovl x_i+1,\cd, \ovl x_1)
&\text{if }
x_i+\ovl x_i<x_{i+1}+\ovl x_{i+1},
\end{cases}(i=1,\cd,n-2),\\
&&\hspace{-15pt}\til e_{n-1}(x)=
\begin{cases}
(x_1,\cd,x_{n-1}+1,x_n,\ovl x_{n-1},
\cd,\ovl x_1)
&\text{if }
x_{n-1}+\ovl x_{n-1}\geq 2x_n,\\
(x_1,\cd,x_{n-1},x_n,\ovl x_{n-1}+1,
\cd,\ovl x_1)
&\text{if }
x_{n-1}+\ovl x_{n-1}< 2x_n,\\
\end{cases}\\
&&\hspace{-15pt}
\til e_n(x)=(x_1\cd, x_{n-1},x_n+1,\ovl x_{n-1}
,\cd, \ovl x_1),\\
&&\hspace{-15pt}
\wt_0(x)=-(x_2+\ovl x_2), \q
\wt_1(x)=2(x_1+\ovl x_1)-(x_2+\ovl x_2),\\
&&\hspace{-15pt}
\wt_i(x)=-(x_{i-1}+\ovl x_{i-1})+2(x_i+\ovl x_i)
-(x_{i+1}+\ovl x_{i+1})\q(i=2,\cd,n-2),\\
&&\hspace{-15pt}
\wt_{n-1}(x)=-(x_{n-2}+\ovl x_{n-2})
+2(x_{n-1}+\ovl x_{n-1})-2x_n,\q
\wt_n(x)=2x_n-(x_{n-1}+\ovl x_{n-1}),\\
&&\hspace{-15pt}\vep_0(x)={\rm max}(x_1+\ovl x_1,
x_2+\ovl x_2)-x_1
=\ovl x_1+(x_2+\ovl x_2-x_1-\ovl x_1)_+,\\
&&\hspace{-15pt}\vep_1(x)=(x_2+\ovl x_2-x_1-\ovl x_1)_+-x_1,\\
&&\hspace{-15pt}\vep_i(x)=x_{i-1}-x_i
+(x_{i+1}+\ovl x_{i+1}-x_i-\ovl x_i)_+\,\,
(i=2,\cd,n-2),\\
&&\hspace{-15pt}\vep_{n-1}(x)=x_{n-2}-x_{n-1}
+(2x_n-x_{n-1}-\ovl x_{n-1})_+,
\q\vep_n(x)=x_{n-1}-x_n.
\end{eqnarray*}
Comparing the result in \ref{aon-ify} and the above 
formulae, it is easy to see  that the following 
map gives an isomorphism of 
$\TY(A,2,2n-1)$-crystals: 
\begin{eqnarray*}
\mu\cl\cB(\TY(B,1,n))&\longrightarrow&
B_\ify(\TY(A,2,2n-1))\\
(x_1,\cd,x_n,\cd,\ovl x_1)&\mapsto&
(b_1,\cd,b_n,\ovl b_n,\cd,\ovl b_1),
\end{eqnarray*}
where 
\begin{eqnarray*}
&&b_1=\ovl x_1, \q b_i=\ovl x_i-\ovl x_{i-1}\,\,
(i=2,\cd,n-1), 
b_n=x_n-\ovl x_{n-1},\\
&&\ovl b_i=x_{i-1}-x_i\,\,
(i=2,\cd,n),\
\ovl b_1=-x_1.
\end{eqnarray*}
We have proved 
Theorem \ref{ud} for $\TY(B,1,n)$.

\subsection{$\TY(C,1,n)$ ($n\geq 2$)}

Let us see the crystal structure of 
${\cal UD}(\cV(\TY(C,1,n)))
=(\cB(\TY(C,1,n)),\{\eit\}, 
\{\wt_i\},\{\vep _i\})$. 
We have ${\cB}(\TY(C,1,n))=\bbZ^{2n}$ and for 
$x=(x_0,x_1,\cd,x_n,\ovl x_{n-1},
\cd,\ovl x_1)\in
{\cB}(\TY(C,1,n))$, we have
\begin{eqnarray*}
&&\hspace{-15pt}\til e_0(x)=
\begin{cases}
(x_0+1,x_1,\cd,x_n,\ovl x_{n-1},
\cd,\ovl x_2,\ovl x_1)
&\text{if }
x_0\geq x_1+\ovl x_1,\\
(x_0-1,\cd,x_{n-1}-1,x_n-2,\ovl x_{n-1}-1,
\cd,\ovl x_2-1,\ovl x_1-1)
&\text{if }
x_0<x_1+\ovl x_1,
\end{cases}\\
&&\hspace{-15pt}\til e_i(x)=
\begin{cases}
(x_0,\cd,x_i+1,\cd,\ovl x_i,\cd,\ovl x_1)
&\text{if }
x_i+\ovl x_i\geq x_{i+1}+\ovl x_{i+1},\\
(x_0,\cd,x_i,\cd,\ovl x_i+1,\cd, \ovl x_1)
&\text{if }
x_i+\ovl x_i<x_{i+1}+\ovl x_{i+1},
\end{cases}(i=1,\cd,n-2),\\
&&\hspace{-15pt}\til e_{n-1}(x)=
\begin{cases}
(x_0,\cd,x_{n-1}+1,x_n,\ovl x_{n-1},
\cd,\ovl x_1)
&\text{if }
x_{n-1}+\ovl x_{n-1}\geq x_n,\\
(x_0,\cd,x_{n-1},x_n,\ovl x_{n-1}+1,
\cd,\ovl x_1)
&\text{if }
x_{n-1}+\ovl x_{n-1}< x_n,\\
\end{cases}\\
&&\hspace{-15pt}
\til e_n(x)=(x_0,\cd,x_{n-1},x_n+1,\ovl x_{n-1}
,\cd, \ovl x_1),\\
&&\hspace{-15pt}
\wt_0(x)=2x_0-2(x_1+\ovl x_1), \q
\wt_1(x)=2(x_1+\ovl x_1)-(x_0+x_2+\ovl x_2),\\
&&\hspace{-15pt}
\wt_i(x)=-(x_{i-1}+\ovl x_{i-1})+2(x_i+\ovl x_i)
-(x_{i+1}+\ovl x_{i+1})\q(i=2,\cd,n-2),\\
&&\hspace{-15pt}
\wt_{n-1}(x)=-(x_{n-2}+\ovl x_{n-2})
+2(x_{n-1}+\ovl x_{n-1})-x_n,\q
\wt_n(x)=2x_n-2(x_{n-1}+\ovl x_{n-1}),\\
&&\hspace{-15pt}\vep_0(x)
=-x_0+2(x_1+\ovl x_1-x_0)_+,\\
&&\hspace{-15pt}\vep_i(x)=x_{i-1}-x_i
+(x_{i+1}+\ovl x_{i+1}-x_i-\ovl x_i)_+\,\,
(i=1,\cd,n-2),\\
&&\hspace{-15pt}\vep_{n-1}(x)=x_{n-2}-x_{n-1}
+(x_n-x_{n-1}-\ovl x_{n-1})_+,
\q\vep_n(x)=2x_{n-1}-x_n.
\end{eqnarray*}
Comparing the result in \ref{d2n-ify} and the above 
formulae, we see that the following 
map gives an isomorphism of 
$\TY(D,2,n+1)$-crystals: 
\begin{eqnarray*}
\mu\cl\cB(\TY(C,1,n))&\longrightarrow&
B_\ify(\TY(D,2,n+1))\\
(x_0,x_1,\cd,x_n,\cd,\ovl x_1)&\mapsto&
(b_1,\cd,b_n,\ovl b_n,\cd,\ovl b_1),
\end{eqnarray*}
where
\begin{eqnarray*}
&&
b_1=\ovl x_1, \q b_i=\ovl x_i-\ovl x_{i-1}\,\,
(i=2,\cd,n-1), 
b_n=\frac{1}{2}x_n-\ovl x_{n-1},\\
&&\ovl b_n=x_{n-1}-\frac{1}{2}x_n,\q
\ovl b_i=x_{i-1}-x_i\,\,
(i=1,\cd,n-1).
\end{eqnarray*}
We have proved  
Theorem \ref{ud} for $\TY(C,1,n)$.

\subsection{$\TY(D,1,n)$ ($n\geq 4$)}

Let us see the crystal structure of 
${\cal UD}(\cV(\TY(D,1,n)))
=(\cB(\TY(D,1,n)),\{\eit\}, 
\{\wt_i\},\{\vep _i\})$. 
We have ${\cB}(\TY(D,1,n))=\bbZ^{2n-2}$ 
and for 
$x=(x_1,\cd,x_n,\ovl x_{n-1},\cd,\ovl x_1)\in
{\cB}(\TY(D,1,n))$, we have
\begin{eqnarray*}
&&\hspace{-15pt}\til e_0(x)=
\begin{cases}
(x_1,x_2-1,\cd,x_n-1,\cd,\ovl x_2-1,\ovl x_1-1)
&\text{if }
x_1+\ovl x_1\geq x_2+\ovl x_2,\\
(x_1-1,x_2-1,\cd,x_n-1,\cd,\ovl x_2-1,
\ovl x_1)&\text{if }
x_1+\ovl x_1<x_2+\ovl x_2,
\end{cases}\\
&&\hspace{-15pt}\til e_i(x)=
\begin{cases}
(x_1,\cd,x_i+1,\cd,\ovl x_i,\cd,\ovl x_1)
&\text{if }
x_i+\ovl x_i\geq x_{i+1}+\ovl x_{i+1},\\
(x_1,\cd,x_i,\cd,\ovl x_i+1,\cd, \ovl x_1)
&\text{if }
x_i+\ovl x_i<x_{i+1}+\ovl x_{i+1},
\end{cases}(i=1,\cd,n-3),\\
&&\hspace{-15pt}\til e_{n-2}(x)=
\begin{cases}
(x_1,\cd,x_{n-2}+1,\cd,\ovl x_{n-2},
\cd,\ovl x_1)
&\text{if }
x_{n-2}+\ovl x_{n-2}\geq x_{n-1}+x_{n},\\
(x_1,\cd,x_{n-2},\cd,\ovl x_{n-2}+1,\cd, \ovl x_1)
&\text{if }
x_{n-2}+\ovl x_{n-2}<x_{n-1}+x_{n},
\end{cases}\\
&&\hspace{-15pt}\til e_{n-1}(x)=
(x_1,\cd,x_{n-1}+1,x_n,\cd,\ovl x_1),\\
&&\hspace{-15pt}\til e_n(x)
=(x_1\cd, x_{n-1},x_n+1,\cd, \ovl x_1),\\
&&\hspace{-15pt}
\wt_0(x)=-(x_2+\ovl x_2), \q
\wt_1(x)=2(x_1+\ovl x_1)-(x_2+\ovl x_2),\\
&&\hspace{-15pt}
\wt_i(x)=-(x_{i-1}+\ovl x_{i-1})+2(x_i+\ovl x_i)
-(x_{i+1}+\ovl x_{i+1})\q(i=2,\cd,n-3),\\
&&\hspace{-15pt} \wt_{n-2}(x)=-(x_{n-3}+\ovl x_{n-3})+
2(x_{n-2}+\ovl x_{n-2})-(x_{n-1}+x_n),\\
&&\hspace{-15pt}
\wt_{n-1}(x)=-(x_{n-2}+\ovl x_{n-2})
+2(x_{n-1}+\ovl x_{n-1}),\q
\wt_n(x)=2x_n-(x_{n-2}+\ovl x_{n-2}),\\
&&\hspace{-15pt}\vep_0(x)={\rm max}(x_1+\ovl x_1,
x_2+\ovl x_2)-x_1
=\ovl x_1+(x_2+\ovl x_2-x_1-\ovl x_1)_+,\\
&&\hspace{-15pt}\vep_1(x)=(x_2+\ovl x_2-x_1-\ovl x_1)_+-x_1,\\
&&\hspace{-15pt}\vep_i(x)=x_{i-1}-x_i
+(x_{i+1}+\ovl x_{i+1}-x_i-\ovl x_i)_+\,\,
(i=2,\cd,n-3),\\
&&\hspace{-15pt}\vep_{n-2}(x)=
x_{n-3}-x_{n-2}+
+(x_{n-1}+x_n-x_{n-2}-\ovl x_{n-2})_+,\\
&&\hspace{-15pt}\vep_{n-1}(x)=x_{n-2}-x_{n-1},
\q\vep_n(x)=x_{n-2}-x_n.
\end{eqnarray*}
Then the following map gives an isomorphism of 
$\TY(D,1,n)$-crystals: 
\begin{eqnarray*}
\mu\cl\cB(\TY(D,1,n))&\longrightarrow&
B_\ify(\TY(D,1,n))\\
(x_1,\cd,x_{n-1},x_n,\ovl x_{n-2},
\cd,\ovl x_1)&\mapsto&
(b_1,\cd,b_n,\ovl b_{n-1},\cd,\ovl b_1),
\end{eqnarray*}
where
\begin{eqnarray*}
&&
b_1=\ovl x_1, \q b_i=\ovl x_i-\ovl x_{i-1}\,\,
(i=2,\cd,n-2), \q
b_{n-1}=x_{n-1}-\ovl x_{n-2},
\q b_n=x_n-x_{n-1},\\
&&\ovl b_{n-1}=x_{n-2}-x_n,\q
\ovl b_i=x_{i-1}-x_i\,\,
(i=2,\cd,n-2),\
\ovl b_1
=-\sum_{i=1}^nb_i-\sum_{i=2}^{n-1}\ovl b_i
=-x_1.
\end{eqnarray*}
We have proved Theorem \ref{ud} for $\TY(D,1,n)$.
\subsection{$\TY(A,2,2n-1)$ ($n\geq 3$)}

The crystal structure of 
${\cal UD}(\cV(\TY(A,2,2n-1)))
=(\cB(\TY(A,2,2n-1)),\{\eit\}, 
\{\wt_i\},\{\vep _i\})$ is given by:\\
${\cB}(\TY(A,2,2n-1))=\bbZ^{2n-1}$
and for 
$x=(x_1,\cd,x_n,\ovl x_{n-1},\cd,\ovl x_1)\in
{\cB}(\TY(A,2,2n-1))$, we have
\begin{eqnarray*}
&&\hspace{-15pt}\til e_0(x)=
\begin{cases}
(x_1,x_2-1,\cd,x_n-2,\cd,\ovl x_2-1,\ovl x_1-1)
&\text{if }
x_1+\ovl x_1\geq x_2+\ovl x_2,\\
(x_1-1,x_2-1,\cd,x_n-2,\cd,\ovl x_2-1,\ovl x_1)
&\text{if }
x_1+\ovl x_1<x_2+\ovl x_2,
\end{cases}\\
&&\hspace{-15pt}\til e_i(x)=
\begin{cases}
(x_1,\cd,x_i+1,\cd,\ovl x_i,\cd,\ovl x_1)
&\text{if }
x_i+\ovl x_i\geq x_{i+1}+\ovl x_{i+1},\\
(x_1,\cd,x_i,\cd,\ovl x_i+1,\cd, \ovl x_1)
&\text{if }
x_i+\ovl x_i<x_{i+1}+\ovl x_{i+1},
\end{cases}(i=1,\cd,n-2),\\
&&\hspace{-15pt}\til e_{n-1}(x)=
\begin{cases}
(x_1,\cd,x_{n-1}+1,x_n,\ovl x_{n-1},
\cd,\ovl x_1)
&\text{if }
x_{n-1}+\ovl x_{n-1}\geq x_n,\\
(x_1,\cd,x_{n-1},x_n,\ovl x_{n-1}+1,
\cd,\ovl x_1)
&\text{if }
x_{n-1}+\ovl x_{n-1}< x_n,\\
\end{cases}\\
&&\hspace{-15pt}
\til e_n(x)=(x_1\cd, x_{n-1},x_n+1,\ovl x_{n-1}
,\cd, \ovl x_1),\\
&&\hspace{-15pt}
\wt_0(x)=-(x_2+\ovl x_2), \q
\wt_1(x)=2(x_1+\ovl x_1)-(x_2+\ovl x_2),\\
&&\hspace{-15pt}
\wt_i(x)=-(x_{i-1}+\ovl x_{i-1})+2(x_i+\ovl x_i)
-(x_{i+1}+\ovl x_{i+1})\q(i=2,\cd,n-2),\\
&&\hspace{-15pt}
\wt_{n-1}(x)=-(x_{n-2}+\ovl x_{n-2})
+2(x_{n-1}+\ovl x_{n-1})-x_n,\q
\wt_n(x)=2x_n-2(x_{n-1}+\ovl x_{n-1}),\\
&&\hspace{-15pt}\vep_0(x)={\rm max}(x_1+\ovl x_1,
x_2+\ovl x_2)-x_1
=\ovl x_1+(x_2+\ovl x_2-x_1-\ovl x_1)_+,\\
&&\hspace{-15pt}\vep_1(x)=(x_2+\ovl x_2-x_1-\ovl x_1)_+-x_1,\\
&&\hspace{-15pt}\vep_i(x)=x_{i-1}-x_i
+(x_{i+1}+\ovl x_{i+1}-x_i-\ovl x_i)_+\,\,
(i=2,\cd,n-2),\\
&&\hspace{-15pt}\vep_{n-1}(x)=x_{n-2}-x_{n-1}
+(x_n-x_{n-1}-\ovl x_{n-1})_+,
\q\vep_n(x)=2x_{n-1}-x_n.
\end{eqnarray*}
Then the following map gives an isomorphism of 
$\TY(B,1,n)$-crystals: 
\begin{eqnarray*}
\mu\cl\cB(\TY(A,2,2n-1))&\longrightarrow&
B_\ify(\TY(B,1,n))\\
(x_1,\cd,x_n,\cd,\ovl x_1)&\mapsto&
(b_1,\cd,b_n,\ovl b_n,\cd,\ovl b_1),
\end{eqnarray*}
where
\begin{eqnarray*}
&&
b_1=\ovl x_1, \q b_i=\ovl x_i-\ovl x_{i-1}\,\,
(i=2,\cd,n-1), 
b_n=\frac{1}{2}x_n-\ovl x_{n-1},\\
&&\ovl b_n=x_{n-1}-\frac{1}{2}x_n,\q
\ovl b_i=x_{i-1}-x_i\,\,
(i=2,\cd,n-1),\
\ovl b_1=-x_1.
\end{eqnarray*}
We have proved
Theorem \ref{ud} for $\TY(A,2,2n-1)$.

\subsection{$\TY(D,2,n+1)$ ($n\geq 2$)}

Let us see the crystal structure of 
${\cal UD}(\cV(\TY(D,2,n+1)))
=(\cB(\TY(D,2,n+1)),\{\eit\}, 
\{\wt_i\},\{\vep _i\})$. 
We have ${\cB}(\TY(D,2,n+1))
=\bbZ^{2n}$ and for 
$x=(x_0,x_1,\cd,x_n,\ovl x_{n-1},
\cd,\ovl x_1)\in
{\cB}(\TY(D,2,n+1))$, we have
\begin{eqnarray*}
&&\hspace{-35pt}\til e_0(x)=
\begin{cases}
(x_0+1,x_1,\cd,x_n,\ovl x_{n-1},
\cd,\ovl x_2,\ovl x_1)
&\text{if }
2x_0\geq x_1+\ovl x_1,\\
(x_0,x_1-1\cd,x_n-1,\ovl x_{n-1}-1,
\cd,\ovl x_1-1)&\text{if }
2x_0+1=x_{1}+\ovl x_1,\\
(x_0-1,\cd,x_{n-1}-2,x_n-2,\ovl x_{n-1}-2,
\cd,\ovl x_2-2,\ovl x_1-2)
&\text{if }
2x_0+1<x_1+\ovl x_1,
\end{cases}
\end{eqnarray*}
where we use the formula 
\[
 {\rm max}(2+x,y)-{\rm max}(x,y)
=\begin{cases}2& x\geq y,\\
1&x+1=y,\\
0&x+1<y,
\end{cases}
(x,y\in \bbZ).
\]
\begin{eqnarray*}
&&\hspace{-15pt}\til e_i(x)=
\begin{cases}
(x_1,\cd,x_i+1,\cd,\ovl x_i,\cd,\ovl x_1)
&\text{if }
x_i+\ovl x_i\geq x_{i+1}+\ovl x_{i+1},\\
(x_1,\cd,x_i,\cd,\ovl x_i+1,\cd, \ovl x_1)
&\text{if }
x_i+\ovl x_i<x_{i+1}+\ovl x_{i+1},
\end{cases}(i=1,\cd,n-2),\\
&&\hspace{-15pt}\til e_{n-1}(x)=
\begin{cases}
(x_1,\cd,x_{n-1}+1,x_n,\ovl x_{n-1},
\cd,\ovl x_1)
&\text{if }
x_{n-1}+\ovl x_{n-1}\geq 2x_n,\\
(x_1,\cd,x_{n-1},x_n,\ovl x_{n-1}+1,
\cd,\ovl x_1)
&\text{if }
x_{n-1}+\ovl x_{n-1}< 2x_n,\\
\end{cases}\\
&&\hspace{-15pt}
\til e_n(x)=(x_1,\cd,x_{n-1},x_n+1,\ovl x_{n-1}
,\cd, \ovl x_1),\\
&&\hspace{-15pt}
\wt_0(x)=2x_0-(x_1+\ovl x_1), \q
\wt_1(x)=2(x_1+\ovl x_1)-(2x_0+x_2+\ovl x_2),\\
&&\hspace{-15pt}
\wt_i(x)=-(x_{i-1}+\ovl x_{i-1})+2(x_i+\ovl x_i)
-(x_{i+1}+\ovl x_{i+1})\q(i=2,\cd,n-2),\\
&&\hspace{-15pt}
\wt_{n-1}(x)=-(x_{n-2}+\ovl x_{n-2})
+2(x_{n-1}+\ovl x_{n-1})-2x_n,\q
\wt_n(x)=2x_n-(x_{n-1}+\ovl x_{n-1}),\\
&&\hspace{-15pt}\vep_0(x)
=-x_0+(x_1+\ovl x_1-2x_0)_+,\q
\vep_1(x)=2x_0-x_1+(x_2+\ovl x_2-x_1-\ovl x_1)_+\\
&&\hspace{-15pt}\vep_i(x)=x_{i-1}-x_i
+(x_{i+1}+\ovl x_{i+1}-x_i-\ovl x_i)_+\,\,
(i=2,\cd,n-2),\\
&&\hspace{-15pt}\vep_{n-1}(x)=x_{n-2}-x_{n-1}
+(2x_n-x_{n-1}-\ovl x_{n-1})_+,
\q\vep_n(x)=x_{n-1}-x_n.
\end{eqnarray*}
Then the following map gives an isomorphism of 
$\TY(C,1,n)$-crystals: 
\begin{eqnarray*}
\mu\cl\cB(\TY(D,2,n+1))&\longrightarrow&
B_\ify(\TY(C,1,n))\\
(x_0,x_1,\cd,x_n,\cd,\ovl x_1)&\mapsto&
(b_1,\cd,b_n,\ovl b_n,\cd,\ovl b_1),
\end{eqnarray*}
where
\[
b_1=\ovl x_1, \q b_i=\ovl x_i-\ovl x_{i-1}\,\,
(i=2,\cd,n), 
\ovl b_i=x_{i-1}-x_i\,\,
(i=2,\cd,n),\q
\ovl b_1=2x_0-x_1.
\]
We have proved Theorem \ref{ud} for $\TY(D,2,n+1)$.
\subsection{$A_{2n}^{(2)\,\dagger}$ ($n\geq 2$)}

Let $\cV(A_{2n}^{(2)\,\dagger})$ be the 
affine $A_{2n}^{(2)\,\dagger}$-geometric crystal 
as in \ref{aen-1}. 
We shall see the crystal structure of 
${\cal UD}(\cV(A_{2n}^{(2)\,\dagger}))
=(\cB(A_{2n}^{(2)\,\dagger}),\{\eit\}, 
\{\wt_i\},\{\vep _i\})$ is given by:
${\cB}(A_{2n}^{(2)\,\dagger})
=\bbZ^{2n}$ and for 
$x=(x_0,x_1,\cd,x_n,\ovl x_{n-1},
\cd,\ovl x_1)\in
{\cB}(A_{2n}^{(2)\,\dagger})$, we have
\begin{eqnarray*}
&&\hspace{-15pt}\til e_0(x)=
\begin{cases}
(x_0+1,x_1,\cd,x_n,\ovl x_{n-1},
\cd,\ovl x_2,\ovl x_1)
&\text{if }
x_0\geq x_1+\ovl x_1,\\
(x_0-1,x_1-1\cd,x_{n-1}-1,x_n-1,
\ovl x_{n-1}-1,
\cd,\ovl x_1-1)&\text{if }
x_0<x_{1}+\ovl x_1,
\end{cases}\\
&&\hspace{-15pt}\til e_i(x)=
\begin{cases}
(x_1,\cd,x_i+1,\cd,\ovl x_i,\cd,\ovl x_1)
&\text{if }
x_i+\ovl x_i\geq x_{i+1}+\ovl x_{i+1},\\
(x_1,\cd,x_i,\cd,\ovl x_i+1,\cd, \ovl x_1)
&\text{if }
x_i+\ovl x_i<x_{i+1}+\ovl x_{i+1},
\end{cases}(i=1,\cd,n-2),\\
&&\hspace{-15pt}\til e_{n-1}(x)=
\begin{cases}
(x_1,\cd,x_{n-1}+1,x_n,\ovl x_{n-1},
\cd,\ovl x_1)
&\text{if }
x_{n-1}+\ovl x_{n-1}\geq 2x_n,\\
(x_1,\cd,x_{n-1},x_n,\ovl x_{n-1}+1,
\cd,\ovl x_1)
&\text{if }
x_{n-1}+\ovl x_{n-1}< 2x_n,\\
\end{cases}\\
&&\hspace{-15pt}
\til e_n(x)
=(x_1,\cd,x_{n-1},x_n+1,\ovl x_{n-1}
,\cd, \ovl x_1),\\
&&\hspace{-15pt}
\wt_0(x)=2x_0-2(x_1+\ovl x_1), \q
\wt_1(x)=2(x_1+\ovl x_1)-(x_0+x_2+\ovl x_2),\\
&&\hspace{-15pt}
\wt_i(x)=-(x_{i-1}+\ovl x_{i-1})+2(x_i+\ovl x_i)
-(x_{i+1}+\ovl x_{i+1})\q(i=2,\cd,n-2),\\
&&\hspace{-15pt}
\wt_{n-1}(x)=-(x_{n-2}+\ovl x_{n-2})
+2(x_{n-1}+\ovl x_{n-1})-2x_n,\q
\wt_n(x)=2x_n-(x_{n-1}+\ovl x_{n-1}),\\
&&\hspace{-15pt}\vep_0(x)
=-x_0+2(x_1+\ovl x_1-x_0)_+,\\
&&\hspace{-15pt}\vep_i(x)=x_{i-1}-x_i
+(x_{i+1}+\ovl x_{i+1}-x_i-\ovl x_i)_+\,\,
(i=1,\cd,n-2),\\
&&\hspace{-15pt}\vep_{n-1}(x)=x_{n-2}-x_{n-1}
+(2x_n-x_{n-1}-\ovl x_{n-1})_+,
\q\vep_n(x)=x_{n-1}-x_n.
\end{eqnarray*}
Then the following map gives an isomorphism of 
$\TY(A,2,2n)$-crystals: 
\begin{eqnarray*}
\mu\cl \cB(A_{2n}^{(2)\,\dagger})&\longrightarrow&
B_\ify({\TY(A,2,2n)})\\
(x_0,x_1,\cd,x_n,\cd,\ovl x_1)&\mapsto&
(b_1,\cd,b_n,\ovl b_n,\cd,\ovl b_1),
\end{eqnarray*}
where
\[
b_1=\ovl x_1, \q b_i=\ovl x_i-\ovl x_{i-1}\,\,
(i=2,\cd,n-1), \,\,
b_n=x_n-\ovl x_{n-1},\,\,
\ovl b_i=x_{i-1}-x_i\,\,
(i=1,\cd,n).
\]
We have proved 
Theorem \ref{ud} for $A_{2n}^{(2)\,\dagger}$.



\end{document}